\documentclass[reqno,11pt]{amsart}
\usepackage{geometry}
\usepackage{mathtools,amssymb,amsthm,mathrsfs,color,lineno,paralist,graphicx,float}
\usepackage[colorlinks,
linkcolor=blue,
anchorcolor=green,
citecolor=blue, 
]{hyperref}

\usepackage[T1]{fontenc}
\usepackage[utf8]{inputenc}
\usepackage{tikz}

\usepackage{calc}
\definecolor{bleu1}{RGB}{0,57,128}
\def\bleu1{\color{bleu1}}

\usepackage{etoolbox}
\patchcmd{\section}{\normalfont}{\normalfont \bleu1}{}{}
\patchcmd{\subsection}{\normalfont}{\normalfont \bleu1}{}{}
\patchcmd{\subsubsection}{\normalfont}{\normalfont \bleu1}{}{}
\renewcommand{\proofname}{\it \bleu1 Proof}

\newcommand{\SL}{\mathrm{SL}(2,\mathbb{R})}

\def\a{\alpha}

\def\e{\varepsilon}

\let\newpf\proof \let\proof\relax 
\newenvironment{pf}{\newpf[\proofname]}{\qed\endtrivlist}

\renewcommand{\sl}{{\mathrm{sl}}}

\def\u{{\mathbb U}}

\def\a{{\alpha}}

\def\SL{{\mathrm{SL}}}
\def\PSL{{\mathrm{PSL}}}

\def\0{{\mathbf 0}}

\newtheorem{Theorem}{Theorem}[section]
\newtheorem{Lemma}{Lemma}[section]
\newtheorem{Proposition}{Proposition}[section]
\newtheorem{Corollary}{Corollary}[section]
\newtheorem{Remark}{Remark}[section]

\newtheorem{Definition}{Definition}[section]
\newtheorem{Claim}{Claim}

\numberwithin{equation}{section}

\theoremstyle{definition}

\def\ssm{\smallsetminus}

\renewcommand{\setminus}{\ssm}

\renewcommand{\mod}{\operatorname{mod}}

\newcommand{\id}{\operatorname{id}}

\def\tF{\widetilde{F}}

\newcommand{\A}{{\mathbb A}}
\newcommand{\C}{{\mathbb C}}

\newcommand{\N}{{\mathbb N}}
\newcommand{\Q}{{\mathbb Q}}
\newcommand{\R}{{\mathbb R}}
\newcommand{\T}{{\mathbb T}}

\newcommand{\Z}{{\mathbb Z}}

\newcommand{\la}{\langle}
\newcommand{\ra}{\rangle}

\def\B0{{\bold{0}}}

\begin{document}
\pagestyle{plain}
\author{Lingrui Ge}
\address{
Beijing International Center for Mathematical Research, Peking University, Beijing, China
}

 \email{gelingrui@bicmr.pku.edu.cn}

\author {Jiangong You}
\address{
Chern Institute of Mathematics and LPMC, Nankai University, Tianjin 300071, China} \email{jyou@nankai.edu.cn}

\author{Qi Zhou}
\address{
Chern Institute of Mathematics and LPMC, Nankai University, Tianjin 300071, China
}

 \email{qizhou@nankai.edu.cn}

\title[]{Sharp Spectral Gaps, Arithmetic Localization, and Reducibility via Resonance Analysis}

\maketitle

\begin{abstract}
This paper establishes several sharp spectral results for analytic quasiperiodic Schr\"odinger operators. Key contributions include: (1) exact exponential decay rates for spectral gaps of the almost Mathieu operator, addressing a question raised by Goldstein; (2) sharp arithmetic Anderson localization for a class of quasiperiodic operators on higher-dimensional lattices, which in particular resolves and generalizes Jitomirskaya's phase transition conjecture; and (3) stratified growth patterns for extended eigenfunctions revealing universal partial hierarchical structures for subcritical quasiperiodic Schr\"odinger operators.

The proofs are based on novel frameworks—\textit{structured quantitative almost reducibility} and \textit{sharp quantitative duality}—to overcome the longstanding challenge of taming infinitely many rotation-number resonances, which enables us to obtain optimal arithmetic reducibility results for analytic $\SL(2,\mathbb{R})$-cocycles, thereby solving another Jitomirskaya's conjecture. These methods enable a first comprehensive treatment of resonance-driven dynamical asymptotics.
\end{abstract}

\section{Introduction}

This paper addresses several sharp spectral results for one-dimensional analytic quasiperiodic Schrödinger operators (QPS), defined as a self-adjoint operator on \( \ell^2(\mathbb{Z}) \):
\[
(H_{V,\alpha,\theta}u)(n) = u(n+1) + u(n-1) + V(\theta + n\alpha)u(n),
\]
where \( V: \mathbb{T}^d \to \mathbb{R} \) is the potential, \( \alpha \in \mathbb{R}^d  \) is the frequency (which is rationally independent), and \( \theta\in\mathbb{T}^d \) is the phase. QPS has been the subject of sustained interest in condensed matter physics \cite{aos, oa}, numerical analysis \cite{JLZ, JZ}, dynamical systems \cite{avila2015global, AK06}, and spectral theory \cite{avila, avila1, aj, ALYZ, ayz, L1} for decades. When \( V(\theta) = 2 \lambda \cos(2\pi \theta) \), the operator reduces to the well-known almost Mathieu operator (AMO):
\[
(H_{\lambda,\alpha,\theta}u)(n) = u(n+1) + u(n-1) + 2\lambda \cos(2\pi(\theta + n\alpha))u(n).
\]
The AMO is also known as the Aubry-André model or the Hofstadter Hamiltonian in physical literature. Initially introduced by Peierls \cite{Peierls} as a model for an electron in a 2D lattice with a homogeneous magnetic field \cite{harper, R}, as pointed out by Marcel den Nijs in the colloquium celebrating Thouless' Nobel Prize,  ``the  genius of Thouless was to focus on the Hofstadter Hamiltonian as a model for
the quantum Hall effect''\cite{Av17}. 

\subsection{Exact Exponential Decay Rate of the Spectral Gaps}
Our first sharp spectral result is related to Thouless et al.'s theory of the quantum Hall effect. To explain this, let us first introduce the spectrum of the Schrödinger operator \( H_{V,\alpha,\theta} \), denoted by \( \Sigma_{V,\alpha} \). It is compact and independent of \( \theta \) for rationally independent \( \alpha \). The intervals in \( \mathbb{R} \setminus \Sigma_{V,\alpha} \) are known as spectral gaps, and their labeling by integers is governed by the Gap-Labelling Theorem (GLT) \cite{JM}. 

In 1981, at the AMS annual meeting, Mark Kac asked if the ``all gaps are there" held for the almost Mathieu operator \( H_{\lambda,\alpha,\theta} \) \cite{Kac,Barry}, offering ten martinis as a prize for a proof. Later, Barry Simon popularized it as the Ten Martini Problem and the Dry Ten Martini Problem. The Ten Martini Problem asks whether the spectrum is a Cantor set, while the Dry Ten Martini Problem asks more precisely whether all gaps allowed by the Gap Labelling Theorem are open. 

Since then, these two problems have attracted significant attention, and important progress has been made in \cite{aj, aj1, AK06, bs, cey, hs, last, PuigC, Sinai}. The Ten Martini Problem was completely solved by Avila and Jitomirskaya \cite{aj}. Compared to the Ten Martini Problem, the Dry Ten Martini Problem is more challenging. It was solved by Avila-You-Zhou \cite{ayz2} for non-critical AMO with all irrational frequencies, and one can consult \cite{cey, aj1, PuigC} for partial advances.

To clarify GLT and the precise meaning of ``all gaps are there", let's first recall that the integrated density of states $N_{V,\alpha}(E)$ of the operator $H_{V,\alpha,\theta}$, which is defined as $$ N_{V,\alpha}(E) = \int_{\R/\Z} \mu_{V,\alpha,\theta}(-\infty, E] d\theta, $$ where $\mu_{V,\alpha,\theta}$ is the universal spectral measure of $H = H_{V,\alpha,\theta}$. According to GLT, every spectral gap $ G $ has a unique integer $ k $ such that $ N_{V,\alpha}|_G = k\alpha \mod \Z $. We denote these gaps by $$ G_k(V) = \begin{cases} (E_k^-, E_k^+) & k \neq 0 \\ (-\infty, \inf \Sigma_{V,\alpha}) \cup (\sup\Sigma_{V,\alpha}, \infty) & k = 0 \end{cases}, $$ and $ |G_k(V)| $ denotes the length of $ G_k(V) $.  For AMO, i.e.,  $V(\theta)=2 \lambda \cos 2\pi \theta$,  we  abbreviate it as $|G_k(\lambda)|$. 

The Dry Ten Martini Problem also has relevance in physics. After Von Klitzing's discovery of the quantum Hall effect \cite{KDP}, Thouless and his coauthors \cite{ntw} provided a theoretical explanation by showing that the Hall conductance at the plateaus is related to a topological invariant known as the Chern number, which leads to its quantization. This is part of the theory for which Thouless was awarded the 2016 Nobel Prize in Physics. 
In the language of physics, the labels \( k \) are identified with Chern numbers, and the statement ``all spectral gaps are open" means that all topological phases will appear for the quantum Hall system governed by the AMO. As explained by Thouless et al. \cite{ntw}, their arguments assume that all gaps \( G_k(\cdot) \) of the AMO are open when \( \lambda \) varies.

In this context, Thouless et al.'s theory of the integer quantum Hall effect \cite{ntw} and the ``Dry Ten Martini Problem'' primarily concern whether \( |G_k(\lambda)| > 0 \) for any \( k \).
A more challenging task is to obtain quantitative estimates on the size of these gaps. Specifically, are there bounds on the upper and lower limits of each spectral gap \( G_k(\lambda) \)? 
These questions are not only  interesting in spectral theory but also crucial in the study of the homogeneity of the spectrum \cite{binder2018almost, LYZZ} and Deift's conjecture \cite{binder2018almost, Deift1, Deift2, LYZZ}.  Indeed, after the third author announced the result of \cite{ayz2} at the conference  ``Almost Periodic and Other Ergodic Problems'' in 2015, Goldstein asked whether any quantitative lower bound on the size of \( G_k(\lambda) \) could be derived in the non-critical case \( |\lambda| \neq 1 \). Our paper answers this question by providing the exact exponential decay rate.

Throughout the paper, we assume that \( \alpha \in \mathbb{R}^d \) is Diophantine, i.e., there exist \( \kappa > 0 \) and \( \tau > d + 1 \) such that
\[
\| \langle n, \alpha \rangle \|_{\mathbb{R}/\mathbb{Z}} \geq \frac{\kappa}{|n|^{\tau}}, \quad \forall n \in \mathbb{Z}^d \setminus \{0\}.
\]
Denote by \( \text{DC}_d(\kappa, \tau) \) the set of \( (\kappa, \tau) \)-Diophantine vectors, and \( \text{DC}_d := \bigcup_{\kappa > 0, \, \tau > d-1} \text{DC}_d(\kappa, \tau) \). In particular, when \( d = 1 \), we simplify the notations to \( \text{DC}(\kappa, \tau) \) and \( \text{DC} \).
Then, our main result is as follows:

\begin{Theorem}\label{optimal gap}
Assume \( \alpha \in \text{DC} \) and \( |\lambda| \neq 1 \), then we have
\[
\lim_{k \to \infty} \frac{\ln |G_k(\lambda)|}{|k|} = -|\ln |\lambda||.
\]
\end{Theorem}

%
%
\begin{Remark}
The upper bound was proven in \cite{LYZZ}, and we include it here for completeness.
\end{Remark}
Let us review some recent works in connection with the question of gap estimates for quasiperiodic Schr\"odinger operators. The study of upper bounds for quasiperiodic Schr\"odinger operators began with Moser and P\"oschel \cite{mp}, who showed that for small analytic $ V: \T^d \rightarrow \R $, the \( k \)-th gap  of  the continuous quasiperiodic Schr\"odinger operator on $L^2(\R)$:
\begin{equation*}
(\mathcal{L}_{V,\alpha}y)(t)=-y''(t)+V(\alpha t)y(t)
\end{equation*} 
 decays exponentially with respect to \( |k| \), provided that 
$k$ satisfies certain arithmetic conditions. Later Amor \cite{amor} extended this result to show sub-exponential decay for any 
$k\in\Z^d\backslash\{0\}$, while Damanik and Goldstein \cite{dg} confirmed an exponential decay upper bound.

The results in \cite{amor,mp} also apply to the discrete quasiperiodic Schr\"odinger operator $H_{V,\alpha,\theta}$, as their proof, based on KAM theory, is applicable to both discrete and continuous settings. In contrast, the proof in \cite{dg} relies on localization techniques and is not directly transferable to the discrete case. Recently, Leguil-You-Zhao-Zhou \cite{LYZZ} developed a new KAM scheme, enabling them to obtain a sharp exponential upper bound for spectral gaps with any $k\in\Z^d\backslash\{0\}$, applicable to both discrete and continuous models.


The lower bound estimation is more challenging, particularly for general analytic potentials. Indeed, it's already very difficult to prove the Cantor structure of the spectrum \cite{abd,gs,mp,eli}. In this context, we focus on the AMO case. Leguil-You-Zhao-Zhou \cite{LYZZ} demonstrated that if $\alpha \in DC$, for any $0 < \lambda < 1$, there exists $\xi \gg 1$ such that: $$ |G_k(\lambda)| \geq C(\lambda, \alpha) \lambda^{\xi |k|}, \quad \forall k \in \Z \backslash \{0\}.$$
Fedotov \cite{Fe} recently used the monodromization method to show that for $0 < \xi < 1/2$ and $\lambda$ sufficiently small, $\lambda < e^{-c/\xi}$, the lower bound is: $$ |G_k(\lambda)| \geq \lambda^k e^{O(1/h)}, \quad 1 \leq k \leq O(\ln(1/\lambda)).$$
Our contribution is that, if $\alpha \in DC$, we prove that for all $0 < \lambda < 1$: \begin{equation}\label{ges} |G_k(\lambda)| \geq C(\lambda, \alpha) \lambda^{|k|}, \qquad \forall k \in \Z \backslash \{0\},\end{equation} which provides a lower bound  for  $\lambda \neq 1$ (not only small $\lambda$) and {\it all}  gaps.  We remark that Fedotov's result \cite{Fe} is applicable to any irrational $\alpha$, but one cannot expect \eqref{ges} for any irrational $\alpha$ \footnote{The authors would like to thank A. Avila for useful discussions on this.}.

\subsection{Optimal arithmetic Anderson localization and optimal reducibility} Our second result is related to sharp phase transition for quasiperiodic operators.

\subsubsection{Optimal arithmetic Anderson localization}
We consider the following quasiperiodic long-range operator on $\ell^2(\Z^d)$:
\begin{equation}\label{1.2}
(L_{V,\lambda,\alpha,x}u)_n=\sum\limits_{k\in\Z^d}\hat{V}_ku_{n-k}+2 \lambda \cos 2\pi(x+\langle n,\alpha\rangle)u_{n}, \ \ n\in \Z^d,
\end{equation}
where $V(\theta)=\sum_{k\in\Z^d}\hat{V}(k)e^{2\pi i\langle k,\theta\rangle}\in C^\omega(\T^d,\R)$, $x\in \T$ is called the phase and $\alpha\in\R^d$ is called the frequency.

Operator \eqref{1.2} has garnered significant interest since the 1980s \cite{aj1, b1, bj02, CD, gyz, jk2, CSZ1, CSZ2}. On one hand, its spectral properties are closely tied to those of its Aubry dual, multi-frequency quasiperiodic Schrödinger operators. On the other hand, \eqref{1.2} includes several well-known quasiperiodic models. In particular, when $V(\theta) = \sum_{i=1}^d 2\lambda^{-1}\cos(2\pi \theta_i)$, \eqref{1.2} reduces to the quasiperiodic Schrödinger operator on $\ell^2(\Z^d)$:
\begin{equation}\label{highschrodinger}
L_{\lambda,\alpha,x}= \Delta +  2\lambda\cos2\pi(x+ \langle n,\alpha\rangle)\delta_{nn'},
\end{equation}
where $\Delta$ is the standard Laplacian on the $\Z^d$ lattice. For $d=1$, this is the celebrated almost Mathieu operator.

In the case $d=1$, Bourgain and Jitomirskaya \cite{bj02} established that for any fixed $\alpha \in DC$ and fixed $\alpha$-Diophantine phase, $L_{V,\lambda,\alpha,x}$ has Anderson localization (AL) for sufficiently large $\lambda$. This result is non-perturbative in the sense that the largeness of $\lambda$ is independent of the Diophantine constant of $\alpha$. However, as they argued: ``current techniques do not extend to operators \eqref{1.2} with multi-frequency potentials ". 
If $d \geq 2$, Jitomirskaya and Kachkovskiy \cite{jk2} proved that for fixed $\alpha \in DC_d(\kappa,\tau)$, $L_{V, \lambda, \alpha, x}$ has pure point spectrum for large enough $\lambda$ and almost every $x$. Under the same assumption, Ge-You-Zhou \cite{gyz} proved that $L_{V, \lambda, \alpha, x}$ has \textit{exponential dynamical localization in expectation} (EDL). An ergodic family $\{H_\omega\}$ is said to have  EDL, if
\begin{equation}\label{edl}
\int_{\omega}\sup\limits_{t\in\R}|\langle\delta_k, e^{-itH_\theta}\delta_{\ell}\rangle|d\mu \leq Ce^{-\gamma|k-\ell|}.
\end{equation}
It is apparent that EDL implies AL for almost every phase. Ge-You \cite{gy} showed that actually $L_{V,\lambda,\alpha,x}$ has AL for any $\alpha$-Diophantine phase. This generalizes Bourgain and Jitomirskaya's result \cite{bj02} to all dimensions. Quite recently, Cao, Shi, and Zhang \cite{CSZ1,CSZ2} provided another proof using multiscale analysis. It is natural to ask what the optimal arithmetic condition on $x$ is. We answer this question as follows:

\begin{Theorem}\label{long}
Assume that $\alpha\in DC_d(\kappa,\tau)$, $V\in C_h^{\omega}(\T^d,\R)$. If  $2\pi h>\delta(\alpha,x)$, then  
$L_{V,\lambda,\alpha,x}$ has Anderson localization for sufficiently large $\lambda\geq \lambda_0(\alpha,V,d,\delta)$. 
\end{Theorem}

\begin{Remark}
If $d=1$, Theorem \ref{long} is non-perturbative, i.e. the largeness of $\lambda$ doesn't depend on $\alpha$. 
\end{Remark}

The optimality of this result is rooted in the well-known arithmetic phase transition of the almost Mathieu operator (AMO). Let's explain this more clearly. One crucial aspect in proving Anderson localization lies in the treatment of resonances eigenvalues of box restrictions that are too close to each other in relation to the distance between the box, leading to the so-called ``small denominators". There are  two distinct types of resonances. One is associated with the frequency: the strength of frequency resonances is measured by the arithmetic parameter
\begin{equation}
\beta(\alpha)=\limsup\limits_{k\rightarrow\infty}-\frac{\ln\|\langle k,\alpha\rangle\|_{\R/\Z}}{|k|},
\end{equation}
where $\|x\|_{\R/\Z}=\inf_{\ell\in\Z}|x-\ell|$. 
 Another is determined by the resonances caused by phases, measured by \begin{equation}
\delta(\alpha,\theta)=\limsup\limits_{k\rightarrow\infty}-\frac{\ln\|2\rho+\langle k,\alpha\rangle\|_{\R/\Z}}{|k|}.
\end{equation}
  Around 1995, Jitomirskaya \cite{J} conjectured the following sharp phase transition result:

\begin{itemize}
\item {\bf Sharp phase transition in frequency:}
\begin{enumerate}
	\item If $|\lambda|>e^{\beta(\alpha)}$, $H_{\lambda,\alpha,\theta}$ has Anderson localization  for $\alpha$-Diophantine $\theta$.

	\item If $1<|\lambda|<e^{\beta(\alpha)}$, $H_{\lambda,\alpha,\theta}$ has purely singular continuous spectrum for all $\theta$.
\end{enumerate}
\item {\bf Sharp phase transition in phase:}
\begin{enumerate}
	\item If $|\lambda|>e^{\delta(\alpha,\theta)}$, $H_{\lambda,\alpha,\theta}$ has Anderson localization for Diophantine $\alpha$.

	\item If $1<|\lambda|<e^{\delta(\alpha,\theta)}$, $H_{\lambda,\alpha,\theta}$ has purely singular continuous spectrum for all irrational $\alpha$.
\end{enumerate}
\end{itemize}

One can consult \cite{aubryandre, as, gor, JS94, simon} for the history of this conjecture. Concerning its solution, consult \cite{aj, ly, yzhou} for partial results. The measure-theoretic version of the frequency part conjecture was resolved by Avila, You, and Zhou \cite{ayz}, who proved singular continuous spectrum for $1<|\lambda|<e^\beta$ and Anderson localization for almost every $\theta$ when $|\lambda|>e^\beta$. The complete solution to part (1) of the frequency conjecture was given by Jitomirskaya and Liu \cite{JLiu} (consult \cite{GYZh2} for another proof based on reducibility). In the phase part, singular continuous spectrum was initially established for $1<|\lambda|<e^{c\delta(\alpha,\theta)}$ with a small constant $c$ \cite{JS94}. If $\delta(\alpha,\theta) = 0$, localization was first proved by Jitomirskaya \cite{J99}, while the full conjecture was ultimately resolved by Jitomirskaya and Liu \cite{JLiu2}, demonstrating the optimality of Theorem \ref{long} since for AMO we have $2\pi h=\ln|\lambda|$.

\subsubsection{Sharp phase trasition}

Theorem \ref{long} requires $\lambda$ to be large, indeed our methods also allowed non-perturbative results for almost Mathieu operator, and even for more general type I operators introduced in \cite{gjy}. Note that recently it became possible to prove pure point spectrum in a non-constructive way, employing reducibility for the dual model, an idea initially developed in \cite{yzhou} and first realized by Avila-You-Zhou \cite{ayz} (see \cite{jk2} for an alternate proof). 
 Indeed, when studying the eigenequation of a quasiperiodic Schrödinger operator $ H_{V,\alpha,\theta}u=Eu $, it yields the Schrödinger cocycle ($\alpha, A_E$), where $$ A_E(x) = \begin{pmatrix} E - V(x) & -1 \\ 1 & 0 \end{pmatrix}. $$ Thus, dynamical methods naturally come into play.  Here, given  any $A \in C^\omega(\mathbb{T}^d, \SL(2,\R))$ and rationally independent $\alpha \in \mathbb{R}^d$, a quasiperiodic cocycle $(\alpha, A)$ is defined as follows:
$$
(\alpha,A)\colon \left\{
\begin{array}{rcl}
	\T^d \times \R^2 &\to& \T^d \times \R^2\\[1mm]
	(\theta,v) &\mapsto& (\theta+\alpha,A(\theta)\cdot v)
\end{array}
\right.  .
$$
Recall that $(\alpha,A_1)$ is conjugated to $(\alpha,A_2)$ if there exists $B\in C^\omega(\T^d,\PSL(2,\R))$ such that $$B^{-1}(\theta+\alpha)A_1(\theta)B(\theta) = A_2(\theta).$$ Then $(\alpha,A)$ is  reducible if  it is conjugated to  the constant.

The shortcoming of this approach is that it does not provide a description of the localization phases. In contrast, the proof in \cite{JLiu, JLiu2} is constructive, based on new developments of localization techniques for the AMO, yielding sharp arithmetic phase transitions. More recently, an arithmetic version of Aubry duality was established in \cite{gy}, leading to a new proof of Jitomirskaya's sharp phase transition in frequency via reducibility \cite{gyzh}. As pointed out by Jitomirskaya-Liu \cite{JLiu2}: ``It appears to present a potential for an alternative proof of sharp transition in phase, which would be quite interesting''. In this paper, we answer their question, thus providing a different proof of the sharp transition conjecture in phase.

\begin{Theorem}\label{main}
	Let $\alpha \in DC$. If $|\lambda|>e^{\delta(\alpha,\theta)}$, then $H_{\lambda,\alpha,\theta}$ has Anderson localization.
\end{Theorem}
\begin{Remark}
Involving recent progress on ``structured analysis of the dual cocycle'' \cite{gj1,gjyz}, Theorem \ref{main} can be also extended to any type I operators, thus beyond the almost Mathieu, for simplicity, we only present the almost Mathieu result here.
\end{Remark}

\subsubsection{Optimal arithmetic reducibility}
Theorem \ref{main} depends on the optimal reducibility result, which we now try to explain. 
There are two key milestones in the reducibility of quasiperiodic $\SL(2,\mathbb{R})$ cocycles $(\alpha, A)$ in the perturbative regime, i.e., $A(\cdot) \in C_h^\omega(\mathbb{T}, \SL(2,\mathbb{R}))$ is sufficiently close to constant. Dinaburg and Sinai \cite{DS} showed that if $\alpha\in DC_d(\kappa,\tau)$ and the fiber rotation number $\rho=\rho(\alpha, A)$ belongs to $DC_{\alpha} (\gamma', \tau')$, where \begin{equation}\label{diorho} DC_{\alpha}(\gamma',\tau'):=\left\{\rho \in \R: \| 2\rho- \la n,\alpha \ra \|_{\R/\Z} > \frac{\gamma'}{|n|^{\tau'}},\quad \forall\, n\in\Z^d\backslash\{0\} \right\}, \end{equation} then $(\alpha, A)$ is reducible. Later, based on the resonances cancellation developed by M\"oser and P\"oschel \cite{mp}, Eliasson \cite{eli} showed that if $\alpha\in DC_d(\kappa,\tau)$, and $\rho \in \cup_{\gamma'>0} DC_{\alpha} (\gamma', \tau')$ or $\rho$ is rational with respect to $\alpha$, then $(\alpha, A)$ is reducible. Therefore, the question remains: what is the optimal arithmetic condition for reducibility of $\rho$?
During the 2015 conference ``Almost Periodic and Other Ergodic Problems", Jitomirskaya conjectured that the optimal arithmetic condition for reducibility is $2\pi h > \delta(\alpha, \rho)$. In this paper, we confirm her conjecture with the following theorems:

\begin{Theorem}\label{reducibility main}
Assume $\alpha\in DC_d$,  $R\in \SL(2,\R)$, $A \in C_h^\omega(\mathbb{T}, \SL(2,\mathbb{R}))$ with $ 2\pi h>2\pi \tilde{h}>\delta(\alpha,\rho)$, Then there exists $\epsilon(\alpha,R,h,\tilde{h})$ such that if $\|A(\cdot)-R\|_h<\epsilon$, then $(\alpha,A)$ is reducible.
\end{Theorem}

Before explaining its optimality, let's make some remarks. This theorem, based on KAM techniques, is inherently perturbative. For \( d \geq 2 \), the smallness requirement may not be independent of \( \alpha \) due to a counterexample by Bourgain \cite{B02}. However, in the one-dimensional case (\( d = 1 \)), we can derive more results. Recall that \( (\alpha, A) \) is subcritical if there is a uniform subexponential bound on the growth of \( \|A_n(z)\| \) through some band \( | \Im z | < \delta \). Then, we have the following:

\begin{Theorem}\label{reducibility main2}
Assume \( \alpha \in \text{DC} \) and \( (\alpha, A) \) is subcritical in the strip \( \{ |\Im z| < h \} \). If 
$
2 \pi h > \delta(\alpha, \rho),
$
then \( (\alpha, A) \) is reducible.
\end{Theorem}

As an immediate corollary, we have  the  following arithmetic  transition result for almost Mathieu cocycles $S_E^{\lambda}(\theta)=\begin{pmatrix} E-2\lambda\cos 2\pi (\theta)&-1\\ 1&0\end{pmatrix}$:
 \begin{Corollary}\label{amo reducibility}
Assume $\alpha\in DC$ and $0<\lambda<1$. Then the following holds:
\begin{enumerate}
\item If $-\ln \lambda >\delta(\alpha,\rho)$, then  $(\alpha, S_E^{\lambda})$  is analytically reducible;
\item If $-\ln \lambda <\delta(\alpha,\rho)$, then  $(\alpha, S_E^{\lambda})$ is not analytically  reducible.
\end{enumerate}
\end{Corollary}

Corollary \ref{amo reducibility} (2) was essentially proved by Jitomirskaya-Liu \cite{JLiu2} using duality, while we will give a self-contained new proof here. Then Corollary \ref{amo reducibility} immediately establishes the optimality of Theorem \ref{reducibility main}. Indeed, we can understand its optimality in another way. By the fibred Anosov-Katok construction \cite{KXZ,LYZ}, one can construct dense sets of $f \in C^\omega(\T^d,sl(2,\R))$ with $\|f\|_h < \epsilon$ and $0 < 2\pi h < \delta(\alpha,\rho)$, such that $(\alpha,Ae^f)$ exhibits sub-linear growth, implying irreducibility, illustrating the necessity of the condition $2\pi h > \delta(\alpha,\rho)$.

We also remark that Fayad-Krikorian \cite{FK} generalized Eliasson's result \cite{eli} to the $C^{\infty}$ setting, i.e., the cocycle map $A\in C^\infty(\T^d, SL(2,\R))$. Considering Herman-Yoccoz's optimal linearization result for circle diffeomorphisms \cite{He79,Yo84}, their work was also regarded as optimal in the $C^{\infty}$ setting \footnote{The authors would like to thank B. Fayad for useful discussions on this. }. However, their focus was on finitely many rotation number resonances, whereas our Theorem \ref{reducibility main} addresses infinite number of resonances.

\subsection{Stratified, optimal asymptotic  growth of extended eigenfunction} Our final result focuses on the asymptotic behavior of extended eigenfunctions (denoted by $u_E(n)$), which corresponds to the generalized eigenfunction in the absolutely continuous spectrum, as found in the physics literature.

These eigenfunctions exhibit a connection to the rotation number. According to Eliasson's reducibility theory \cite{eli}, if the rotation number $\rho(E) := \rho(\alpha, S_E^\lambda)$ is rationally related to the frequency, the growth of $u_E(n)$ is linear. Conversely, when $\rho(E)$ is irrational, the eigenfunction's growth is sub-linear, a property that is shown to be optimal \cite{KXZ}. Theorem \ref{reducibility main} further states that if the condition $2\pi h > \delta(\alpha,\rho(E))$ is met, where $\delta$ is an arithmetic invariant, the growth is bounded. We now show that the growth of $u_E(n)$ is actually stratified, reflecting the arithmetic strength of the rotation number.

\begin{Theorem}\label{high}
Assume $\alpha\in DC_d$ and $V\in C_h^\omega(\T^d,\R)$.  Then there exists $\epsilon_0(\alpha,V)$ such that if $\|V\|_h\leq\epsilon_0$, then 
\begin{equation}\label{upper-bound}
\limsup\limits_{n\rightarrow \infty}\frac{\ln\|u_E(n)\|}{\ln |n|}\leq 1-\min\{\frac{2\pi h}{\delta},1\}.
\end{equation}

\end{Theorem}

The natural question then arises: does the upper bound in \eqref{upper-bound} constitute an optimal growth rate for extended eigenfunctions, or does the arithmetic quantity $\delta(\alpha,\rho)$ truly govern their growth?

On the other hand, as described in \cite{JLiu}, a captivating question in solid-state physics is to comprehend the hierarchical structure of spectral features for operators that represent 2D Bloch electrons in a perpendicular magnetic field, linked to the continued fraction expansion of the magnetic flux. Earlier results, like the hierarchical structure, were derived in the works of Sinai \cite{Sinai}, Helffer-Sjostrand \cite{hs}, and Buslaev-Fedotov \cite{f}.  Recently,  for supercritical AMO, Jitomirskaya-Liu \cite{JLiu,JLiu2} have uncovered a universal hierarchical structure in these localized eigenfunctions, not only determining the optimal asymptotics but also revealing the hierarchical pattern of local maxima through a universal function $f(k)$ derived from resonance phases.

The parallel question in the absolutely continuous spectrum regime is whether a similar universal hierarchical structure exists for extended eigenfunctions. We will address these two questions in  the following:

Assume that $0 < \lambda < 1$ and $h_\lambda = -\ln \lambda$. For any $0<\e_0<h_{\lambda}$, we let $\{\ell_i\}_{i=1}^\infty$ be a sequence such that
$$
\|2\rho(E)-\ell_i\alpha\|_{\R/\Z}\leq e^{-\e_0|\ell_i|}.
$$ 
Let $\eta_i\in (0,\infty]$ be such that
$$
\sin2\pi(2\rho(E)+\ell_i\alpha)=e^{-\eta_i|\ell_i|}.
$$
Fix $\e>0$, for $e^{\frac{1}{256}\e h_{\lambda} |\ell_j|}\leq|n|\leq e^{\frac{1}{256}\e h_{\lambda} |\ell_{j+1}|}$,  we define
\begin{equation*}
f(n)=\begin{cases}\max\{1-\frac{|{\ell_j}|h_\lambda}{\ln|n|},0\}& e^{\frac{1}{256}\e h_{\lambda}|\ell_i|} \leq |n|<e^{\eta_{j}|\ell_j|}\\
\max\{\frac{\ln |\sin 2\pi n(\rho(E)+\ell_j\alpha/2)|+|\ell_j|(\eta_{j}-h_\lambda)}{\ln|n|},0\}& e^{\eta_{j}|\ell_j|}<|n|<e^{\frac{1}{256}\e h_{\lambda} |\ell_{j+1}|}.
\end{cases}
\end{equation*}

\begin{center}
  \begin{tikzpicture}[scale=0.6]\label{Fig2}
  \draw [->](2,0)--(23.6,0);
  \draw [->](2,0)--(2,5.5);
  \draw [thick] plot [smooth] coordinates {(4,1.5)(4.5,4)(5.5,2)(7,1.65)(8,1.55)(9,1.5)};
  \draw [thick] plot [smooth] coordinates {(11,1.5)(12,3.5)(15,2)(17,1.65)(22,1.5)};
  \draw [dashed] (4,5)--(4,1.5);
  \draw [dashed] (4.5,5)--(4.5,1.5);
  \draw (4.5,4)--(6,4);
  \node [right] at (6,4){$1-\frac{h_\lambda}{\eta_s}$};
  \draw [dashed] (2,5)--(4,5);
  \draw [dashed] (4,5)--(9,5);
  \draw [dashed] (9,5)--(11,5);
  \node  at (10,3.25){ $ \cdots\cdots $};
  \draw [dashed] (11,5)--(22.5,5);
  \draw [dashed] (9,5)--(9,1.5);
  
  \draw [dashed] (11,5)--(11,1.5);
  \draw [dashed] (12,5)--(12,1.5);
  \draw (12,3.5)--(16,3.5);
  \node [right] at (16,3.5){$1-\frac{h_\lambda}{\eta_{s'}}$};
  \draw [dashed] (22,5)--(22,1.5);
  \draw [dashed] (2,1.5)--(4,1.5);
  \draw [dashed] (4,1.5)--(9,1.5);
  \draw [dashed] (9,1.5)--(11,1.5);
  \draw [dashed] (11,1.5)--(22.5,1.5);
  \node [above] at (4,0){$e^{2\pi h_\lambda|\ell_{s}|}$};
  \node [below] at (5,0){$e^{2\pi \eta_{s} |\ell_{i}|}$};
  \node [above] at (8.2,0){$e^{2\pi h_\lambda|\ell_{s+1}|}$};
  \node [above] at (11,0){$e^{2\pi h_\lambda|\ell_{s'}|}$};
  \node [below] at (13,0){$e^{2\pi \eta_{s'} |\ell_{s'}|}$};
  \node [above] at (21.2,0){$e^{2\pi h_\lambda|\ell_{s'+1}|}$};
  \node [above] at (24,0){$n$};
  \node [left] at (2,5){$1$};

  \node [left] at (2,1.5){$0$};
  \node [left] at (2,5.8){$f(n)$};
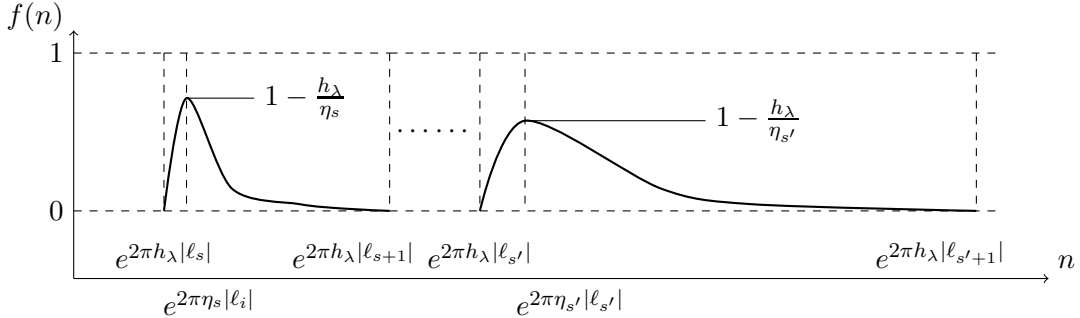
  \node [below] at (9,-2){Figure 1. Picture of $ f(n) $ in two resonant windows};
  \end{tikzpicture}
  \end{center}

 Let $\phi$ and $\psi$ be two generalized eigenfunctions  of the eigenequation $ H_{V,\alpha,\theta}u=Eu $,  having initial conditions $\begin{pmatrix}\phi(1)&\psi(1)\\ \phi(0)&\psi(0)\end{pmatrix}=\begin{pmatrix}1&0\\ 0&1\end{pmatrix}$. 
Let $U_E(n)=\begin{pmatrix}\phi(n) & \psi(n) \\ \phi(n-1) &\psi(n-1)\end{pmatrix}$.  Then we have the following: 

\begin{Theorem}\label{amoasy} Let $\alpha \in DC$, $0<\lambda<1$. 
For any $E\in \Sigma_{\lambda,\alpha}$ and $\e>0$, there exists $N$, such that the associated extended eigenfuntion  matrix $U_E(n)$ satisfies
$$
|n|^{ f(n)-\e}\leq \|U_E(n)\|_{HS}\leq |n|^{f(n)+ \e}.
$$
provided $n>N$. In particular,  we have
\begin{eqnarray}
\label{op-gr1} \limsup\limits_{n\rightarrow \infty}\frac{\ln\|U_E(n)\|_{HS}}{\ln |n|}&=& 1-\min\{-\frac{\ln\lambda}{\delta},1\},\\
\label{op-gr2} \liminf\limits_{n\rightarrow \infty}\frac{\ln\|U_E(n)\|_{HS}}{\ln |n|}&=&0.
\end{eqnarray}

\end{Theorem}

\begin{Remark}
In the analysis of lower bounds for the solution, the Hilbert-Schmidt norm $\|U_E(n)\|_{HS}$ is preferred. This selection is justified by the distinct characteristics of gap edges, where one solution displays linear growth while the other remains bounded.
\end{Remark}

In conclusion, though we didn't acquire the full universal hierarchical structure, we managed to obtain a partial structure for $U_E(n)$. To simplify, let's assume $0 < \delta(\alpha, \rho(E)) < \infty$. By definition, there exists a sequence of resonances $\{\ell_i\}_{i=1}^\infty$ such that: $$
e^{-2\pi(\delta-\epsilon)|\ell_i|} \leq \|2\rho(E) - \ell_i\alpha\|_{\R/\Z} \leq e^{-2\pi(\delta+\epsilon)|\ell_i|}.
$$ As depicted in Figure 1, for $e^{\frac{1}{256}\epsilon h_{\lambda} |\ell_j|} \leq |n| \leq e^{\frac{1}{256}\epsilon h_{\lambda} |\ell_{j+1}|}$, the function $f(n)$ increases from 0 to $1 + \frac{\ln\lambda}{\delta(\alpha,\rho(E))}$, then decreases back to 0. Theorem \ref{amoasy} highlights the optimal stratified sub-linear growth pattern and reveals that the extended eigenfunction $U_E(n)$ exhibits self-similarity around each resonance $\ell_i$.

\subsection{Difficulty, Novelty, ideas of the proof}\label{nove}

This paper focuses on the role of resonances in quasiperiodic dynamical systems. A defining characteristic of quasiperiodic dynamics is the profound impact of arithmetic on its behavior, famously manifesting as the ``small denominators" problem in various models, such as circle diffeomorphisms, disc dynamics, quasiperiodic cocycles, surface flows, and Hamiltonian systems in finite and infinite dimensions \cite{arnold, AKN, S14, TP, YouICM}. The primary challenge in managing ``small denominators'' issues lies in addressing various types of resonances.

\subsubsection{Infinitely many rotation number resonances}

We mainly concerned on a significant class of quasiperiodic dynamical systems: analytic quasiperiodic $\SL(2,\mathbb{R})$ cocycles. Though our main analysis applies to any dimension $d$, we concentrate on the one-frequency case, where Avila's groundbreaking global theory for $\SL(2,\mathbb{R})$-cocycles is established \cite{avila2015global}. A central aspect of his theory is the ``Almost Reducibility Conjecture'' (ARC), which posits that $(\alpha,A)$ is almost reducible if it is subcritical \cite{avila, avila1, avila2015global}. Recall that $(\alpha,A)$ is almost reducible if its analytic conjugate class contains a constant. Almost reducibility crucially depends on two types of resonances, or ``small denominators'': one is frequency resonances, measured by $\beta(\alpha)$; the other is determined by resonances caused by the fibred rotation number $\rho$ \footnote{See Section 2 for definition.}, measured by $\delta(\alpha,\rho)$.

Almost reducibility captures the properties of cocycles for which local theory (KAM theory) applies. Fundamentally, classical KAM theory avoids resonances by assuming Diophantine-type conditions on the frequency; the underlying philosophy is that polynomial growth of the ``small denominators'' can be controlled by exponential decay of the Fourier coefficients. Around the 2010s, Avila-Fayad-Krikorian \cite{afk} and Hou-You \cite{houyou} independently developed non-standard KAM techniques to handle infinitely many frequency resonances (i.e., $\beta(\alpha) > 0$), but still managed only a finite number of rotation number resonances.

However, obtaining optimal reducibility results—such as in Theorem \ref{reducibility main}—presents a challenge due to the presence of infinitely many rotation resonances, where \( \delta(\alpha, \rho) > 0 \). This question has remained open for several decades. From the spectral perspective, quantitative almost reducibility has proven to be a powerful tool for investigating the spectral theory of quasiperiodic Schrödinger operators \cite{YouICM, JICM}, both in the subcritical regime \cite{avila, aj1, ayz2, gk, LYZZ, puig} and the supercritical regime \cite{ayz, gy, GYZh2, gyzh, gyz}. However, as we will explore, the main difficulty in achieving optimal spectral results lies in whether \textit{infinitely many} rotation resonances can be precisely managed.

\subsubsection{Structured quantitative almost reducibility}

Hence, the central problem is how to handle these \textit{infinitely many} rotation resonances. We must emphasize that dealing with infinitely many rotation resonances requires a completely new understanding of the dynamics, whereas handling infinitely many frequency resonances depends crucially on understanding continued fractions. To address this challenge, we introduce a new approach called \textit{structured quantitative almost reducibility}, which constitutes the first novelty of this paper.

 Recall that $(\alpha, A)$ is almost reducible means there exist $B_j \in C_{h_j}^\omega(\T^d, \PSL(2, \R))$, $A_j \in \SL(2, \R)$, and $f_j \in C_{h_j}^\omega(\T^d, \sl(2, \R))$ such that $$ B_j^{-1}(\theta + \alpha) A(\theta) B_j(\theta) = A_j e^{f_j(\theta)}  $$ with $\|f_j\|_{h_j}\to 0$.  
 The minimal $n_j$ for which
 $\|2\rho_j -\langle n_j, \alpha \rangle\|_{\mathbb{R}/\mathbb{Z}} \leq \epsilon_j^\frac1{15}$
 where $\epsilon_j=:\|f_j\|_{h_j}$ is denoted as the KAM resonances.  In  \cite{eli}, it was proved that $n_j$ grow fast than any polynomial.  These resonances are of crucial importance in KAM iterations. Given that $A_j$ is symplectic, we have $A_j= M^{-1} \exp\begin{pmatrix} i\rho_j & \nu_j \\ \bar{\nu}_j & -i\rho_j \end{pmatrix} M $ where $M = \frac{1}{2i}\begin{pmatrix} 1 & -i \\ 1 & i \end{pmatrix}$ is the isomorphism between $\SL(2, \R)$ and $SU(1,1)$.

 In practical applications, precise structures and estimates of 
$A_j, B_j$ and $F_j$, 
 are required, which motivates the establishment of the structured quantitative almost reducibility.
  We say $(\alpha, A)$ is  structured quantitative almost reducibility, if we simultaneously have the following: \\
$\mathbf{SQ1}:$ Structure of the conjugations, i.e.,
$$
B_j(\cdot)=\tilde{B}_j(\cdot)R_{\frac{\langle n,\cdot\rangle}{2}}e^{Y_j(\cdot)},
$$
where $R_{\frac{\langle n,\cdot\rangle}{2}}$ is the rotation,   and $e^{Y_j(\cdot)}$ is close to identity, $\|\tilde{B}_j\| \sim e^{o(n)}$.  \\
$\mathbf{SQ2}:$  Structure and quantitative estimates of the  perturbations  $f_j$, i.e. $\|f_j\| \sim e^{-O(N)}$.\\
$\mathbf{SQ3}:$  Structure and quantitative estimates of the constants $A_j$, i.e. $|\nu_j | \leq e^{-O(n)}$\footnote{Sharp exponential decay rate of $\nu_j$ is also obtained.}.

This is a rough version. For the detailed statements, refer to Proposition \ref{reducibility} (local) and Proposition \ref{Amo} (global).
 Structured quantitative almost reducibility plays an important role in addressing an infinite array of rotation number resonances. Our methodology entails investigating the interdependence of the conjugacy $B_j(\cdot)$, constant $A_j$, and perturbation $f_j(\cdot)$ on the resonances $n_j$.
As shown in Section \ref{kam-res}, the KAM resonances are closely linked to the resonances of the rotation number $\rho$. For any $\e_0 > 0$, we define an $\e_0$-resonance of the rotation number $\rho$ as follows: \begin{equation*} \|2\rho -\langle l, \alpha \rangle\|_{\mathbb{R}/\mathbb{Z}} \leq e^{-\e_0|\ell|}. \end{equation*} By arranging these $\epsilon_0$-resonances in increasing order of their magnitudes, $0 < |\ell_1| < |\ell_2| < \cdots$, we will demonstrate that a subset of these rotation number resonances aligns with a subsequence of KAM resonances.

 Before discussing its powerful consequences, let's review the history of almost reducibility that led to our definition.
The first weak almost reducibility result (where the analytic radius 
$h_j \rightarrow 0$) was proved by Eliasson \cite{eli}, without emphasizing its dependence on the KAM resonances 
$n_j$  or the structure. Using almost localization of the dual operator \cite{aj1}, Avila \cite{avila} provided the first strong almost reducibility result (where $h_j \rightarrow h_*>0$) satisfying conditions  ($SQ2$) and ($SQ3$), but the conjugacy was constructed via the Corona Theorem (with Uchiyama estimates), thus lacking a good structure.
Recently, another strong almost reducibility scheme \cite{CCYZ,LYZZ} was developed, which grants a better structure of the conjugacy
 (fulfilling $SQ1$) \cite{gyz}, and an improved exponential decay rate of  $\nu_j$ compared to \cite{avila}. However, 
$(SQ2)$ remains unaddressed in this approach.

Our focus lies in determining the simultaneous satisfaction of $(SQ1)$-$(SQ3)$, which yields significant implications. For instance, by leveraging $(SQ1)$, we can introduce a novel KAM scheme (Theorem \ref{reducibility main}) that guarantees optimal arithmetic reducibility, with the convergence guaranteed by $(SQ2)$ and $(SQ3)$. This, in turn, leads to the optimal arithmetic localization result encapsulated in Theorem \ref{main}.

\subsubsection{Sharp version of quantitative Aubry duality.}
The exponential decay of \( \nu_j \) and \( f_j \) with respect to the resonances is of paramount importance, as it ensures that the ``small denominators" remain exponentially small. However, no optimal lower bound for \( \nu_j \) has been achieved in any of the aforementioned schemes. This kind of lower bound estimate turns out to be crucial in spectral applications.

It was first realized by Puig that the positivity of \( |\nu_j| \) was crucial to the Cantor spectrum. He then used Aubry duality (without requiring quantitative estimates) to show the positivity of \( |\nu_j| \). This is the main idea behind the Dry Ten Martini Problem in the Diophantine case \cite{aj1, PuigC}. Later, as introduced in \cite{ayz2}, quantitative Aubry duality provided a crucial connection between quantitative estimates of \( B_j(\theta) \) and localized eigenfunctions of dual operators. This was pivotal in resolving the non-critical Dry Ten Martini Problem \cite{ayz2}. 

Our contribution lies in developing a refined yet sharp quantitative version of this duality, tailored to the specific structure of \( B_j(\theta) \). The structured properties \( (\text{SQ1}) \) enable an in-depth analysis of dual localized eigenfunctions, ultimately leading to optimal lower bounds for \( |\nu_j| \) (similar to \( (\text{SQ3}) \)). These bounds are fundamentally important for achieving the optimal exponential decay rate of spectral gaps (as demonstrated in Theorem \ref{optimal gap}) and understanding the asymptotic growth patterns of extended eigenfunctions (Theorem \ref{amoasy}).

%
%
%
%
%
%

\section{Preliminaries}

%

\subsection{Quasiperiodic cocycles}

Given $A \in C^\omega(\T^d,{\rm SL}(2,\C))$ and rationally independent $\alpha\in\R^d$, we define the \textit{quasiperiodic cocycle} $(\alpha,A)$:
$$
(\alpha,A)\colon \left\{
\begin{array}{rcl}
\T^d \times \C^2 &\to& \T^d \times \C^2\\[1mm]
(\theta,v) &\mapsto& (\theta+\alpha,A(\theta)\cdot v)
\end{array}
\right.  .
$$
The iterates of $(\alpha,A)$ are of the form $(\alpha,A)^n=(n\alpha,  \mathcal{A}_n)$, where
$$
\mathcal{A}_n(\theta):=
\left\{\begin{array}{l l}
A(\theta+(n-1)\alpha) \cdots A(\theta+\alpha) A(\theta),  & n\geq 0\\[1mm]
A^{-1}(\theta+n\alpha) A^{-1}(\theta+(n+1)\alpha) \cdots A^{-1}(\theta-\alpha), & n <0
\end{array}\right.    .
$$
The {\it Lyapunov exponent} is defined by
$\displaystyle
L(\alpha,A):=\lim\limits_{n\to \infty} \frac{1}{n} \int_{\T^d} \ln \|\mathcal{A}_n(\theta)\| d\theta
$.

The cocycle $(\alpha,A)$ is {\it uniformly hyperbolic} if, for every $\theta \in \T^d$, there exists a continuous splitting $\C^2=E^s(\theta)\oplus E^u(\theta)$ such that for every $n \geq 0$,
$$
\begin{array}{rl}
|\mathcal{A}_n(\theta) \, v| \leq C e^{-cn}|v|, &  v \in E^s(\theta),\\[1mm]
|\mathcal{A}_n(\theta)^{-1}   v| \leq C e^{-cn}|v|, &  v \in E^u(\theta+n\alpha),
\end{array}
$$
for some constants $C,c>0$.
This splitting is invariant by the dynamics, i.e.,
$$A(\theta) E^{*}(\theta)=E^{*}(\theta+\alpha), \quad *=``s"\;\ {\rm or} \;\ ``u", \quad \forall \  \theta \in \T^d.$$

Assume that $A \in C (\T^d, {\rm SL}(2, \R))$ is homotopic to the identity. It induces the projective skew-product $F_A\colon \T^d \times \mathbb{S}^1 \to \T^d \times \mathbb{S}^1$ with
$$
F_A(\theta,w):=\left(\theta+\a,\, \frac{A(\theta) \cdot w}{\|A(\theta) \cdot w\|}\right),
$$
which is also homotopic to the identity.
Thus we can lift $F_A$ to a map $\tF_A\colon \T^d \times \R \to \T^d \times \R$ of the form $\tF_A(\theta,y)=(\theta+\alpha,y+\psi_\theta(y))$, where for every $\theta \in \T^d$, $\psi_\theta$ is $\Z$-periodic.
The map $\psi\colon\T^d \times \T  \to \R$ is called a {\it lift} of $A$. Let $\mu$ be any probability measure on $\T^d \times \R$ which is invariant by $\widetilde{F}_A$, and whose projection on the first coordinate is given by Lebesgue measure.
The number
$$
\rho(\alpha,A):=\int_{\T^d \times \R} \psi_\theta(y)\ d\mu(\theta,y) \ {\rm mod} \ \Z
$$
 depends  neither on the lift $\psi$ nor on the measure $\mu$, and is called the \textit{fibered rotation number} of $(\alpha,A)$ (see \cite{H,JM} for more details).

Given $\phi\in\R$, let $
R_\phi:=
\begin{pmatrix}
\cos2 \pi\phi & -\sin2\pi\phi\\
\sin2\pi\phi & \cos2\pi\phi
\end{pmatrix}$.
If $A\colon \T^d\to{\rm PSL}(2,\R)$ is homotopic to $\theta \mapsto R_{\frac{\la n, \theta\ra}{2}}$ for some $n\in\Z^d$,
then we call $n$ the {\it degree} of $A$ and denote it by $\deg A$.
The fibered rotation number is invariant under real conjugacies which are homotopic to the identity. More generally, if $(\alpha,A_1)$ is conjugated to $(\alpha, A_2)$, i.e., $B(\cdot+\alpha)^{-1}A_1(\cdot)B(\cdot)=A_2(\cdot)$, for some $B \colon \T^d\to{\rm PSL}(2,\R)$ with ${\rm deg} B=n$, then
\begin{equation}\label{rotation number}
\rho(\alpha, A_1)= \rho(\alpha, A_2)+ \frac{\la n,\alpha \ra}2.
\end{equation}
Moreover, immediately from the definition of rotation number, we have
\begin{Lemma}\label{rol}
For any $A\in SL(2,\R)$, there exists a numerical constant $C_2>0$ such that
$$
|\rho(\alpha,B)-\rho(\alpha,A)|<C_2\|B(\cdot)-A\|^{\frac{1}{2}}_0.
$$
In particular, if $A=R_{\phi}$, then we have 
$$
|\rho(\alpha,B)-\phi|<C_2\|B(\cdot)-R_\phi\|_0.
$$
\end{Lemma}

%

\subsection{Almost Mathieu cocycle}

Note that a sequence $(u_n)_{n \in \Z}$ is a formal solution of the
eigenvalue equation $H_{\lambda,\alpha,\theta} u=Eu$ if and only if
it satisfied $$\begin{pmatrix}
u_{n+1}\\u_n\end{pmatrix}=S_{E}^{\lambda}(\theta+n\alpha) \cdot
\begin{pmatrix} u_n\\u_{n-1} \end{pmatrix},$$ where we donote
\begin{eqnarray*}
S_{E}^{\lambda}(\theta)=\left( \begin{array}{ccc}
 E-2\lambda\cos2\pi(\theta) &  -1\cr
  1 & 0\end{array} \right)\in SL(2,\mathbb{R}).
\end{eqnarray*}
 then $(\alpha,S_{E}^{\lambda} )$ can be seen as a  quasiperiodic cocycle, and we
 call it almost Mathieu cocycle.

Denote the spectrum of  $H_{\lambda,\alpha,\theta}$ by
$\Sigma_{\lambda,\alpha}$,  which  is independent of $\theta$ when
$\alpha\in \R\backslash \Q$. The spectral properties of $H_{\lambda,\alpha,\theta}$ and the dynamics of $(\alpha,S_E^{\lambda})$ are closely related by the well-known fact:
 $E\in \Sigma_{\lambda,\alpha}$ if and only if $(\alpha,S_E^{\lambda})$ is \textit{not} uniformly hyperbolic. If $E \in \Sigma_{\lambda,\alpha}$,
then the Lyapunov exponent of almost Mathieu cocycle can be computed
explicitly.

\begin{Theorem}\cite{BJ}\label{bj-formula} 
If $\alpha \in \R\backslash\Q$, $E \in \Sigma_{\lambda,\alpha}$, then  we have
$$L(\alpha,S_E^{\lambda})=\max
\{0,\ln |\lambda|\}.$$
\end{Theorem}

Denote $(n\alpha, \mathcal{A}_{E}^n):= (\alpha, S_{E}^{\lambda})^n$, the iterates of the almost Mathieu cocycle. Then as a consequence of Theorem \ref{bj-formula}, we have the following:
\begin{Lemma}\label{lem:upperbdd}\cite{ayz}
Let $\alpha\in \R\backslash\Q$, $\lambda>1$. For any small 
$\epsilon>0$, there exist $N(\epsilon,\lambda,\alpha)<\infty$,
 such that for any $|m|>N(\epsilon,\lambda, \alpha)$, $E \in \Sigma_{\lambda,\alpha}$, we have
$$\sup_{\theta\in \T}\frac{1}{|m|}\ln\| \mathcal{A}_{E}^m(\theta)\|< \ln\lambda +\epsilon.$$
\end{Lemma}

\subsection{Aubry duality}
Consider the fiber direct integral,
$$
\mathcal{H}:=\int_{\T}^{\bigoplus}\ell^2(\Z)d\theta,
$$
which, as usual, is defined as the space of $\ell^2(\Z)$-valued, $L^2$-functions over the measure space $(\T,d\theta)$.  The extensions of the 
Sch\"odinger operators  and their long-range duals to  $\mathcal{H}$ are given in terms of their direct integrals, which we now define. 
Let $\alpha\in\T$ be fixed. Interpreting $H_{V,\alpha,\theta}$ as fibers of the decomposable operator,
$$
H_{V,\alpha}:=\int_{\T}^{\bigoplus}H_{V,\alpha,\theta}d\theta,
$$
then the family $\{H_{V,\alpha,\theta}\}_{\theta\in\T}$ naturally induces an operator on the space $\mathcal{H}$, i.e. , 
$$
(H_{V,\alpha} \Psi)(\theta,n)= \Psi(\theta,n+1)+ \Psi(\theta,n-1) +  V(\theta+n\alpha) \Psi(\theta,n).
$$

Similarly,  the direct integral of long-range operator  $L_{V,\alpha,x}$,
denote as $L_{V,\alpha}$, is given by 
$$
(L_{V,\alpha}  \Psi)(x,n)=  \sum\limits_{k\in\Z^d} \hat{V}_k \Psi(x,n+k)+  2\cos2\pi (x+\langle n,\alpha\rangle) \Psi(x,n),
$$
where $\hat{V}_k$ is the $k$-th Fourier coefficient of $V(\theta)$.
Let  $U$ be the following operator on $\mathcal{H}:$
$$
(\mathcal{U}\phi)(\eta,m)=\sum_{n\in\Z}\int_{\T^d}e^{2\pi i\langle m,\theta\rangle}e^{2\pi in(\langle m,\alpha\rangle+\eta)}\phi(\theta,n)d\theta.
$$
Then direct computations show that $U$ is unitary and satisfies 
$U H_{V,\alpha} U^{-1}=L_{V,\alpha},$ and the quasiperiodic  long-range operator  $L_{V,\alpha,x}$ is called the dual operator of $H_{V,\alpha,\theta}$ \cite{gjls}.

\subsection{Algebraic lemma}

Recall $sl(2,\R)$ is the set of $2\times 2$ matrices  of the form
$$\left(
\begin{array}{ccc}
 x &  y+z\\
 y-z &  -x
 \end{array}\right)$$
 where $x,y,z\in \R.$ It is
 isomorphic to $su(1,1)$, the group of matrices of the form
$$\left(
\begin{array}{ccc}
 i t &  \nu\\
\bar{ \nu} &  -i t
 \end{array}\right)$$
 with $t\in \R$, $\nu\in \C$. The isomorphism between $sl(2,\R)$ and
 $su(1,1)$ is given by $B\rightarrow M B M^{-1}$ where
$$M=\frac{1}{2i}\left(
\begin{array}{ccc}
 1 &  -i\\
 1 &  i
 \end{array}\right).$$
 Direct calculation shows that
 $$M\left(
\begin{array}{ccc}
 x &  y+z\\
 y-z &  -x
 \end{array}\right)M^{-1}=\left(
\begin{array}{ccc}
 i z &  x-iy\\
x+iy &  -i z
 \end{array}\right).$$
We state the following
lemma concerning the diagonalization of elliptic matrices in $\mathrm{su}(1,1)$. The
diagonalizing conjugation given by the lemma is of optimal norm.

\begin{Lemma}\cite{KXZ}\label{nf}
Let the matrix
$$\tilde{A}=\begin{pmatrix} it & \nu \\ \bar{\nu} & -it \end{pmatrix} \in \mathrm{su}(1,1)$$
satisfy $\det \tilde{A}> 0$. Then, calling $\rho=\sqrt{\det \tilde{A}}$, we have
$$
U^{-1}\tilde{A}U=\begin{pmatrix} i\rho & 0 \\ 0 & -i\rho\end{pmatrix},
$$
where
\begin{equation*} 
\begin{array}{r@{}l}
U &=
(\cos 2\varphi)^{-\frac{1}{2}}
\begin{pmatrix}
e^{i\phi} & 0 \\ 0 & e^{- i\phi}
\end{pmatrix}
\begin{pmatrix}\cos \varphi & \sin \varphi \\ \sin\varphi & \cos\varphi
\end{pmatrix}
\begin{pmatrix}
e^{-i\phi} & 0 \\ 0 & e^{i\phi}
\end{pmatrix}
\\
&= (\cos 2\varphi)^{-\frac{1}{2}}
\begin{pmatrix}\cos\varphi &e^{2i\phi} \sin\varphi \\ e^{-2i\phi}\sin\varphi & \cos\varphi
\end{pmatrix}.
\end{array}
\end{equation*}
Here $2\phi = \arg \nu-\frac{\pi}{2}$  and $\varphi \in (-\frac{\pi}{2}, \frac{\pi}{2})$ satisfies
\begin{equation*}
2\varphi = - \arctan \frac{|\nu|}{\sqrt{t^{2}-|\nu|^{2}}}
\end{equation*}
In addition we have
\begin{eqnarray}\label{Dn1}
\|U\|^{2} &=& \frac{(1-\tan \varphi)^2}{1-\tan^2 \varphi } =\frac{|t|+|\nu|}{\rho}.
\end{eqnarray}
\end{Lemma}

As a consequence, we have the following:

\begin{Lemma}\label{Dominate matrix}
Assume that $$MAM^{-1}=exp \left(
\begin{array}{ccc}
 i t &  \nu\\
\bar{ \nu} &  -i t
 \end{array}\right)\in SU(1,1)$$ with $spec\{A\}=\{e^{ i\xi},e^{- i\xi}\}$, $\xi\in \R$, then there exists $U\in SL(2,\R)$, such that
$$
U^{-1}AU=R_{\xi},
$$
with the following estimates:
\begin{enumerate}
\item If  $|\frac{2\nu}{\xi}|\leq 1$, then 
$ \|U-id\|\leq |\frac{\nu}{\xi}|.$
\item otherwise, we have 
$
\|U\|^2\leq \frac{4|\nu|}{\xi}.
$
\end{enumerate}
Moreover, there exists $U' \in SL(2,\C)$ such that 
$$
U'^{-1}AU'=\begin{pmatrix}e^{ i \xi}&\nu' \\ 0&e^{- i \xi}\end{pmatrix},
$$
with
\begin{equation}\label{unitary}
\|U'\|\leq 2, \ \ |\nu'|\leq 4|\nu|.
\end{equation}
\end{Lemma}
\begin{pf}
Note $(1)$ follows direct construction of $U$ in Lemma \ref{nf}, one can consult  Corollary 2.1 in \cite{KXZ} (see also Lemma 4.1 in \cite{gyz}) for details. 
Now if $|\frac{2\nu}{\xi}|>1$, it follows that 
\begin{equation}\label{est}
|t| + |\nu| <4 |\nu|,
\end{equation}
then (2) follows from \eqref{est} and \eqref{Dn1} of Lemma \ref{nf}, while \eqref{unitary} follows from  \eqref{est} and Schur's lemma.

\end{pf}

\section{Local structured quantitative  almost reducibility}

In this section, we introduce the concept of structured quantitative almost reducibility within the perturbated regime. By emphasizing the ``structure" and providing precise estimates, we aim to offer a novel approach to manage resonances. This method will enable us to achieve various sharp dynamical and spectral applications.

First we recall the following strong almost reducibility result:
\begin{Proposition}\label{iteration}\cite{CCYZ, LYZZ}
Let $\alpha\in DC_d(\kappa,\tau)$. Suppose that $A\in \SL(2,\R)$, $f\in C^\omega_h(\T^d,$ $\sl(2,\R))$. Then for any $h_+<h$, there exists numerical constant $C_0$,  and constant $D_0=D_0(\kappa,\tau,d)$ such that if
\begin{align*}
\|f\|_h\leq \epsilon\leq \frac{D_0}{\|A\|^{C_0}}(h-h_+)^{C_0\tau},
\end{align*}
then there exists $B\in C_{h_+}^\omega(\T^d,\PSL(2,\R))$, $A_+\in \SL(2,\R)$ and $f_+\in C_{h_+}^\omega(\T^d,\sl(2,\R))$ such that
$$
B^{-1}(\theta+\alpha)Ae^{f(\theta)}B(\theta)=A_+e^{f_+(\theta)}= M^{-1}exp \left(
\begin{array}{ccc}
 i t^{+} &  \nu^{+}\\
\bar{ \nu}^{+} &  -i t^{+}
 \end{array}\right)M e^{f_+(\theta)}.
$$
More precisely, let $spec(A)=\{e^{ i\xi},e^{- i\xi}\}$, $N=\frac{2}{h-h_+} | \ln \epsilon |$, then we can distinguish two cases:
\begin{itemize}
\item (Non-resonant case)   if for any $n\in \Z^{d}$ with $0< |n| \leq N$, we have
$$
\| 2\xi - <n,\alpha> \|_{\R/\Z} \geq \epsilon^{\frac{1}{15}},
$$
then
$$\| B-id\|_{h_+}\leq \epsilon^{\frac{1}{2}} , \quad   \|f_{+}\|_{h'}\leq \epsilon e^{-2 \pi N(h-h')}, \quad \forall h'\leq h_+.$$
Moreover, $\|A_+-A\|<2\epsilon$.
\item (Resonant case) if there exists $n_\ast$ with $0< |n_\ast| \leq N$ such that
\begin{equation}\label{reso-condi}
\| 2\xi- <n_\ast,\alpha> \|_{\R/\Z}< \epsilon^{\frac{1}{15}},
\end{equation}
then there exists $P\in SL(2,\R)$ such that  $P^{-1}AP= R_{\xi}$ with estimate
$$\|P\|\leq \frac{|n_*|^{\tau}}{\kappa},$$
and the conjugacy $B(\theta)$ takes the form $B(\theta)=B'(\theta) R_{\frac{\langle n_*,\theta\rangle}{2}}= Pe^{Y(\theta)}R_{\frac{\langle n_*,\theta\rangle}{2}}$ with 
\begin{equation}\label{res1}
\deg B=n_*, \qquad  \|Y\|_{h_+}\leq \epsilon^{\frac{1}{2}}.
\end{equation} Moreover, we have 
$$ \|f_{+}\|_{h_+}< \epsilon e^{-h_+\epsilon^{-\frac{1}{18\tau}}},   \qquad |\nu^{+}|\leq \epsilon^{\frac{15}{16}}e^{-2\pi|n_*|h} .$$
\end{itemize}
\end{Proposition}

\subsection{Local structured quantitative  almost reducibility} \label{lar}

By Proposition \ref{iteration}, we can easily prove that for any cocycle  $(\alpha, A_0e^{f_0(\cdot)})$, if $f_0$ is small enough, then $(\alpha, A_0e^{f_0(\cdot)})$ is almost reducible \cite{CCYZ,LYZZ}, i.e.  there exist  sequence of $B_j\in C_{\tilde{h}}^\omega(\T^d, \PSL(2,\R))$, $A_j\in \SL(2,\R)$ and $f_j\in C^\omega_{\tilde{h}}(\T^d,\sl(2,\R))$ such that
$$
B_j^{-1}(\theta+\alpha)A_0e^{f_0(\theta)}B_j(\theta)=A_je^{f_j(\theta)}.$$
Here, $n_*$ satisfying \eqref{reso-condi} is called the KAM resonance. Indeed, by \eqref{res1}, we can label all the resonances  as 
$$
I_j=\{n_{i-1}=\deg{B_{i}}-\deg{B_{i-1}}|\deg{B_{i}}\neq \deg{B_{i-1}},1\leq i\leq j\}:=\{n_{i_1},\cdots,n_{i_j}\},
$$ 
i.e. $i_j \leq j-1$ is the last resonance up to  $j$-th KAM step. By  the construction
we have  
\begin{equation}\label{resos}
\deg B_j=n_{i_1}+\cdots+n_{i_j}.
\end{equation} If there are only finite number KAM resonances, we formally set $|n_{i_{j+1}}|=\infty$.  
To prove the structured quantitative  almost reducibility, the first step is to  relate directly the estimate of $B_j$ to the KAM resonances $n_{i_j}$, and explore the structure of $B_j$, which is also directly related to  $n_{i_j}$.  Now we state the result as follows:

\begin{Proposition}\label{reducibility}
Let $\alpha\in DC_d(\kappa,\tau)$,  $A\in C_{h}^\omega(\T^d,\SL(2,\R))$ with $ h> \tilde{h}>0$, $R\in \SL(2,\R)$. 
If 
 $$\|A(\cdot)-R\|_{h}\leq \epsilon \leq \frac{D_0(\kappa,\tau,d)}{\|A\|^{C_0}}(h-h_+)^{C_0\tau},  $$
then there exist $B_j\in C_{\tilde{h}}^\omega(\T^d, \PSL(2,\R))$,  $f_j\in C_{\tilde{h}}^\omega(\T^d, \sl(2,\R))$  such that 
$$
B_j^{-1}(\theta+\alpha)A(\theta)B_j(\theta)=A_je^{f_j(\theta)}=M^{-1}exp \left(
\begin{array}{ccc}
 i t^j &  \nu^j\\
\bar{\nu}^j &  -i t^j
 \end{array}\right)Me^{f_j(\theta)}
$$
with estimates
\begin{eqnarray}
 \label{es2}
|\nu^j| &\leq& 2 \epsilon_{i_j}^{\frac{15}{16}}e^{-2\pi|n_{i_j}|\tilde{h}},\\
   \label{es0} \|f_j\|_{\tilde{h}} &\leq&  \epsilon_{i_j}e^{- \frac{\pi}{8} |n_{i_{j+1}}|(h-\tilde{h})}.
\end{eqnarray}
 where $\epsilon_j\leq \epsilon^{2^j}$.   Moreover,  the conjugation $B_j(\cdot)$ takes the form 
  $$B_j(\theta)=\tilde{B}_j(\theta)R_{\frac{\langle n_{i_j},\theta\rangle}{2}}e^{Y_j(\theta)},$$
 with estimate 
  \begin{eqnarray} \label{es1}
 \|Y_j\|_{\tilde{h}} &<&   e^{-2\pi|n_{i_j}|^4\tilde{h}},\\
  \label{es3}
\|\tilde{B}_j\|_{\tilde{h}}&<&C(\alpha)|n_{i_{j}}|^{\tau}e^{2 \epsilon_{i_{j-1}}^{\frac{1}{18\tau}} |n_{i_j}|},\\
   \label{es4}
\|B_j\|_0&<&C(\alpha)|n_{i_j}|^{2\tau},\\
\label{es5} 
|\deg B_j - n_{i_j} | &\leq&  2\epsilon_{i_{j-1}}^{\frac{1}{18\tau}} |n_{i_j}|.
\end{eqnarray} 
\end{Proposition}

\begin{pf} We prove the result by iteration, and divide the proof into three steps: \\

\textbf{Step 1: One step of KAM:}

First  by simple implicit functional theorem, we can rewrite $A(\theta)=Re^{f(\theta)}$ with $\|f\|_{h}\leq\epsilon$. 
Let $A_0=R$, $f_0(\theta)=f(\theta)$, $h_0=h$, $\epsilon_0=\epsilon$,  assume that we are at the $(j+1)$-th KAM step, i.e. we already construct $B_j\in C^\omega_{h_{j}}(\T^d,\PSL(2,\R))$ such that
$$
B_{j}^{-1}(\theta+\alpha)A_0e^{f_0(\theta)}B_{j}(\theta)=A_{j}e^{f_{j}(\theta)},
$$
where $A_j\in \SL(2,\R)$ with two eigenvalues $e^{\pm i\xi_j}$ and
\begin{eqnarray}\label{estimate2}
\|B_{j}\|_{h_j}&\leq& C(\alpha)e^{4\pi h(1-\frac{1}{4^j})|n_{i_j}|},\\ \label{estimate3}
\|B_{j}\|_0&\leq& C(\alpha)\prod_{k\leq j}(1+\frac{1}{4^k})|n_{i_j}|^{2\tau},\\ \nonumber
\|f_j\|_{h_j}&\leq& \epsilon_j\leq \epsilon_0^{2^j},
\end{eqnarray}
then we define
$$
h_j-h_{j+1}=\frac{h-\tilde{h}}{4^{j+1}}, \ \ N_j=\frac{2|\ln\epsilon_j|}{h_j-h_{j+1}}.
$$

By our selection of  $\epsilon_0$, one can check that
\begin{equation}\label{itera1}
\epsilon_j \leq \frac{D_0}{\|A_j\|^{C_0}}(h_j-h_{j+1})^{C_0\tau}.
\end{equation}
 Indeed, $\epsilon_j$ on the left side of the inequality decays at least super-exponentially with $j$, while $(h_j-h_{j+1})^{C_0\tau}$ on the right side decays exponentially with $j$.

Note that $\eqref{itera1}$ implies that  Proposition \ref{iteration} can be applied iteratively, consequently one can construct
$$
\bar{B}_j\in C^\omega_{h_{j+1}}(\T^d,\PSL(2,\R)),\ \ A_{j+1}\in \SL(2,\R),\ \ f_{j+1}\in C_{h_{j+1}}(\T^d,\sl(2,\R))
$$
such that
$$
\bar{B}_j^{-1}(\theta+\alpha)A_je^{f_j(\theta)}\bar{B}_j(\theta)=A_{j+1}e^{f_{j+1}(\theta)}.
$$
Let $B_{j+1}=B_j(\theta)\bar{B}_j(\theta)$, then we have
$$
B_{j+1}^{-1}(\theta+\alpha)A_0e^{f_0(\theta)}B_{j+1}(\theta)=A_{j+1}e^{f_{j+1}(\theta)},
$$
More precisely, we can distinguish two cases:\\

\noindent \textbf{Non-resonant case:}  If for any $n\in \Z^{d}$ with $0< |n| \leq N_j$, we have
$$
\| 2\xi_j - <n,\alpha> \|_{\R/\Z}\geq \epsilon_j^{\frac{1}{15}},
$$
then by Proposition \ref{iteration}, we have
\begin{equation}\label{con-d}
\|A_{j+1}-A_j\|\leq 2\epsilon_j.
\end{equation}
\begin{align}\label{ge11}
\| \bar{B}_j-id\|_{h_{j+1}}\leq \epsilon_j^{\frac{1}{2}} ,\   \ \| f_{j+1}\|_{h_{j+1}}\leq \epsilon_j^2:=  \epsilon_{j+1}, \ \ 
\end{align}
As a consequence,
 \begin{equation}\label{deg1}
\deg{B_{j+1}}=\deg{B_{j}},\end{equation}
 since $\bar{B}_j(\theta)$ is close to the identity.
By \eqref{estimate2} and \eqref{estimate3}, we have 
\begin{eqnarray*}
\|B_{j+1}\|_{h_{j+1}} &\leq&  C(\alpha)(1+\epsilon_j^{\frac{1}{2}})e^{4\pi h(1-\frac{1}{4^j})|n_{i_j}|}\leq C(\alpha)e^{4\pi h(1-\frac{1}{4^{j+1}})|n_{i_j}|},
\end{eqnarray*}
\begin{eqnarray*}
\|B_{j+1}\|_0 &\leq& C(\alpha)(1+\epsilon_j^{\frac{1}{2}})\prod_{k\leq j}(1+\frac{1}{4^k})|n_{i_j}|^{2\tau}\leq C(\alpha)\prod_{k\leq j+1}(1+\frac{1}{4^k})|n_{i_j}|^{2\tau}, 
\end{eqnarray*}
where we use the simple fact that $\epsilon_j\leq 16^{-j-1}$. 
Moreover, by the selection $ \epsilon_{j+1}=\epsilon_j^2$, we can easily check in this case the truncation satisfies 
\begin{equation}\label{nj}
N_{j+1}=8N_{j}.
\end{equation}
%

\medskip
\noindent \textbf{Resonant case:} If there exists $n_j$  with $0<| n_j| \leq N_j$ such that
$$
\| 2\xi_j- <n_j,\alpha> \|_{\R/\Z}< \epsilon_j^{\frac{1}{15}},
$$
again by Proposition \ref{iteration}, the conjugacy $\bar{B}_j(\theta)$ takes the form as  $\bar{B}_j(\theta)=\bar{B}_j'(\theta)R_{\frac{\langle n_j,\theta\rangle}{2}}$ with estimates
\begin{eqnarray}\label{esti1}
\|\bar{B}_j\|_{h_{j+1}} &\leq& 2\frac{|n_{j}|^{\tau}}{\kappa}e^{2\pi h_{j+1}|n_j|}, \\
\label{esti2}
 \|\bar{B}'_j\|_{h_{j+1}} &<& 2\frac{|n_j|^{\tau}}{\kappa},\\
\label{esti3}
\| f_{j+1}\|_{h_{j+1}}  &\leq& \epsilon_j e^{-h_{j+1}\epsilon_j^{-\frac{1}{18\tau}}} := \epsilon_{j+1}.
\end{eqnarray}
By the construction, 
\begin{equation}\label{deg2}
\deg{B_{j+1}} =\deg{B_{j}}+\deg \bar{B}_j= \deg{B_{j}}  + n_j,
\end{equation}
thus by the definition of $i_{j+1}$, we have $n_{i_{j+1}}=n_j$.  

Furthermore,  the resonance condition implies 
$$
\left\|2\xi_{i_{j+1}}-\langle n_{i_{j+1}},\alpha\rangle\right\|_{\R/\Z}\leq \epsilon_{i_{j+1}}^{\frac{1}{15}},
$$
hence by the assumption that $\alpha\in DC_{d}(\kappa,\tau)$, we have
$$
|\xi_{i_{j+1}}|\geq \frac{\kappa}{2|n_{i_{j+1}}|^{\tau}}-\epsilon_{i_{j+1}}^{\frac{1}{15}}\geq \frac{\kappa}{3|n_{i_{j+1}}|^{\tau}}.
$$
On the other hand, according to Proposition \ref{iteration}, after the $(i_j+1)^{th}$-step, $|\xi_{i_j+1}|\leq \epsilon_{i_j}^{\frac{1}{16}}$, and between $(i_j+1)^{th}$-step and $(i_{j+1}+1)^{th}$-step, there are no resonant steps, it follows from \eqref{con-d} that $|\xi_{i_{j+1}}|\leq 2\epsilon_{i_j}^{\frac{1}{16}}$, consequently we have 
\begin{equation}\label{different resonances}
|n_{i_{j+1}}|\geq \epsilon_{i_j}^{-\frac{1}{17\tau}}|n_{i_j}|.
\end{equation}

 By \eqref{estimate2}, \eqref{esti1} and \eqref{different resonances}, we have
\begin{align*}
\|B_{j+1}\|_{h_{j+1}}&\leq e^{4\pi h(1-\frac{1}{4^j})|n_{i_j}|}\frac{|n_{i_{j+1}}|^{\tau}}{\kappa}e^{2\pi h_{j+1}|n_{i_{j+1}}|}
\leq e^{4\pi h(1-\frac{1}{4^{j+1}})|n_{i_{j+1}}|}.
\end{align*}
By \eqref{estimate3}, \eqref{esti2} and \eqref{different resonances}, we have
\begin{align*}
\|B_{j+1}\|_0&\leq \prod_{k\leq j}(1+\frac{1}{4^k})|n_{i_j}|^{2\tau}\frac{|n_{i_{j+1}}|^{\tau}}{\kappa}
\leq \prod_{k\leq j+1}(1+\frac{1}{4^k})|n_{i_{j+1}}|^{2\tau}.
\end{align*}
We thus finish one step of iteration, and prove the cocycle is almost reducible. \\

\textbf{Step 2: Structure of the conjugacy $B_j(\theta)$:}

Next we prove that the conjugacy $B_j(\theta)$ can be written in the desired form with good estimates. 
Indeed, we may assume that the last resonance happens at some step $0<i_j<j-1$ (if $i_j=j-1$ then the proof would be much simpler). By the above iteration process, there exist $\{\bar{B}_\ell\}_{\ell=i_j}^{j-1}$ and $\{A_\ell\}_{\ell=i_j+1}^{j}$ such that
\begin{equation*}
\bar{B}^{-1}_{\ell}(\theta+\alpha)A_\ell e^{f_\ell(\theta)}\bar{B}_{\ell}(\theta)=A_{\ell+1}e^{f_{\ell+1}(\theta)}.
\end{equation*}

For $\ell=i_j$, since this is a resonant step, one has  $\bar{B}_{i_j}(\theta)=\bar{B}_{i_j}'(\theta)R_{\frac{\langle n_{i_j},\theta\rangle}{2}}$. By \eqref{esti2} and \eqref{esti3}, we have
\begin{eqnarray}\label{esti9}
\| \bar{B}'_{i_j} \|_{h_{i_j+1}} &\leq& 2\frac{|n_{i_{j}}|^{\tau}}{\kappa},\\
\label{esti5}
\| f_{i_j+1}\|_{h_{i_j+1}} &\leq&  e^{-\epsilon_{i_j}^{-\frac{1}{18\tau}}\tilde{h}}=\epsilon_{i_{j+1}}. 
\end{eqnarray}
Furthermore, by Proposition \ref{iteration}, we have
 $$MA_{i_j+1}M^{-1}= exp \left(
\begin{array}{ccc}
 i t^{i_j+1} &  \nu^{i_j+1}\\
\bar{ \nu}^{i_j+1} &  -i t^{i_j+1}
 \end{array}\right)
 $$ with estimate
\begin{equation}\label{esti-nu}|\nu^{i_j+1}|< \epsilon_{i_j}^{\frac{15}{16}}e^{-2\pi|n_{i_j}|\tilde{h}}.\end{equation}

For $i_j+1\leq \ell\leq j-1$, there is no resonance since $i_j$ is the last resonant step before $j$. By the iteration process, we have
\begin{equation}\label{esti6}
\|\bar{B}_\ell-id\|_{h_{\ell+1}}\leq \epsilon_\ell^{\frac{1}{2}}, \ \ \|A_{\ell+1}-A_{\ell}\|\leq 2\epsilon_\ell.
\end{equation}
Now if we denote by
\begin{equation*}
\tilde{B}_j(\theta)=B_{i_j-1}(\theta)\bar{B}'_{i_j}(\theta), \qquad
e^{Y_j(\theta)}=\bar{B}_{i_j+1}(\theta)\cdots \bar{B}_{j-1}(\theta),
\end{equation*}
then by our construction $B_j(\theta)$ can be rewritten as 
\begin{eqnarray*}
B_j(\theta)
 = B_{i_j-1}(\theta) \bar{B}_{i_j}'(\theta)R_{\frac{\langle n_{i_j},\theta\rangle}{2}} \bar{B}_{i_j+1}(\theta)\cdots \bar{B}_{j-1}(\theta)= \tilde{B}_j(\theta)R_{\frac{\langle n_{i_j},\theta\rangle}{2}}e^{Y_j(\theta)}.
\end{eqnarray*}

Note by  \eqref{esti5} and \eqref{esti6},  we have
\begin{equation*}
\|Y_j\|_{\tilde{h}}\leq 2  \|\bar{B}_{i_{j}+1}-id\|_{\tilde{h}}  \leq e^{- \frac{1}{2} \epsilon_{i_j}^{-\frac{1}{18\tau}}\tilde{h}} < e^{-2\pi|n_{i_j}|^4 \tilde{h}},
\end{equation*}
by  \eqref{estimate2},  \eqref{different resonances} and  \eqref{esti9}, one has
\begin{align*}
\|\tilde{B}_j\|_{\tilde{h}} & \leq \|B_{i_j-1}\|_{\tilde{h}}\|B_{i_j}'\|_{\tilde{h}}\leq 2C(\alpha)e^{4\pi h(1-\frac{1}{4^j})|n_{i_{j-1}}|}\frac{|n_{i_{j}}|^{\tau}}{\kappa} \leq C(\alpha) |n_{i_{j}}|^{\tau}e^{2\epsilon_{i_{j-1}}^{\frac{1}{18\tau}}|n_{i_j}|},
\end{align*}
which  proves \eqref{es3}, while  \eqref{es4} follows directly from  \eqref{estimate3}.   By  \eqref{resos} and \eqref{different resonances}, one can estimate
$$|\deg B_j - n_{i_j}| \leq \sum\limits_{m=1}^{j-1}|n_{i_m}|\leq 2\epsilon_{i_{j-1}}^{\frac{1}{18\tau}}|n_{i_j}|,$$
which gives \eqref{es5}.  These finish the estimates of the conjugacy $B_j$.\\

\textbf{Step 3: Quantitative estimates of $\nu_j$ and $f_j$.}

By \eqref{esti6}, we have
$
 \|A_j-A_{i_j+1}\|\leq  C\epsilon_{i_j+1}^{\frac{1}{2}},
$
combining \eqref{esti-nu}, which imply that 
\begin{eqnarray*}
|\nu^j| &\leq& |\nu^{i_j+1}| +\|A_j-A_{i_j+1}\| \leq \epsilon_{i_j}^{\frac{15}{16}} e^{-2\pi|n_{i_j}|\tilde{h}}+C\epsilon_{i_j+1}^{\frac{1}{2}}<2\epsilon_{i_j}^{\frac{15}{16}}e^{-2\pi|n_{i_j}|\tilde{h}}.
\end{eqnarray*}
Finally,  if we take $j=i_{j+1}$ to be the resonant step, while $(j-1)$-th is non-resonant, then according to Proposition \ref{iteration}, while $f_j$ is well defined in  $ C_{h_{j}}(\T^d,\sl(2,\R))$, we have estimates 
$$ \|f_{j}\|_{ \tilde{h}}\leq \epsilon_{j-1} e^{-2\pi N_{j-1}(h_j-\tilde{h})},$$
by \eqref{nj}, we have $N_{i_{j+1}}=8N_{i_{j+1}-1},$ it follows that
$$ \|f_{j}\|_{ \tilde{h}}\leq \epsilon_{i_j} e^{- \frac{\pi}{8} N_{i_{j+1}}(h-\tilde{h})}  \leq   \epsilon_{i_j} e^{- \frac{\pi}{8} |n_{i_{j+1}}|(h-\tilde{h})} .$$

\end{pf}

\subsection{KAM resonances and rotation number resonances}\label{kam-res}

In this section, we delve into the intricate connection between KAM resonances and rotation number resonances.

For any $\e_0>0$, we say that $k$ is an $\e_0$-resonance if 
\begin{equation}\label{upper+}
\|2\rho-\langle k,\alpha \rangle\|_{\R/\Z}\leq e^{-|k|\e_0}
\end{equation} 
and 
\begin{equation}\label{upper++}
 \|2\rho-\langle k,\alpha\rangle\|_{\R/\Z}=\min_{|j|\leq |k|}\|2\rho-\langle j,\alpha\rangle\|_{\R/\Z}.   
\end{equation}
As pointed out by Avila-Jitomirskaya \cite{aj1}, if $\alpha \in DC_d(\kappa,\tau)$, then \eqref{upper+} implies \eqref{upper++} for $|n|>n(\gamma,\tau)$.
We order the $\e_0$-resonances $0=|\ell_0|<|\ell_1|\leq |\ell_2|\leq \cdots$. We say that $\rho$ is $\e_0$-resonant if the set of resonances is infinite. If $\rho$ is non-resonant, with the set of resonances $\{\ell_0,\cdots, \ell_j\}$ we formally set $\ell_{j+1}\in \Z^d$ be any vector such that $|\ell_{j+1}|=\infty$. 
 \begin{Lemma} \label{relation1}
Under the assumption of Proposition  \ref{reducibility}, for any $\e_0>0$, there exist $B_j\in C_{\tilde{h}}^\omega(\T^d, \PSL(2,\R))$,   such that 
\begin{equation}\label{conjugacy}
B_j^{-1}(\theta+\alpha)A(\theta)B_j(\theta)=A_je^{f_j(\theta)}=M^{-1}exp \left(
\begin{array}{ccc}
 i t^j &  \nu^j\\
\bar{\nu}^j &  -i t^j
 \end{array}\right)Me^{f_j(\theta)}
\end{equation}
with $\ell_i=\deg B_j$ provided $\ell_i$ is large enough (the choice of $j$ depends on $\ell_i$). 
\end{Lemma}
\begin{pf}
First by Proposition \ref{reducibility},  $(\alpha,A)$ is almost reducible, and 
there exist $B_j\in C_{\tilde{h}}^\omega(\T^d,$ $\PSL(2,\R))$,   such that 
$$
B_j^{-1}(\theta+\alpha)A(\theta)B_j(\theta)=A_je^{f_j(\theta)}$$
with estimate $\|f_j\|_{\tilde{h}} \leq \epsilon_j,$ where 
\begin{equation}\label{eqle2}
\epsilon _{k+1}=
\begin{cases}
\epsilon_{k}^2&\text{$k\neq i_m$},\\
\epsilon_{k}e^{-h_{k}\epsilon_{k}^{-\frac{1}{18\tau}}}&\text{$k=i_m$}.
\end{cases}
\end{equation}
and $n_1<n_2 <\cdots < n_{i
_j} \leq j-1$ are the resonances, see \eqref{ge11} and \eqref{esti3}. 

First  by \eqref{rotation number}, we have 
$$  2\rho(\alpha, A_ke^{f_k})=  2\rho(\alpha,A)- \langle \deg B_k, \alpha \rangle,$$
by Lemma \ref{rol}, we have 
\begin{equation*}
|  \rho(\alpha, A_ke^{f_k})- \xi_k|= |\rho(\alpha,A_{k}e^{f_{k}})-\rho(\alpha,A_{k})|\leq C\epsilon^\frac{1}{2}_{k},    
\end{equation*}
it then  follows from   \eqref{upper+}  that
\begin{equation}\label{diff-rot}
\|\langle \ell_i-\deg B_k, \alpha \rangle-   2\xi_k\|_{\R/\Z}\leq  2C \epsilon^\frac{1}{2}_{k}+ e^{- |\ell_i|\e_0}
\end{equation}
Now for any $\ell_{i}\in\{\ell_i\}_{i=1}^\infty$, there exists $j$ such that $N_{j-1}\leq |\ell_i|< N_{j}$, thus we can distinguish the proof into two cases:\\

\smallskip
\textbf{Case I:  $j-1=i_j$.} By the definition of resonances,  there exists $|n_{i_j}|\leq N_{j-1}$  we have
\begin{equation*}
\|2\xi_{j-1}-\langle n_{i_j},\alpha\rangle\|_{\R/\Z}\leq \epsilon_{j-1}^{\frac{1}{15}},
\end{equation*}
combining \eqref{diff-rot} with $k=j-1$, we have 
\begin{equation*}
\| \langle \ell_i -\deg  B_{j}, \alpha \rangle \|_{\R/\Z}= \| \langle \ell_i -\deg  B_{j-1}-n_{i_j}, \alpha \rangle \|_{\R/\Z}\leq 2 \epsilon_{j-1}^{\frac{1}{15}}.
\end{equation*}
By Diophantine condition $\alpha\in DC_d(\kappa,\tau)$, if $\ell_i\neq \deg B_j$ and  if $j$ is large, we have
$$
|\ell_i-\deg B_j|\geq \kappa^{\frac{1}{\tau}} \epsilon_{j-1}^{-\frac{1}{15\tau}}\gg N_{j}=\frac{2|\ln\epsilon_j|}{h_j-h_{j+1}},
$$
the last inequality holds by \eqref{eqle2}. 

However, by  \eqref{es5},  we have 
\begin{equation}\label{lb} |\ell_i-\deg B_j|\leq N_j + 2|n_{i_j}| \leq N_j +2 N_{j-1} \leq 2  N_j ,\end{equation}
this is a  contradiction. Thus in this case $\ell_i=\deg B_j$ for $\ell_i$ sufficiently large.\\

 \textbf{Case II:  $j-1\neq i_j$.} In this case, by \eqref{nj}, we have $N_{j}=8N_{j-1}$. Apply \eqref{diff-rot} with $k=j$, we have 
\begin{equation*}
\|\langle \ell_i-\deg B_j, \alpha \rangle- 2\xi_j\|_{\R/\Z}\leq  2C \epsilon^\frac{1}{2}_{j}+ e^{-  \frac{|N_j|}{8}\e_0} \leq  \epsilon_{j}^{\frac{1}{15}}.
\end{equation*}
If  $\ell_i-\deg B_j=0$, there is nothing to prove. Otherwise,  similar as \eqref{lb},  $|\ell_i-\deg B_j| \leq  2  N_j.$  
We further distinguish the proof into two subcases:\\

\noindent \textbf{Case II-1:} If  $j=i_{j+1}.$   Then by the definition of resonances, there exists $|n_{i_{j+1}}|\leq N_{j}$,  such that 
\begin{equation*}
\|2\xi_{j}-\langle n_{i_{j+1}},\alpha\rangle\|_{\R/\Z}\leq \epsilon_{j}^{\frac{1}{15}},
\end{equation*}
then same as Case I, by the Diophantine condition $\alpha\in DC_d(\kappa,\tau)$, if $\ell_i-\deg B_j \neq n_{i_{j+1}}$, then  
$$
|\ell_i-\deg B_j- n_{i_{j+1}} |\geq \kappa^{\frac{1}{\tau}} \epsilon_{j}^{-\frac{1}{15\tau}}\gg N_{j}
$$
which again contradicts to $|\ell_i-\deg B_j- n_{i_{j+1}} |< 3 N_j$.\\

\noindent \textbf{Case II-2:} If $j\neq i_{j+1}$, i.e. it is still non-resonant, then again \eqref{nj}, we have $N_{j+1}=8N_{j}$, 
Apply \eqref{diff-rot} with $k=j+1$, we have 
\begin{equation*}
\|\langle \ell_i-\deg B_{j+1}, \alpha \rangle- 2\xi_{j+1}\|_{\R/\Z}\leq  2C \epsilon^\frac{1}{2}_{j+1}+ e^{-  \frac{|N_{j+1}|}{64}\e_0} \leq  \epsilon_{j+1}^{\frac{1}{15}}.
\end{equation*}
Since  $|\ell_i-\deg B_{j+1}|= |\ell_i-\deg B_{j}|< 2  N_j < N_{j+1}$, by the definition of resonances, this forces  $j+1 =i_{j+1}$,  and consequently there exists 
$|n_{i_{j+1}}|< N_{j+1}$, such that  $\ell_i-\deg B_{j+1}=n_{i_{j+1}}$, which just means $\ell_i=\deg B_{j+2}$.
\end{pf}

\section{Optimal arithmetic reducibility and Anderson localization}\label{sec4}

In this section, we demonstrate the practical applications of structured quantitative almost reducibility when the condition $2\pi h > \delta(\alpha, \rho)$ holds. By doing so, we extend Eliasson's reducibility theorem \cite{eli} in an optimal manner. This advancement allows us to exploit the power of Aubry duality to achieve the optimal arithmetic Anderson localization.

\subsection{Optimal arithmetic reducibility-local case}

We will construct a novel KAM scheme for $SL(2, \mathbb{R})$ quasiperiodic cocycles with Liouville rotation numbers, which constitutes a generalization of Eliasson's theorem. To appreciate the challenge, consider Proposition \ref{iteration}, which highlights two key obstacles: the normalization transformation $P_j$, of order $\mathcal{O}(|n_{i_j}|^\tau)$, and the rotation transformation $R_{\frac{\langle n_{i_j},\theta\rangle}{2}}$, with a size of $\mathcal{O}(e^{\pi h|n_{i_j}|})$. When dealing with an infinite number of resonances, these transformations can potentially lead to a divergent conjugacy. However, the strength of quantitative structured almost reducibility, as encapsulated in Proposition \ref{reducibility}, enables us to overcome these hurdles. We wish to clarify that our strategy  is to offset the effect of each large conjugacy
by immediately applying the inverse of the conjugacy after a KAM
step, essentially "undoing" the conjugacy. This technique, which was first developed  on the
linearization of quasiperiodically forced circle flows  \cite{KJYZ}.

\begin{Theorem}\label{reducibilitylocal}
Let $\alpha\in DC_d(\kappa,\tau)$,  $A\in C_{h}^\omega(\T^d,SL(2,\R))$ with $ 2\pi h>2\pi \tilde{h}>\delta>0$, $R\in SL(2,\R)$. If\begin{itemize}
\item $$\|A(\cdot)-R\|_{h}\leq \bar{\epsilon}:=\frac{D_0(\kappa,\tau,d)}{\|R\|^{C_0}}(\frac{h-\tilde{h}}{2})^{C_0\tau},$$
\item $\rho(\alpha,A)\in LC_\alpha(\gamma,\delta)$, where
$$
LC_\alpha(\gamma,\delta)=\{\phi\in\R|\|2\phi-\langle m,\alpha\rangle\|_{\R/\Z}\geq \gamma e^{- |m|\delta},\ \ m\in\Z^d\}.
$$
\end{itemize}
then $(\alpha,A)$ is analytically reducible. Moreover, there exist $C=C(\gamma,\delta,\alpha,\|R\|)$ and $B_A\in C^\omega_{\tilde{h}-\frac{\delta}{2\pi}}(\T^d,\SL(2,\R))$ with $\deg B_A=0$, such that
$$
B_A(\theta+\alpha)^{-1}A(\theta)B_A(\theta)=R_{\rho(\alpha,A)},
$$
with estimate $\|B_A\|_{\tilde{h}-\frac{\delta}{2\pi}}\leq C$, and $B_A$ is continuous on $\mathcal{A}_\alpha(\gamma,\delta)$ in the sense
$$
\lim\limits_{A',A\in \mathcal{A}_\alpha(\gamma,\delta)\atop \|A'-A\|\rightarrow 0}\|B_{A'}-B_A\|_{\tilde{h}-\frac{\delta}{2\pi}}=0
$$
where
$$
\mathcal{A}_\alpha(\gamma,\delta)=\{A:\rho(\alpha,A)\in L_\alpha(\gamma,\delta)\}.
$$
\end{Theorem}
\begin{Remark}
If one replace $ LC_\alpha(\gamma,\delta)$ by $DC_\alpha(\gamma,\tau)$ where
$$
DC_\alpha(\gamma,\tau)=\{\phi\in\R|\|2\phi-\langle m,\alpha\rangle\|_{\R/\Z}\geq \frac{\gamma}{(|m|+1)^\tau},\ \ m\in\Z^d\}.
$$
Then, it is the classical Eliasson's reducibility theorem \cite{eli}, the main difference here is that one replaces the polynomial small lower bound by the exponential small lower  bound which will lead to essential difficulties in the proof. 
\end{Remark}

%

\begin{pf} 
We will first prove the cocycle is analytically reducible, then give the continuity argument. 

\subsubsection{Part 1: Reducibility}

First by Proposition  \ref{reducibility}, there exist  $Z_p\in C_{\frac{h+\tilde{h}}{2}}^\omega(\T^d, PSL(2,\R))$ and $i_p\in \N$ such that 
$$
Z_p^{-1}(\theta+\alpha)A(\theta)Z_p(\theta)=\widetilde{A}_pe^{\widetilde{f}_p(\theta)}= M^{-1}exp \left(
\begin{array}{ccc}
 i \widetilde{t}^p &  \widetilde{\nu}^p\\
\bar{ \widetilde{\nu}}^p &  -i \widetilde{t}^p
 \end{array}\right)Me^{\widetilde{f}_p(\theta)}
$$
 with estimates
 \begin{eqnarray}
   \nonumber \| \widetilde{f_p}\|_{\frac{h+\tilde{h}}{2}} &\leq&  \tilde{\epsilon}_{i_p} e^{- \frac{\pi}{16} |n_{i_{p+1}}|(h-\tilde{h})},\\
\nonumber   |\widetilde{\nu}^p| &\leq& 2  \tilde{\epsilon}_{i_p} ^{\frac{15}{16}} e^{-\pi|n_{i_p}|(h+\tilde{h})},\\
\label{es5-11} 
 |\ell_p|:= |\deg Z_p| &\leq& (1+2 \tilde{\epsilon}_{i_{p-1}} ^{\frac{1}{18\tau}} )|n_{i_p}|
 \end{eqnarray}
 where $\tilde{\epsilon}_p\leq \bar{\epsilon}^{2^p}$.

If there are finite resonant steps, i.e. there are at most finitely many $n_i$, Propsition  \ref{reducibility} actually implies that the cocycle is reducible (actually in this case $|n_{i_{j+1}}|=\infty$), then there is nothing to say. In the following, we always assume that there are infinite many resonant steps, i.e $|n_i| \rightarrow \infty$.
Just note that Proposition \ref{reducibility} implies that the constant matrix $\widetilde{A}_p$ are uniformly bounded,  one can thus select 
$p_0\in\Z$ (independent of $\widetilde{A}_p$) to be the smallest integer $p$ (or equivalently to say, to select the smallest $i_p$) satisfying
\begin{eqnarray}\label{reducibility condition}
\label{reducibility condition-31}  \tilde{\epsilon}_{i_p}  e^{- \frac{\pi}{16} |n_{i_{p+1}}|(h-\tilde{h})} &\leq&   \frac{D_0\gamma^{4} e^{-4 |n_{i_p}|(1+2 \tilde{\epsilon}_{i_{p-1}} ^{\frac{1}{18 \tau }} )\delta} }  {\|  \widetilde{A}_p \|^{C_0}}(\frac{h-\tilde{h}}{2})^{C_0\tau},\\
\label{reducibility condition-2}  \tilde{\epsilon}_{i_p} &\leq& \gamma^{4} ,\\
 \label{reducibility condition-1}
\pi (h+\tilde{h})&>& (1+2 \tilde{\epsilon}_{i_{p-1}} ^{\frac{1}{18\tau}} ) \delta,
\end{eqnarray}
where \eqref{reducibility condition-31} follows from  \eqref{different resonances}, \eqref{reducibility condition-1}  follows from our selection that  $2\pi h>2\pi \tilde{h}>\delta$.

%
%

Once we have this, we can prove the result by induction. Take $(\alpha,A_{0}e^{f_{0}}):=(\alpha,\widetilde{A}_{p_0}e^{\widetilde{f}_{p_0}})$,
 $$h_0= \frac{h+\tilde{h}}{2}, \quad  \epsilon_0 =  \tilde{\epsilon}_{i_p}  e^{- \frac{\pi}{16} |n_{i_{p+1}}|(h-\tilde{h})}.$$
Assume that we are at the $(j+1)^{th}$ KAM step,  i.e. we already construct $B_j\in C^\omega_{h_{j}}(\T^d,SL(2,\R))$ such that
\begin{equation}\label{estima18}
B_{j}^{-1}(\theta+\alpha)A_0e^{f_0(\theta)}B_{j}(\theta)=A_{j}e^{f_{j}(\theta)}=M^{-1}exp \left(
\begin{array}{ccc}
 i t^j &  \nu^j\\
\bar{\nu}^j &  -i t^j
 \end{array}\right)Me^{f_{j}(\theta)},
\end{equation}
where $A_j\in SL(2,\R)$ has two eigenvalues $e^{\pm i\xi_j}$ and
\begin{equation}\label{estima2-1}
 \|f_j\|_{h_j} \leq \epsilon_j\leq \epsilon_0^{2^j}.
\end{equation}
 Moreover,  if we denote $j_1<j_2 <\cdots < j_{\ell} <\cdots$ to be all the KAM resonant steps,  we have the following estimates:
 \begin{eqnarray}\label{estima1}
\|B_{j}\|_{h_j- \frac{\delta}{2\pi}} \leq   (1+\tilde{\epsilon}_{i_p}^{\frac{1}{4}} ) \prod_{ \ell \atop j_{\ell}\leq j}(1+\epsilon_{j_\ell}^{\frac{1}{4}} )\prod_{k\leq j}(1+\frac{1}{4^k}),
 \end{eqnarray}
 \begin{equation}\label{add1}
 \|B_j\|_{h_j}\leq \epsilon_{j}^{-\frac{1}{100}},
 \end{equation}
 \begin{eqnarray}\label{estima3-1}
 |\nu^j| \leq \left\{ \begin{array}{ccc}2 \tilde{\epsilon}_{i_p} ^{\frac{15}{16}}e^{-\pi|n_{i_p}|(h+\tilde{h})} + \sum_{k=0}^{j-1} \epsilon_k & j< j_1\\
0 & j=j_i \\
\sum_{k={j_{i-1}}}^{j-1} \epsilon_k &  j_{i-1} < j < j_i.
\end{array} \right.
\end{eqnarray}

We define the induction sequence
$$ 
h_j-h_{j+1}=\frac{h_0-\tilde{h}}{4^{j+1}},\ \ N_j=\frac{2|\ln\epsilon_j|}{h_j-h_{j+1}}.
$$
By our selection of  $\epsilon_0$, one can check that
\begin{equation}\label{estima2}
\epsilon_j \leq \frac{D_0 \gamma^4 e^{-4 |n_{i_p}|(1+2 \tilde{\epsilon}_{i_p} ^{\frac{1}{18 \tau }} )\delta} }{\|A_j\|^{C_0}}(h_j-h_{j+1})^{C_0\tau}.
\end{equation}
 Indeed, $\epsilon_j$ on the left side of the inequality decays at least super-exponentially with $j$, while $(h_j-h_{j+1})^{C_0\tau}$ on the right side decays exponentially with $j$.

Note that \eqref{estima2} implies that one can apply  Proposition \ref{iteration}, consequently 
there exists $\bar{B}_{j}\in C^\omega_{h_{j+1}}(\T^d,\PSL(2,\R))$, $\bar{A}_{j+1}\in \SL(2,\R)$ and $\bar{f}_{j+1}\in C_{h_{j+1}}^\omega(\T^d,\sl(2,\R))$ such that
\begin{equation}\label{estima9}
\bar{B}_{j}^{-1}(\theta+\alpha)A_je^{f_j(\theta)}\bar{B}_{j}(\theta)=\bar{A}_{j+1}e^{\bar{f}_{j+1}(\theta)}=M^{-1}exp \left(
\begin{array}{ccc}
 i \overline{t^{j+1}} & \overline{ \nu^{j+1}}\\
\overline{\bar{ \nu}^{j+1}} &  -i \overline{ t^{j+1}}
 \end{array}\right)M e^{\bar{f}_{j+1}(\theta)}
\end{equation}
Moreover, we can distinguish the following two cases:\\
\smallskip

\noindent \textbf{Non-resonant case:}  If the $(j+1)^{th}$ is non-resonant, i.e.  for any $n\in \Z^{d}$ with $0< |n| \leq N_j$, we have
$$
\| 2\xi_j - <n,\alpha> \|_{\R/\Z}\geq \epsilon_j^{\frac{1}{15}},
$$
by Proposition \ref{iteration}, the conjugacy $\bar{B}_j$ takes the form  $\bar{B}_j(\theta)=e^{\bar{Y}_j(\theta)}$ with 
\begin{equation}\label{estima31}
\|\bar{Y}_j\|_{h_{j+1}}< \epsilon_j^{\frac{1}{2}},\qquad \|\bar{A}_{j+1}-A_j\|\leq  2\epsilon_j.
\end{equation}
$$\| \bar{f}_{j+1}\|_{h_{j+1}}\leq \epsilon_j^2:=  \epsilon_{j+1}.$$

Let $B_{j+1}=B_j(\theta)\bar{B}_j(\theta)$, we have
$$
B_{j+1}^{-1}(\theta+\alpha)A_0e^{f_0(\theta)}B_{j+1}(\theta)=\bar{A}_{j+1}e^{\bar{f}_{j+1}(\theta)}:=A_{j+1}e^{f_{j+1}(\theta)}.
$$
By \eqref{estima31},  $B_{j+1}$  satisfy the estimate \eqref{estima1} and \eqref{add1} with $j+1$ in place of $j$, and 
 $ |\nu^{j+1}|$ satisfy the estimate \eqref{estima3-1} with $j+1$ in place of $j$.\\

%

\smallskip
\noindent \textbf{Resonant case:} If  the $(j+1)^{th}$ is resonant, i.e. there exists $n_j$ with $0<| n_j| \leq N_j$ such that
$$
\| 2\xi_j- <n_j,\alpha> \|_{\R/\Z}< \epsilon_j^{\frac{1}{15}},
$$
by Proposition \ref{iteration}, the conjugacy $\bar{B}_j$ takes the form  $$\bar{B}_j(\theta)=\bar{P}_je^{\bar{Y}_j(\theta)}R_{\frac{\langle n_j,\theta\rangle}{2}}$$
where $P_j$ is given by
\begin{equation}\label{pjes} \bar{P}_j^{-1} A_j \bar{P}_j= R_{\xi_j}.\end{equation}
Moreover, we have the following estimates:
\begin{eqnarray}\label{estima7}
\|\bar{Y}_j\|_{h_{j+1}}&<&\epsilon_j^{\frac{1}{2}},\\
\label{estima8}
\| \bar{f}_{j+1}\|_{h_{j+1}}  &\leq& \epsilon_j e^{-h_{j+1}\epsilon_j^{-\frac{1}{18\tau}}}: =\bar{\epsilon}_{j+1},\\
\label{est2}
|\overline{\nu^{j+1}}|&<&\epsilon_j^{\frac{15}{16}} e^{-2\pi|n_j|h_j}.
\end{eqnarray}

 To make the whole scheme converges, we need to add two new ingredients 
to the KAM machine. First we need a refined estimate of $\bar{P}_j$, which will be taken care in the next step.\\

\smallskip
\textbf{Step 1: Refined estimate of normalization  of $A_j$.}

\begin{Lemma}\label{liouass}
Assume that  $\rho(\alpha,A)\in LC_\alpha(\gamma,\delta)$, then 
\begin{align}\label{rotation}
\rho(\alpha,A_{0}e^{f_{0}})\in LC_\alpha(\gamma e^{-  |n_{i_p}|(1+2 \tilde{\epsilon}_{i_{p-1}} ^{\frac{1}{18 \tau }} )\delta},\delta).
\end{align}
As a consequence, we have 
\begin{align}\label{estima41}
\|\bar{P}_j-id\| \leq \left\{ \begin{array}{ccc} \tilde{\epsilon}_{i_p} ^{\frac{1}{2}} & j +1= j_1\\
 \epsilon_{j_{i-1}} ^{\frac{1}{2}} &   j+1 =  j_i >j_1.
\end{array} \right.\end{align}
\end{Lemma}
\begin{pf}

 Note by \eqref{rotation number}, we have 
\begin{equation}\label{new123}
2\rho(\alpha,A_{0}e^{f_{0}})=2\rho(\alpha,A)-\langle\deg{Z_{p}},\alpha\rangle=2\rho(\alpha,A)-\langle\ell_p,\alpha\rangle, 
\end{equation}
thus   if $\rho(\alpha,A)\in LC_{\alpha}(\gamma,\delta)$,   then \eqref{new123}  implies that
\begin{eqnarray*}
\|2\rho(\alpha,A)-\langle m,\alpha\rangle-\langle \ell_{p},\alpha\rangle\|_{\R/\Z} \geq \gamma e^{-(|m+\ell_{p}|) \delta}\geq   \gamma e^{- |n_{i_p}|(1+2 \tilde{\epsilon}_{i_{p-1}} ^{\frac{1}{18 \tau }} )\delta}e^{-\delta|m|},
\end{eqnarray*}
which directly implies \eqref{rotation}.

Similarly,  notice that $\deg{B_{j}}=0$ by \eqref{estima1}, then by \eqref{rotation number}, we have
 \begin{equation}\label{rinv}
\rho(\alpha,A_je^{f_j})=\rho(\alpha,A_0e^{f_0}),
\end{equation}
 consequently by   \eqref{reducibility condition} and \eqref{rotation}, we have
\begin{equation}\label{estima30}
|\rho(\alpha,A_je^{f_j})|\geq \gamma e^{- |n_{i_p}|(1+2 \tilde{\epsilon}_{i_{p-1}} ^{\frac{1}{18 \tau }} )\delta} \geq \epsilon_0^{\frac{1}{4}}.
\end{equation}
By Lemma \ref{rol}, we have
\begin{equation}\label{estima32}
|\rho(\alpha,A_j)-\rho(\alpha,A_je^{f_j})|\leq C\epsilon_j^{\frac{1}{2}},
\end{equation}
which implies that
\begin{equation}\label{estima33}
|\xi_j|=|\rho(\alpha,A_j)|\geq \gamma e^{- |n_{i_p}|(1+2 \tilde{\epsilon}_{i_{p-1}} ^{\frac{1}{18 \tau }} )\delta} -C\epsilon_j^{\frac{1}{2}}\geq \frac{\gamma}{2} e^{- |n_{i_{p}}|(1+2 \tilde{\epsilon}_{i_{p-1}} ^{\frac{1}{18 \tau }} )\delta}.
\end{equation}

On the other hand, by \eqref{estima3-1}, one can estimate 
 \begin{eqnarray*}
 |\nu^j| \leq \left\{ \begin{array}{ccc}4 \tilde{\epsilon}_{i_{p}} ^{\frac{15}{16}}e^{-\pi|n_{i_p}|(h+\tilde{h})}  & j< j_1\\
2 \epsilon_{ j_{i-1}} &     j_{i-1} \leq j < j_i
\end{array} \right.
\end{eqnarray*}
Therefore, if $j +1= j_1$, by  \eqref{reducibility condition-2} and \eqref{reducibility condition-1}, we have
$$
\left|\frac{2\nu^{j}}{\xi_{j}}\right|\leq \frac{16 \tilde{\epsilon}_{i_p} ^{\frac{15}{16}}e^{-\pi|n_{i_p}|(h+\tilde{h})} }{\gamma e^{- |n_{i_p}|(1+2 \tilde{\epsilon}_{i_{p-1}} ^{\frac{1}{18 \tau }} )\delta}}\leq  \tilde{\epsilon}_{i_p}^{\frac{1}{2}} <1.
$$
If $ j+1 =  j_i >j_1$, by  \eqref{estima30} and \eqref{estima33}, we have
$$
\left|\frac{2\nu^{j}}{\xi_{j}}\right|\leq \frac{4\epsilon_{ j_{i-1}}}{\epsilon_0^{\frac{1}{4}}}\leq \epsilon^{\frac{1}{2}}_{ j_{i-1}} < 1.
$$
By Lemma \ref{Dominate matrix}, $\bar{P}_j$ satisfy \eqref{pjes} with estimate \eqref{estima41}.
\end{pf}


\smallskip

\textbf{Step 2: Normalization of $\bar{A}_{j+1}$.} The rest difficulty  is to deal with the effect of rotation transformation $R_{\frac{\langle n_{j},\theta\rangle}{2}}$. The key observation here is that if the rotation number $\rho(\alpha,A)\in LC_\alpha(\gamma,\delta)$, one can  normalize $\bar{A}_{j+1}$ with a conjugacy $\bar{U}_{j+1}$ which is close to the identity, consequently one can rotation backward $R_{\frac{\langle n_{j},\theta\rangle}{2}}$, and the conjugacy $R_{\frac{\langle n_j,\theta\rangle}{2}}\bar{U}_{j+1}R_{-\frac{\langle n_j,\theta\rangle}{2}}$ is close to the identity.

Denote $spec(\bar{A}_{j+1})=\{e^{ i \overline{\xi_{j+1}}},e^{- i \overline{\xi_{j+1}}} \}$, then one can normalize  $\bar{A}_{j+1}$ as follows:
\begin{Lemma}
Assume that  $\rho(\alpha,A)\in LC_\alpha(\gamma,\delta)$, then 
\begin{align*}
| \overline{\xi_{j+1}}|\geq \frac{1}{8}\epsilon_j^{\frac{1}{8}}e^{- \delta|n_j|}.
\end{align*}
Consequently, there exists $\bar{U}_{j+1} \in SL(2,\R)$ such that
\begin{equation}\label{estima17}
\bar{U}_{j+1}^{-1}\bar{A}_{j+1}\bar{U}_{j+1}=R_{\overline{\xi_{j+1}}}
\end{equation}
with estimate
\begin{equation}\label{estima16}
\|\bar{U}_{j+1}-id\|\leq \epsilon_j^{\frac{1}{2}}e^{-(2\pi h_j-\delta)|n_j|}.
\end{equation}
\end{Lemma}

\begin{pf}
Notice that $\deg{\bar{B}_{j}}=n_{j}$, then by \eqref{rotation number}, \eqref{estima9} and \eqref{rinv}, we have
\begin{equation}\label{n1}
2\rho(\alpha,\bar{A}_{j+1}e^{\bar{f}_{j+1}})=2\rho(\alpha,A_je^{f_j})-\langle n_j,\alpha\rangle= 2\rho(\alpha,A_0e^{f_0})-\langle n_j,\alpha\rangle,
\end{equation}
by \eqref{estima8} and Lemma \ref{rol}, we have
\begin{equation*}
|\rho(\alpha,\bar{A}_{j+1}e^{\bar{f}_{j+1}})-\rho(\alpha,\bar{A}_{j+1})|\leq C \bar{\epsilon}^\frac{1}{2}_{j+1},
\end{equation*}
meanwhile by \eqref{estima9} and \eqref{est2}, we have
\begin{equation*}
| \overline{t^{j+1}}-\rho(\alpha,\bar{A}_{j+1})|\leq 2| \overline{\nu^{j+1}}|\leq 2\epsilon_j^{\frac{15}{16}} e^{-2\pi|n_j|h_j}.
\end{equation*}
It thus follows from \eqref{n1} that 
\begin{equation}\label{n2}
|2 \overline{t^{j+1}}-2\rho(\alpha,A_0e^{f_0})+\langle n_j,\alpha\rangle|\leq 4\epsilon_j^{\frac{15}{16}} e^{-2\pi|n_j|h_j}+2 C\bar{\epsilon}_{j+1}^{\frac{1}{2}}.
\end{equation}

On the other hand, by Lemma \ref{liouass}  and  our selection of $\epsilon_0$, 
we have 
\begin{equation}\label{estima13}
\|2\rho(\alpha,A_0e^{f_0})-\langle n_j,\alpha\rangle\|_{\R/\Z}\geq \epsilon_0^{\frac{1}{4}}e^{- \delta|n_j|},
\end{equation}
thus \eqref{estima13} and \eqref{n2} imply that
\begin{align}\label{estima14}
| \overline{t^{j+1}}|&\geq \frac{1}{2}\epsilon_0^{\frac{1}{4}}e^{- \delta|n_j|}-2\epsilon_j^{\frac{15}{16}} e^{-2\pi|n_j|h_j}- 2C\bar{\epsilon}_{j+1}^{\frac{1}{2}}\geq \frac{1}{4}\epsilon_j^{\frac{1}{4}}e^{-\delta|n_j|},
\end{align}
where in the final inequality, we use the fact that $ 2\pi h_{j}>2 \pi \tilde{h}>\delta$ and \eqref{estima8}.

Notice that $spec(\bar{A}_{j+1})=\{e^{ i \overline{\xi_{j+1}}},e^{-i  \overline{\xi_{j+1}}} \}$.  By \eqref{est2} and \eqref{estima14},  we have
\begin{equation*}
\begin{split}
| \overline{\xi_{j+1}}|^2 &=| \overline{t^{j+1}}|^2-| \overline{\nu^{j+1}}|^2 \geq \frac{1}{16}\epsilon_j^{\frac{1}{2}}e^{-2\delta|n_j|}-\epsilon_j^{\frac{5}{8}} e^{-4\pi h_j|n_j|} \geq  \frac{1}{64}\epsilon_j^{\frac{1}{2}}e^{-2\delta|n_j|}
\end{split}
\end{equation*}
which implies that 
$$
\left|\frac{2 \overline{\nu^{j+1}}}{ \overline{\xi_{j+1}}}\right|\leq \frac{2\epsilon_j^{\frac{15}{16}} e^{-2\pi h_j|n_j|}}{\frac{1}{8}\epsilon_j^{\frac{1}{4}}e^{-\delta|n_j|}}\leq 16\epsilon_j^{\frac{11}{16}}e^{-(2\pi h_j-\delta)|n_j|}\leq 1,
$$
then the result follows from Lemma \ref{Dominate matrix} directly. \end{pf}

\textbf{Step 3: Rotation backward.} Now we can make a rotation backward,
which makes the conjugacy close to identity in a shrinking strip. More precisely,  let
\begin{equation}\label{estima19}
\tilde{B}_{j}(\theta)=\bar{B}_j(\theta)\bar{U}_{j+1}R_{-\frac{\langle n_j,\theta\rangle}{2}}=\bar{P}_{j}e^{\bar{Y}_{j}(\theta)}R_{\frac{\langle n_j,\theta\rangle}{2}}\bar{U}_{j+1}R_{-\frac{\langle n_j,\theta\rangle}{2}},
\end{equation}
then by \eqref{estima9} and \eqref{estima17}, we have
\begin{equation*}
\tilde{B}_{j}^{-1}(\theta+\alpha)A_je^{f_j(\theta)}\tilde{B}_{j}(\theta) = R_{\overline{\xi_{j+1}}+\langle n_j,\alpha\rangle}e^{f_{j+1}(\theta)} := A_{j+1}e^{f_{j+1}(\theta)},
\end{equation*}
where 
\begin{equation} \label{newf}f_{j+1}(\theta)= R_{\frac{\langle n_j,\theta\rangle}{2}} \bar{U}_{j+1}^{-1}  \bar{f}_{j+1}(\theta)  \bar{U}_{j+1}  R_{-\frac{\langle n_j,\theta\rangle}{2}},\end{equation}
consequently let $B_{j+1}=B_j(\theta)\tilde{B}_j(\theta)$,  then we have
\begin{equation}\label{almost1}
B_{j+1}^{-1}(\theta+\alpha)A_0e^{f_0(\theta)}B_{j+1}(\theta)=A_{j+1}e^{f_{j+1}(\theta)}.
\end{equation}
Now we give the desired estimate and finish one step of the iteration.

First the construction $A_{j+1}= R_{\overline{\xi_{j+1}}+\langle n_j,\alpha\rangle}$ implies that if $j+1=j_i$ is the resonant step, then 
$\nu^{j+1}=0.$
By \eqref{estima8} and \eqref{newf}, we have 
\begin{equation}\label{add2}
\|f_{j+1}\|_{h_{j+1}}  \leq e^{2\pi h|n_j|}\| \bar{f}_{j+1}\|_{h_{j+1}}\leq\epsilon_j e^{-h_{j+1}\varepsilon_j^{-\frac{1}{20\tau}}}:=\epsilon_{j+1}.
\end{equation}
By  \eqref{estima19} and \eqref{add1}, we have 
\begin{align*}
\|B_{j+1}\|_{h_{j+1}}\leq \|\tilde{B}_j\|_{h_{j+1}}\|B_j\|_{h_{j+1}}\leq 2\frac{|n_j|^\tau}{\kappa}e^{2\pi h_j|n_j|}\epsilon_{j}^{-\frac{1}{100}}\leq \epsilon_{j+1}^{-\frac{1}{100}}.
\end{align*}
Moreover, \eqref{estima7} and \eqref{estima16}  imply that
\begin{eqnarray*}
&& \|\tilde{B}_j-id\|_{h_{j+1}- \frac{\delta}{2\pi}} \\
&\leq& \|\bar{P}_j-id\|+\|\bar{Y}_{j}\|_{h_{j+1}}+\|R_{\frac{\langle n_j,\theta\rangle}{2}}\bar{U}_{j+1}R_{-\frac{\langle n_j,\theta\rangle}{2}}-id\|_{h_{j+1}-\frac{\delta}{2\pi}}\\
&\leq& \|\bar{P}_j-id\|+\|\bar{Y}_{j}\|_{h_{j+1}}+ \|R_{-\frac{\langle n_j,\theta\rangle}{2}}\|^2_{h_{j+1}-\frac{\delta}{2\pi}}\|\bar{U}_{j+1}-id\|\\
&\leq& \|\bar{P}_j-id\|+ \epsilon_j^{\frac{1}{2}}+ \epsilon_j^{\frac{1}{2}}e^{-2\pi (h_j-h_{j+1})|n_j|}.
\end{eqnarray*}
consequently Lemma \ref{liouass} implies that
 \begin{align*}
\|\tilde{B}_j-id\|_{h_{j+1}-\frac{\delta}{2\pi}} \leq \left\{ \begin{array}{ccc} \tilde{\epsilon}_{i_p} ^{\frac{1}{4}} & j +1= j_1\\
 \epsilon_{j_i-1} ^{\frac{1}{4}} &   j+1 =  j_i >j_1
\end{array} \right.\end{align*}
which means  $B_{j+1}$  satisfies the estimate \eqref{estima1} with $j+1$ in place of $j$.

Once we have this, we finished the whole step of the iteration, \eqref{estima18}, \eqref{estima2-1} and \eqref{estima1} imply the reducibility of $(\alpha,A)$.

\subsubsection{Part 2: Continuity}

Finally, assume $A,A'\in L_\alpha(\gamma,\delta)$ satisfy $\|A'-A\|_h\leq \eta$ where $\eta$ is arbitrary small such that \begin{equation}\label{initial-1}
    \|Z_{p_0}\|^2_{\frac{h+\tilde{h}}{2}}\|A-A'\|_h\leq \epsilon_0,
\end{equation} then there is $j\geq 0$ such that $$\epsilon_{j+1}<\|A-A'\|_h\leq \epsilon_{j}.$$ By the above iteration, there exist $B_{j}\in C^\omega_{\tilde{h}}(\T^d,\SL(2,\R))$ and $f_{j}\in C^\omega_{\tilde{h}}(\T^d,\sl(2,\R))$, such that
\begin{equation}\label{equation1}
B_{j}^{-1}(\theta+\alpha)Z_{p_0}^{-1}(\theta+\alpha)A(\theta)Z_{p_0}(\theta)B_{j}(\theta)=R_{\phi_j}e^{f_{j}(\theta)}
\end{equation}
with 
\begin{equation}\label{equation2}
\|B_{j}\|_{\tilde{h}}\leq \epsilon_{j}^{-\frac{1}{100}},\ \ \|f_{j}\|_{\tilde{h}}\leq \epsilon_{j},
\end{equation}
Hence \eqref{initial-1} and \eqref{equation2} allow us to conclude
\begin{equation*}
B_{j}^{-1}(\theta+\alpha)Z_{p_0}^{-1}(\theta+\alpha)A'(\theta)Z_{p_0}(\theta) B_{j}(\theta)=R_{\phi_j}e^{f'_{j}(\theta)}
\end{equation*}
with $\|f'_{j}\|_{\tilde{h}}\leq \epsilon^{\frac{1}{2}}_{j}$. Note that \eqref{reducibility condition-31}-\eqref{reducibility condition-1} are always satisfied with $(\alpha, \tilde{A}_{p_0}e^{\tilde{f}_{p_0}})$ replaced by $(\alpha,R_{\phi_j}e^{f_j})$ or $(\alpha,R_{\phi_j}e^{f'_j})$. Thus, there are $Y,Y'\in C^\omega_{\tilde{h}-\frac{\delta}{2\pi}}(\T,sl(2,\R))$ such that
$$
e^{-Y(\theta+\alpha)}R_{\phi_j}e^{f_{j}(\theta)}e^{Y(\theta)}=R_{\phi},\ \ e^{-Y'(\theta+\alpha)}R_{\phi_j}e^{f'_{j}(\theta)}e^{Y'(\theta)}=R_{\phi'}
$$
with estimate
$$
\|Y'\|_{\tilde{h}-\frac{\delta}{2\pi}},\ \ \|Y\|_{\tilde{h}-\frac{\delta}{2\pi}}\leq \epsilon^{\frac{1}{8}}_{j}
$$
Denote $\deg{Z_{p_0}B_j}=k_j$,  and  let 
\begin{eqnarray*}
   B_{A}(\theta)&=&Z_{p_0}(\theta)B_j(\theta)e^{Y(\theta)}R_{-\frac{\langle k,\theta\rangle}{2}},\\
   B_{A'}(\theta)&=&Z_{p_0} (\theta)B_j(\theta)e^{Y'(\theta)}R_{-\frac{\langle k,\theta\rangle}{2}},
\end{eqnarray*}
then $\deg B_{A}=\deg B_{A'}=0$. By \eqref{add2}, we have
$$
\|B_A-B_{A'}\|_{\tilde{h}-\frac{\delta}{2\pi}}\leq C\|B_j\|^4_{\tilde{h}}\epsilon_j^{\frac{1}{8}}\leq \epsilon_j^{\frac{1}{18\tau}}\leq \frac{1}{-\ln|\epsilon_{j+1}|}\leq \frac{1}{-\ln\|A-A'\|_h}.
$$
Thus, we finish the whole proof.

\end{pf}

\subsection{Optimal arithmetic reducibility-global case}
We denote 
$$
\mathcal{AR}_h=\{E:(\alpha.A_E) \text{ is almost reducible in the strip $|\Im z|<h$}\}.
$$
As a direct consequence of Theorem \ref{reducibilitylocal}, we have the following:
\begin{Corollary}\label{Amoredu}
Let $\alpha\in DC_d(\kappa,\tau)$, $\Sigma_{V,\alpha} =\mathcal{AR}_h$, $\rho(E)\in L_\alpha(\gamma,\delta)$ with $2\pi h>2\pi \tilde{h}>\delta$. Then $(\alpha, A_E)$ is analytically reducible. Moreover, there exist $C=C(\gamma,\delta,\alpha,V)$ and $B_E\in C^\omega_{\tilde{h}-\frac{\delta}{2\pi}}(\T,SL(2,\R))$, such that
$$
B_E(\theta+\alpha)^{-1}A_E(\theta)B_E(\theta)=R_{\rho(E)},
$$
with estimates $\|B_E\|_{\tilde{h}-\frac{\delta}{2\pi}}\leq C$, moreover $B_E$ is continuous on $\mathcal{E}_\alpha(\gamma,\delta)$ in $\|\cdot\|_{\tilde{h}-\frac{\delta}{2\pi}}$ where
$$
\mathcal{E}_\alpha(\gamma,\delta)=\{E:\rho(E)\in L_\alpha(\gamma,\delta)\}.
$$
\end{Corollary}
\begin{pf}
By the assumption that $\Sigma_{V,\alpha}=\mathcal{AR}_h$. One can conclude for any $E\in \Sigma_{V,\alpha}$, and for any  $\eta>0$, there exist $\Phi_{E}\in C^{\omega}(\T^d, \PSL(2,\R))$ such that
\begin{equation*}
\Phi_{E}(\theta+\alpha)^{-1}A_E(\theta)\Phi_{E}(\theta)=C_{E}e^{f_E(\theta)},
\end{equation*}
with estimates
\begin{equation}\label{estim2}
\|\Phi_E\|_{h}\leq \Gamma=\Gamma(V,\alpha,\eta), \qquad \|f_E\|_{h} \leq \eta.
\end{equation}
Indeed, $(\alpha, A_E)$ is  almost reducible but not uniformly hyperbolic, one can always take $C_{E}=R_{\phi(E)} \in SO(2,\R)$ \cite[Corollay 4.2]{yzhou}.

Indeed, one can always take $\eta$ small enough such that
\begin{equation}\label{small}
\eta  \leq \bar{\epsilon}_0\leq  \frac{D_0(\kappa,\tau,1)}{\|R_{\phi(E)}\|^{C_0}}(\frac{h- \tilde{h}}{2})^{C_0\tau},
\end{equation}
where $D_0(\kappa,\tau,1)$ is the constant defined in Theorem  \ref{reducibilitylocal}. 
Just note the smallness  of $\bar{\epsilon}_0$ can be taken uniform with respect to $R_{\phi} \in {\rm SO}(2,\R)$.
By footnote 5 of \cite{avila1} and \eqref{estim2}, we have \begin{equation}\label{deg}|\deg{\Phi_{E}}|\leq  C |\ln \Gamma|:=\Gamma_1\end{equation} for some constant $C=C(V,\alpha)>0$. Moreover, by \eqref{rotation number},
\begin{equation}\label{royyy}
  2\rho(\alpha,R_{\phi(E)}e^{f_{E}(\theta)})=2\rho(E)-\langle \deg{\Phi_E},\alpha\rangle.  
\end{equation}

Therefore, if $|E-E'|\leq\Gamma^{-4}\bar{\epsilon}^4_{0}$, then
\begin{equation}\label{conp}
\Phi_{E}(\theta+\alpha)^{-1}A_{E'}(\theta)\Phi_{E}(\theta)=R_{\phi(E)}e^{f_{E'}(\theta)},
\end{equation}
with estimate  $\|f_{E'}(\theta)\|_{h}\leq \bar{\epsilon}_{0}$.
On the other hand, if $\rho(E)=\rho(E')\in L_\alpha(\gamma,\delta)$, by \eqref{royyy} and \eqref{conp}, this implies that  there is $\gamma'>0$ such that
$$
\rho(\alpha,R_{\phi(E)}e^{f_{E}(\theta)}) = \rho(\alpha,R_{\phi(E)}e^{f_{E'}(\theta)})   \in LC(\gamma',\delta).
$$
By Theorem \ref{reducibilitylocal}, there are $\tilde{B}_E, \tilde{B}_{E'}\in C^\omega_{h-\frac{\delta}{2\pi}}(\T^d,SL(2,\R))$ such that
\begin{equation*}
\tilde{B}_{E}^{-1}(\theta+\alpha)R_{\phi(E)} e^{f_{E}(\theta)}\tilde{B}_{E}(\theta)=R_{\rho(E)-\frac{\langle\deg{\Phi_E},\alpha\rangle}{2}},
\end{equation*}
\begin{equation*}
\tilde{B}_{E'}^{-1}(\theta+\alpha)R_{\phi(E)} e^{f_{E'}(\theta)}\tilde{B}_{E'}(\theta)=R_{\rho(E')-\frac{\langle\deg{\Phi_E},\alpha\rangle}{2}}
\end{equation*}
with 
\begin{equation*}
\lim\limits_{E\rightarrow E'}\|\tilde{B}_E-\tilde{B}_{E'}\|_{\tilde{h}-\frac{\delta}{2\pi}}=0,
\end{equation*}
\begin{equation*}
\|\tilde{B}_{E}\|_{\tilde{h}-\frac{\delta}{2\pi}}\leq C(\gamma,\delta,\alpha,V).
\end{equation*}
Let $B_E(\theta)=\Phi_E(\theta) \tilde{B}_E(\theta)R_{-\frac{\langle \Phi_E,\theta\rangle}{2}}$ and $B_{E'}(\theta)=\Phi_E(\theta) \tilde{B}_{E'}(\theta)R_{-\frac{\langle \Phi_E,\theta\rangle}{2}}$, we finish the whole proof.
\end{pf}
\noindent
{\bf Proof of Theorem \ref{reducibility main2}:} Note that $\alpha\in DC$ and $(\alpha,A)$ is subcritical in the strip $\{|\Im z|<h\}$. By Avila's solution of almost reducibility conjecture \cite{avila1,avila2}, for any $\e>0$, $(\alpha,A)$ is almost reducible in the strip $\{|\Im z|<h-\e\}$.  Note by the definition of $\delta(\alpha, \rho(\alpha,A))$, for any $\e>0$, there exists $\gamma>0$, such that
$$\rho(\alpha,A)\in LC_\alpha(\gamma,\delta(\alpha,\rho(\alpha,A))+\e),$$then the result follows from  Corollary \ref{Amoredu}. \qed
\\

\subsection{Optimal arithmetic localization}

For any $E\in \mathcal{E}_\alpha(\gamma,\delta)$, by Corollary \ref{Amoredu}, there exists $B_E\in C^\omega(\T^d,SL(2,\R))$ such that
\begin{equation}\label{redsch}
B_{E}^{-1}(\theta+\alpha)A_{E}(\theta)B_{E}(\theta)=R_{\rho(E)}.
\end{equation}
We define a vector-valued function $u_E:\mathcal{E}_\alpha(\gamma,\delta)\rightarrow \ell^2(\Z)$ as the following,
\begin{equation}\label{def1}
u_E(n)=\frac{\hat{b}_E(n)}{\|b_E\|_{L^2}}=\frac{\int_{\T}b_E(\theta)e^{2\pi i\langle n,\theta\rangle}d\theta}{\|b_E\|_{L^2}},
\end{equation}
where $b_E(\theta)=\frac{ib_E^{11}(\theta)-b_{E}^{12}(\theta)}{2i}$ and $B(\theta)=\begin{pmatrix}b_E^{11}(\theta)&b_E^{12}(\theta)\\ b_E^{21}(\theta)&b_E^{22}(\theta)\end{pmatrix}$.

For any fixed $x\in L_\alpha(\delta)=\cup_{\gamma>0}L_\alpha(\gamma,\delta)$, we denote by 
$$
E_m(x)=\begin{cases}
\rho^{-1}(T^mx)&T^mx\in [0,\frac{1}{2})\\
\rho^{-1}(-T^mx+1)&T^mx\in [\frac{1}{2},1)
\end{cases}.
$$ 
We can define the following $\mathcal{R}$-measure,
\begin{Definition}[$\mathcal{R}$-measure]\label{Rm}
\label{defR}
$\nu_{x,\delta_n}:\mathcal{B}\rightarrow \R$ is defined as:
$$
\nu_{x,\delta_n}(B)=\sum\limits_{m\in N_x^B}|u_{E_m(x)}(m+n)|^2,
$$
for all $B$ in the Borel $\sigma$-algebra $\mathcal{B}$ of $\R$, where $N_x^B=\{m| E_m(x)\in B\}$. 
\end{Definition}
The $\mathcal{R}$-measure is well defined, as proved in \cite{gy}. Moreover, following the proof of \cite{gy,gk}, one has
\begin{Lemma}\label{property}
We have $\nu_{x,\delta_n}(\mathcal{E}_\alpha(\gamma,\delta))=|L_\alpha(\gamma,\delta)|$ for $a.e.$ $x$, where $|\cdot|$ is the Lebesgue measure.
\end{Lemma}
We fix $\gamma,\delta$ and $n$ in the following.
\begin{Lemma}\label{tail}
For any $\epsilon>0$, there exists $N_0(\gamma,\delta,V,\alpha,n,\epsilon)>0$ such that for all $x\in L_\alpha(\gamma,\delta)$, we have
$$
\mathcal{R}_{N_0}\nu_{x,\delta_n}(\mathcal{E}_\alpha(\gamma,\delta)):=\sum\limits_{|m|>N_0:T^mx\in L_\alpha(\gamma,\delta)}|u_{E_m(x)}(m+n)|^2\leq \epsilon.
$$
\end{Lemma}
\begin{pf}
Note that for any $T^mx\in L_\alpha(\gamma,\delta)$, by Corollary \ref{Amoredu}, there exists $B_{E_m(x)}\in C^\omega(\T^d,\SL(2,\R))$ such that
\begin{equation*}
B_{E_m(x)}^{-1}(\theta+\alpha)A_{E_m(x)}(\theta)B_{E_m(x)}(\theta)=R_{T^mx},
\end{equation*}
with $\left\|B_{E_m(x)}\right\|_{\tilde{h}-\frac{\delta}{2\pi}}\leq C(\gamma,\delta,\alpha,V)$. We denote by
$$
\widetilde{B}_E(\theta)=B_E(\theta)\frac{1}{2i}\begin{pmatrix}i&i\\ -1&1\end{pmatrix}=\begin{pmatrix}\tilde{b}_E^{11}(\theta)&\tilde{b}_E^{12}(\theta)\\ \tilde{b}_E^{21}(\theta)&\tilde{b}_E^{22}(\theta)\end{pmatrix}
$$
then 
\begin{equation}\label{f2}
\widetilde{B}_{E_m(x)}^{-1}(\theta+\alpha)A_{E_m(x)}(\theta)\widetilde{B}_{E_m(x)}(\theta)=\begin{pmatrix}e^{2\pi i T^mx}&0\\ 0&e^{-2\pi i T^mx}\end{pmatrix}.
\end{equation}

By \eqref{def1}, we have
$$
\left|u_{E_m(x)}(m+n)\right|^2=\frac{\left|\hat{\tilde{b}}_{E_m(x)}^{11}(m+n)\right|^2}{\left\|\tilde{b}_{E_m(x)}^{11}\right\|^2_{L^2}}.
$$
By \eqref{f2}, we have
$$
\tilde{b}_{E_m(x)}^{11}(\theta)=\tilde{b}_{E_m(x)}^{21}(\theta+\alpha)e^{2\pi iT^mx},
$$
By the fact that  $|\det{\widetilde{B}_{E_m(x)}}|=1$, one has
\begin{align}\label{redest12}
2\left\|\tilde{b}_{E_m(x)}^{11}\right\|_{L^2}=\left\|\tilde{b}_{E_m(x)}^{11}\right\|_{L^2}+\left\|\tilde{b}_{E_m(x)}^{21}\right\|_{L^2}\geq\left\|\widetilde{B}_{E_m(x)}\right\|_{C^0}^{-1}.
\end{align}
It follows that if $T^mx\in L_\alpha(\gamma,\delta)$, then
$$
\left|u_{E_m(x)}(m+n)\right|^2\leq Ce^{-(\tilde{h}-\frac{\delta}{2\pi})|m+n|}.
$$
Thus for any $\epsilon>0$, there exists $N_0(\gamma,\delta,V,\alpha,n,\epsilon)>0$ such that for all $x\in L_\alpha(\gamma,\delta)$,
$$
\sum\limits_{|m|>N_0:T^mx\in L_\alpha(\gamma,\delta)}|u_{E_m(x)}(m+n)|^2\leq \epsilon.
$$
\end{pf}
\begin{Lemma}\label{continuity}
For any $N>0$ and $\epsilon>0$, there exists $\eta(\gamma,\delta,V,\alpha,N,\epsilon)>0$ such that
$$
\left|\mathcal{T}_N\nu_{x,\delta_n}(\mathcal{E}_\alpha(\delta))-\mathcal{T}_N\nu_{x',\delta_n}(\mathcal{E}_\alpha(\delta))\right|\leq \epsilon
$$
for any $x,x'\in L_\alpha(\gamma,\delta)$ with $|x-x'|\leq \eta$ where 
$$
\mathcal{T}_N\nu_{x,\delta_n}(\mathcal{E}_\alpha(\delta)):=\nu_{x,\delta_n}(\mathcal{E}_\alpha(\delta))-\mathcal{R}_N\nu_{x,\delta_n}(\mathcal{E}_\alpha(\delta)).
$$
\end{Lemma}
\begin{pf}
For any fixed $N$ and any $x\in L_\alpha(\gamma,\delta)$, we have
$$
\|x+\langle k,\alpha\rangle+\langle n,\alpha\rangle\|_{\R/\Z}\geq \gamma e^{-\delta |k+n|}\geq \gamma e^{-\delta N} e^{-\delta|n|},\ \  |k|\leq N.
$$
Let $\gamma_1=\gamma e^{-\delta N}$. By Corollary \ref{Amoredu}, there exists $\eta(\gamma,\delta,V,\alpha,N,\epsilon)$, such that the following holds: if $T^kx,T^kx'\in L_\alpha(\gamma_1,\delta)$ and $|T^kx-T^kx'|=|x-x'|<\eta$, then there exist $B_{E_k(x)},B_{E_k(x')}\in C_{\tilde{h}-\frac{\delta}{2\pi}}^\omega(\T^d,SL(2,\R))$, such that
\begin{align*}
B_{E_k(x)}^{-1}(\theta+\alpha)A_{E_k(x)}(\theta)B_{E_k(x)}(\theta)=R_{T^kx},
\end{align*}
\begin{align*}
B_{E_k(x')}^{-1}(\theta+\alpha)A_{E_k(x')}(\theta)B_{E_k(x')}(\theta)=R_{T^kx'},
\end{align*}
with estimates
\begin{align}\label{redests1}
\|B_{E_k(x)}\|_{\tilde{h}-\frac{\delta}{2\pi}},\ \ \|B_{E_k(x')}\|_{\tilde{h}-\frac{\delta}{2\pi}}\leq C(\gamma,\delta,V,\alpha,N),
\end{align}
\begin{align}\label{redest3}
\|B_{E_k(x)}-B_{E_k(x')}\|_{\tilde{h}-\frac{\delta}{2\pi}}\leq \frac{\epsilon C^{-4}}{500(2N+1)}.
\end{align}
On the one hand,
\begin{align}\label{redest10}
\left|\|b_{E_k(x)}\|_{L^2}-\|b_{E_k(x')}\|_{L^2}\right|&\leq \|b_{E_k(x)}-b_{E_k(x')}\|_{L^2}\\ \nonumber
&\leq \|B_{E_k(x)}-B_{E_k(x')}\|_{C^0}\leq \frac{\epsilon C^{-4}}{500(2N+1)}.
\end{align}
On the other hand, by \eqref{redest3},
\begin{align}\label{redest11}
\left|\hat{b}_{E_k(x)}(n)-\hat{b}_{E_k(x')}(n)\right|&\leq \left|\int\limits_{\T}\left(b_{E_k(x)}(\theta)-b_{E_k(x')}(\theta)\right)e^{-2\pi i\langle n,\theta\rangle}d\theta\right|\\ \nonumber
&\leq\|B_{E_k(x)}-B_{E_k(x')}\|_{\tilde{h}-\frac{\delta}{2\pi}}\leq \frac{\epsilon C^{-4}}{500(2N+1)}.
\end{align}
\eqref{redest10} and \eqref{redest11} imply for any $n\in\Z^d$,
\begin{align*}
|u_{E_k(x)}(n)-u_{E_k(x')}(n)|&=\left|\frac{\hat{b}_{E_k(x)}(n)}{\|b_{E_k(x)}\|_{L^2}}-\frac{\hat{b}_{E_k(x')}(n)}{\|b_{E_k(x')}\|_{L^2}}\right|\\
&=\frac{\left|\hat{b}_{E_k(x)}(n)\|b_{E_k(x')}\|_{L^2}-\hat{b}_{E_k(x')}(n)\|b_{E_k(x)}\|_{L^2}\right|}{\|b_{E_k(x)}\|_{L^2}\|b_{E_k(x')}\|_{L^2}}\\
&\leq \frac{\epsilon C^{-4}}{500(2N+1)}\frac{\|b_{E_k(x)}\|_{L^2}+C}{\|b_{E_k(x)}\|_{L^2}\|b_{E_k(x')}\|_{L^2}}.
\end{align*}
As we have
$$
\|b_{E_k(x)}\|_{L^2}\|b_{E_k(x')}\|_{L^2}\geq \frac{1}{4C^2},
$$
these imply that 
\begin{align*}
|u_{E_k(x)}(n)-u_{E_k(x')}(n)|&\leq \frac{\epsilon C^{-4}}{500(2N+1)}\frac{\|b_{E(T^kx)}^{11}\|_{L^2}+C}{\|b_{E(T^kx)}^{11}\|_{L^2}\|b_{E(T^kx')}^{11}\|_{L^2}}\\
&\leq \frac{\epsilon}{100(2N+1)}.
\end{align*}
By definition \ref{defR} and \eqref{def1}, one has
\begin{align*}
&\ \ \ \ |\mathcal{T}_N\nu_{x,\delta_n}(\mathcal{E}_\alpha(\delta))-\mathcal{T}_N\nu_{x',\delta_n}(\mathcal{E}_\alpha(\delta))|\\
&=|\sum\limits_{|k|\leq N}|u_{E(T^kx)}(n+k)|^2-\sum\limits_{|k|\leq N}|u_{E(T^kx')}(n+k)|^2|\\
&\leq \frac{\epsilon}{50(2N+1)}(2N+1)\leq \epsilon.
\end{align*}
\end{pf}
\begin{Lemma}\label{diophantine}
For any $x\in L_\alpha(\gamma,\delta)$, we have 
$$
(x-\sigma,x+\sigma)\cap L_\alpha(\gamma,(1+\e)\delta)>0
$$ 
for any sufficiently small $\sigma>0$.
\end{Lemma}
\begin{pf}
Let 
$$
\Theta_k=\{x\in[0,1):\|2x-\langle k,\alpha\rangle\|_{\R/\Z}<\gamma e^{-(1+\e)\delta|k|},
$$ 
then 
$L_\alpha(\gamma,(1+\e)\delta)=[0,1)\backslash\cup_{k\in\Z}\Theta_k$. Thus for any $\sigma>0$ sufficiently small, we have
$$
(x-\sigma,x+\sigma)\cap L_\alpha(\gamma,(1+\e)\delta)=(x-\sigma,x+\sigma)\backslash\cup_{k\in\Z^d}\Theta_k.
$$
Notice that if $\Theta_{k}\cap (x-\sigma,x+\sigma)\neq \emptyset$, then 
$$
\gamma e^{-\delta|k|}\leq\|\langle k,\alpha\rangle+2x\|_{\R/\Z}\leq 2\sigma.
$$
It follows that $|k|\geq \frac{\ln\frac{\gamma}{2\sigma}}{\delta}$, thus $|\Theta_k|\leq C(\gamma)\sigma^{1+\e}$ which implies that
$$
|(x-\sigma,x+\sigma)\cap L_\alpha(\gamma,(1+\e)\delta)|\geq 2\sigma-\sum\limits_{|k|\geq \frac{\ln\frac{\gamma}{2\sigma}}{\delta}}\Theta_k\geq \sigma>0.
$$
\end{pf}

Once we have these preparing lemmas, we arrive at the following result: 
\begin{Theorem}\label{longgeneral}
Assume that $\alpha\in DC_d(\kappa,\tau)$ and $V\in C_h^{\omega}(\T^d,\R)$ with $\Sigma_{V,\alpha}=\mathcal{AR}$. If  $2\pi h>\delta(\alpha,x)$, then  
$L_{V,\alpha,x}$ has Anderson localization.
\end{Theorem}
\begin{pf} For any fixed $x\in L_\alpha(\gamma,\delta)$ with $2\pi h>\delta(\alpha,x)$ and $\e=\frac{2\pi h-\delta}{10}$, by Lemma \ref{diophantine} and Lemma \ref{property},  there exists a sequence $x_k\in L_\alpha(\gamma,(1+\e)\delta)$ such that $x_k\rightarrow x$ and 
$$
\nu_{x_k,\delta_n}(\mathcal{E}_\alpha(\gamma,(1+\e)\delta))=\left|L_\alpha(\gamma,(1+\e)\delta)\right|.
$$
By Lemma \ref{tail}, there exists $N_0(\gamma,(1+\e)\delta,V,\alpha,n,\epsilon)>0$ such that 
$$
\mathcal{R}_{N_0}\nu_{x_k,\delta_n}(\mathcal{E}_\alpha(\gamma,(1+\e)\delta))\leq \gamma.
$$
By Lemma \ref{continuity},
$$
\mathcal{T}_{N_0}\nu_{x,\delta_n}(\mathcal{E}_\alpha((1+\e)\delta))=\lim\limits_{k\rightarrow \infty}\mathcal{T}_{N_0}\nu_{x_k,\delta_n}(\mathcal{E}_\alpha((1+\e)\delta)),
$$
Therefore, we have
\begin{align*}
\nu_{x,\delta_n}(\mathcal{E}_\alpha((1+\e)\delta))&\geq \mathcal{T}_{N_0}\nu_{x,\delta_n}(\mathcal{E}_\alpha((1+\e)\delta))\geq \limsup\limits_{k\rightarrow\infty}\mathcal{T}_{N_0}\nu_{x_k,\delta_n}(\mathcal{E}_\alpha(\gamma,(1+\e)\delta))\\
&\geq\left|L_\alpha(\gamma,(1+\e)\delta)\right|-\gamma\geq 1-2\gamma.
\end{align*}
Let $\gamma\rightarrow 0$, we have 
$$
1\leq \nu_{x,\delta_n}(\mathcal{E}_\alpha((1+\e)\delta))\leq\mu^{pp}_{x,\delta_n}(\mathcal{E}_\alpha((1+\e)\delta))\leq\mu_{x,\delta_n}(\mathcal{E}_\alpha((1+\e)\delta))\leq1,
$$
where $\mu_{x,\delta_n}$ is the spectral measure of $L_{V,\alpha,x}$ defined by
$$
\langle\delta_n,\chi_{B}(L_{V,\alpha,x})\delta_n\rangle=\int_{\R}\chi_{B} d\mu_{\theta,\delta_n}.
$$
Thus we finish the proof. 
\end{pf}

\smallskip
{\bf Proof of Corollary \ref{long}:} By the assumption, if $\lambda\geq \lambda_0(\alpha,V,d,\delta)$ is large enough, the dual Schr\"odinger cocycle is almost reducible by Proposition \ref{reducibility}, then the result follows directly from Theorem \ref{longgeneral}. \qed

\smallskip
{\bf Proof of Corollary \ref{main}:} By Avila's global theory \cite{avila2015global}, if $|\lambda|<1$, $(\alpha,S_E^\lambda)$ is subcritical in the strip $\{|\Im z|<-\frac{\ln \lambda}{2\pi}\}$, 
 and  for any $\e>0$, $(\alpha,S_E^{\lambda})$ is almost reducible in the strip $\{|\Im z|<-\frac{\ln \lambda}{2\pi}(1- \frac{\e}{2})\}$ \cite{avila,avila2}. The result again follows directly from Theorem \ref{longgeneral}. \qed

\section{Growth of the cocycle}
In previous section, we apply the  structured quantitative  almost reducibility (Proposition  \ref{reducibility}) in the case $2\pi h>\delta(\alpha,\rho(\alpha,A))$, to prove the cocycle is reducible. In this section, we apply the  structured quantitative  almost reducibility in the opposite case $2\pi h\leq \delta(\alpha,\rho(\alpha,A))$, to estimate the growth of cocycles.

Recall that Proposition  \ref{reducibility} implies that $(\alpha,A)$ is almost reducible: 
$$
B_j^{-1}(\theta+\alpha)A(\theta)B_j(\theta)=A_je^{f_j(\theta)}=M^{-1}exp \left(
\begin{array}{ccc}
 i t^j &  \nu^j\\
\bar{ \nu}^j &  -i t^j
 \end{array}\right)Me^{f_j(\theta)}
$$
Denote $spec(A_{j})=\{e^{ i \xi_{j}},e^{- i\xi_{j}} \}$, $(\alpha,A)^n=(n\alpha,  \mathcal{A}_n)$, then we have the following:

\begin{Lemma}\label{growth}
Under the assumption of Proposition  \ref{reducibility}, then we have the following:
\begin{enumerate}
\item If $|\frac{\nu_j}{\xi_j}|<\frac{1}{2}$, then for any $e^{ \frac{\pi}{128} |n_{i_{j}}|(h-\tilde{h})}\leq  n<  e^{ \frac{\pi}{8} |n_{i_{j+1}}|(h-\tilde{h})}$, we have 
$$ C^{-1} (\frac{ |\ln n|}{h-\tilde{h}})^{-4\tau}      \leq \|\mathcal{A}_n\|_0\leq C (\frac{ |\ln n|}{h-\tilde{h}})^{4\tau}.$$
\item If $|\frac{\nu_j}{\xi_j}| \geq \frac{1}{2}$, then we have the following:
\begin{itemize}  
\item If $e^{ \frac{\pi}{128} |n_{i_{j}}|(h-\tilde{h})}\leq n< \min\{\frac{1}{\xi_j}, e^{\frac{\pi}{32} |n_{i_{j+1}}|(h-\tilde{h})}\}$, then 
$$
C^{-1} (|n||\nu_j|+1)(\frac{ |\ln n|}{h-\tilde{h}})^{-4\tau}\leq \|\mathcal{A}_n\|_0\leq C (|n||\nu_j|+1) (\frac{ |\ln n|}{h-\tilde{h}})^{4\tau}.
$$
\item If $\frac{1}{\xi_j}  \leq n  < e^{ \frac{\pi}{32} |n_{i_{j+1}}|(h-\tilde{h})},$ then 
$$ C^{-1} \sqrt{1+ 2\sin^2(n\xi_j) |\frac{\nu_j}{\xi_j}|^2 } (\frac{ |\ln n|}{h-\tilde{h}})^{-4\tau}  \leq   \|\mathcal{A}_n\|_0 \leq  C \sqrt{1+ 2\sin^2(n\xi_j) |\frac{\nu_j}{\xi_j}|^2} (\frac{ |\ln n|}{h-\tilde{h}})^{4\tau}. $$
\end{itemize}

\end{enumerate}

\end{Lemma}

\begin{Remark}\label{grow-rem}
The inequality $\sqrt{1+ 2\sin^2(n\xi_j) |\frac{\nu_j}{\xi_j}|^2} \leq  C (|n||\nu_j|+1)$ suggests that according to Lemma \ref{growth}, the quantity $\|\mathcal{A}_n\|_0$ initially experiences growth, followed by a subsequent decay. The maximum value is attained when $n\sim 1/\xi_j$. 
\end{Remark}

\begin{pf}
By Proposition  \ref{reducibility}, there exist $B_j\in C_{\tilde{h}}^\omega(\T^d, \PSL(2,\R))$,  $f_j\in C_{\tilde{h}}^\omega(\T^d, \sl(2,\R))$  such that 
$$
B_j^{-1}(\theta+\alpha)A(\theta)B_j(\theta)=A_je^{f_j(\theta)}=M^{-1}exp \left(
\begin{array}{ccc}
 i t^j &  \nu^j\\
\bar{ \nu}^j &  -i t^j
 \end{array}\right)Me^{f_j(\theta)}
$$
with estimates 
\begin{eqnarray}
 \label{es2-2}
|\nu^j| &\leq& 2 \epsilon_{i_j}^{\frac{15}{16}}e^{-2\pi|n_{i_j}|\tilde{h}},\\
   \label{es0-2} \|f_j\|_{\tilde{h}} &\leq&  \epsilon_{i_j}e^{- \frac{\pi}{8} |n_{i_{j+1}}|(h-\tilde{h})},\\
   \label{es-2}
\|B_j\|_0&<&C(\alpha)|n_{i_j}|^{2\tau}.
\end{eqnarray}

To control the growth of the cocycles, the starting point is the following simple observation:
\begin{Lemma}[\cite{afk, ZhouW12}]\label{afk}
We have that
\begin{equation*}
	M_l(\id+y_l)\cdots M_0(\id+y_0)=M^{(l)}(\id+y^{(l)}),
\end{equation*}
where $M^{(l)}=M_l\cdots M_0$ and
\begin{equation}\label{error-growth}
	\| y^{(l)}\| \leq \mathrm{e}^{\sum_{k=0}^{l}\| M^{(k)}\|^2\| y_k\|}-1.
\end{equation}
\end{Lemma}

In case (1), by  Lemma \ref{Dominate matrix},  there exists $U_j \in \SL(2,\R)$  with  $\|U_j -\id \|\leq  |\frac{\nu_j}{\xi_j}|$   such that 
$
U_j^{-1}A_jU_j= R_{\xi_j},
$
which implies that 
\begin{equation}\label{conju-1}
U_j^{-1}B_j^{-1}(\theta+\alpha)A(\theta)B_j(\theta)U_j= R_{\xi_j}e^{f_j'(\theta)}
\end{equation}
where  $f_j'(\theta)= U_j^{-1}f_j(\theta)U_j$. Consequently, for any $n< e^{\frac{\pi}{8} |n_{i_{j+1}}|(h-\tilde{h})},$ by Lemma \ref{afk}, we have estimate
$$
\frac{1}{2} \leq \|\prod\limits_{k=0}^n R_{\xi_j}e^{f_{j}'(\theta+k\alpha)}\|_0\leq 2,
$$
by \eqref{conju-1},  we further have  estimate
$$
\frac{1}{2} \|B_jU_j \|^{-2}_0\leq \|\mathcal{A}_n\|_0\leq 2  \|B_jU_j \|^2_0,$$
thus by the assumption  $e^{\frac{\pi}{128} |n_{i_{j}}|(h-\tilde{h})} \leq n$, and \eqref{es-2}, we have 
$$ C^{-1} (\frac{ |\ln n|}{h-\tilde{h}})^{-4\tau}   \leq \|\mathcal{A}_n\|_0\leq  C (\frac{ |\ln n|}{h-\tilde{h}})^{4\tau}.$$

In case (2),  we actually can assume (otherwise the argument would be more simpler)
\begin{equation*}
|\xi_j| \geq  e^{- \frac{\pi}{32} |n_{i_{j+1}}|(h-\tilde{h})}.
\end{equation*}
We distinguish the proof into two cases, according to the time scales. \\

\noindent \textbf{Case (2)-1: $e^{ \frac{\pi}{128} |n_{i_{j}}|(h-\tilde{h})}< n< \frac{1}{\xi_j}$.}
First  by  Lemma \ref{Dominate matrix},  there exists $U_j \in SL(2,\C)$ such that 
\begin{equation*}
U_j^{-1}A_jU_j=\begin{pmatrix}e^{ i \xi_j}&\nu'_j \\ 0&e^{- i \xi_j}\end{pmatrix},
\end{equation*}
with estimate
\begin{equation*} \|U_j\|\leq 2, \ \ |\nu'_j|\leq 4|\nu_j|.
\end{equation*}
Consequently, 
\begin{equation}\label{conju-2}
U_j^{-1}B_j^{-1}(\theta+\alpha)A(\theta)B_j(\theta)U_j=  A_j'  e^{f_j'(\theta)}= \begin{pmatrix}e^{ i \xi_j}&\nu'_j \\ 0&e^{- i \xi_j}\end{pmatrix}e^{f_j'(\theta)}
\end{equation}
 To estimate the growth of the cocycle, we recall the following:
\begin{Lemma}[\cite{aj}]\label{comp}
Assume that
$$
T=\begin{pmatrix} e^{2\pi i\theta}&c\\0&e^{-2\pi i\theta}
\end{pmatrix},
$$
then we have 
$$
T^{n}=\begin{pmatrix}
e^{2\pi i(n\theta)}&t_{n}\\
0&e^{-2\pi i(n\theta)},
\end{pmatrix}
$$
where $t_{n}=ce^{2\pi i(n-1)\theta}\frac{e^{-4\pi in\theta}-1}{e^{-4\pi i\theta}-1}$. Thus, for any $n$ with $|n|\|\theta\|_{\R/\Z}<1$, we have 
$$ \frac{1}{2}c|n| \leq  |t_n|\leq  c|n|.$$
\end{Lemma}

Consequently, for any $n< \frac{1}{\xi_j}<  e^{\frac{\pi}{32} |n_{i_{j+1}}|(h-\tilde{h})},$ by Lemma \ref{afk} and Lemma \ref{comp}, we have 
$$
1+\frac{1}{2}|n||\nu_j'|\leq \|\prod\limits_{k=0}^n A_{j}'e^{f_{j}'(\theta+k\alpha)}\|_0\leq 1+2|n||\nu_j'|,
$$
similar as case (1), by  \eqref{es2-2}, \eqref{es-2}, and \eqref{conju-2} we have
$$
C^{-1} (|n||\nu_j|+1)(\frac{ |\ln n|}{h-\tilde{h}})^{-4\tau}\leq \|\mathcal{A}_n\|_0\leq C (|n||\nu_j|+1) (\frac{ |\ln n|}{h-\tilde{h}})^{4\tau}.
$$

\noindent \textbf{Case (2)-2:  $\frac{1}{\xi_j}  \leq n  < e^{ \frac{\pi}{32} |n_{i_{j+1}}|(h-\tilde{h})}$.}
By  Lemma \ref{Dominate matrix},  there exists $U_j \in \SL(2,\R)$  with  $\|U_j\|^2 \leq  |\frac{4\nu_j}{\xi_j}|$   such that 
$
U_j^{-1}A_jU_j= R_{\xi_j}.
$
In this case, the crucial observation is the following:

\begin{Lemma}
Assume that $$A=exp \left(
\begin{array}{ccc}
 i t &  \nu\\
\bar{ \nu} &  -i t
 \end{array}\right)\in SU(1,1)$$ with $spec\{A\}=\{e^{ i\xi},e^{- i\xi}\}$, $0\neq \xi\in \R$, then there exists $U\in SU(1,1)$, such that
$$
U^{-1}AU=\left(\begin{array}{ccc}
 e^{ i\xi} &  0\\
0 &  e^{- i\xi}
 \end{array}\right),
$$
Consequently,  for any $ \theta\in \R$, denoting $\tilde{A}= U \left(\begin{array}{ccc}
 e^{ i\theta} &  0\\
0 &  e^{- i\theta}
 \end{array}\right) U^{-1} ,$ then we have
 $$\|\tilde{A}\|_{HS}^2= 2+ 4\sin^2 \theta |\frac{\nu}{\xi}|^2.$$ 
\end{Lemma}
\begin{pf}
By Lemma \ref{nf}, one can take 
\begin{equation*} 
\begin{array}{r@{}l}
U = (\cos 2\varphi)^{-\frac{1}{2}}
\begin{pmatrix}\cos\varphi &e^{2i\phi} \sin\varphi \\ e^{-2i\phi}\sin\varphi & \cos\varphi
\end{pmatrix},
\end{array}
\end{equation*}
then direct computations shows that 
\begin{equation*} 
\begin{array}{r@{}l}
\tilde{A}= \frac{1}{\cos 2\varphi}
\begin{pmatrix}  e^{i \theta} \cos^2 \varphi - e^{-i \theta} \sin^2 \varphi   &i e^{2i\phi} \sin \theta \sin2\varphi \\- i e^{2i\phi} \sin \theta \sin2\varphi& e^{-i \theta} \cos^2 \varphi - e^{i \theta} \sin^2 \varphi  
\end{pmatrix},
\end{array}
\end{equation*}
therefore one can further compute
\begin{eqnarray*} 
\|\tilde{A}\|_{HS}^2 &= &\frac{1}{ \cos^2 2\varphi } \left( 2\cos^2\varphi+2 \sin^4 \varphi+ (4\sin^2 \theta-1 ) \sin^2 2\varphi  \right)\\
&=&\frac{1}{ \cos^2 2\varphi } \left(  \cos^2 2\varphi+1 -\sin^2 2\varphi  + 4\sin^2 \theta  \sin^2 2\varphi   \right)\\
&=& 2 + 4  \sin^2 \theta  \tan^2 2\varphi
\end{eqnarray*}
then the result follows directly from \begin{equation*}
2\varphi = - \arctan \frac{|\nu|}{\sqrt{t^{2}-|\nu|^{2}}}.
\end{equation*}
\end{pf}

Apply Lemma \ref{afk}, we have 
$$\prod\limits_{k=0}^n A_{j}e^{f_{j}(\theta+k\alpha)}= U_j R_{n\xi_j}U_j^{-1} e^{\tilde{f}_n(\theta)},$$
thus if $n  < e^{ \frac{\pi}{32} |n_{i_{j+1}}|(h-\tilde{h})}$, then by \eqref{es2-2}, \eqref{es0-2} and \eqref{error-growth}, we have  estimate
\begin{eqnarray*} \| \tilde{f}_n\|_0 &\leq&  \sum_{k=0}^{n}\| U_j R_{k\xi_j}U_j^{-1} \|^2 \|f_{j}(\theta+k\alpha)\|_0    \\
&\leq & n \|U_j\|^4  \|f_{j}(\theta)\|_0  \\
&\leq & n   |\frac{4\nu_j}{\xi_j}|^2  \epsilon_{i_j}e^{- \frac{\pi}{8} |n_{i_{j+1}}|(h-\tilde{h})}  < \epsilon_{i_j},
\end{eqnarray*}
Lemma \ref{growth} imply that 
$$C^{-1} \sqrt{1+ 2\sin^2 (n\xi_j) |\frac{\nu_j}{\xi_j}|^2 }  \|B_j\|_0^{-2}  \leq   \|\mathcal{A}_n\|_0 \leq  C \sqrt{1+ 2\sin^2 (n\xi_j) |\frac{\nu_j}{\xi_j}|^2 }  \|B_j\|_0^2, $$
similar  as the first case, the result follows. 
\end{pf}

\noindent
{\bf Proof of Theorem \ref{high}:} 
We only need to consider the case $\delta(\alpha,\rho)>0$, otherwise if $\delta(\alpha,\rho)=0$, then Theorem \ref{reducibilitylocal} implies that $(\alpha,A_E)$ is reducible thus the cocycle is bounded, \eqref{upper-bound} become the trivial  bound. 

In the following, we only consider $0<\delta<\infty$, as  the remaining case is easier. First by Proposition \ref{reducibility},  $(\alpha,A_E)$ is almost reducible. More precisely, take $\tilde{h}=h-\e$,   there exist a sequence of KAM resonances $\{n_{i_j}\}_{j}$,  $B_j\in C_{\tilde{h}}^\omega(\T^d, \PSL(2,\R))$ with 
$$ \deg B_j= n_{i_{1}}+\cdots+n_{i_{j}}$$ 
such that 
$$
B_j^{-1}(\theta+\alpha)A_E(\theta)B_j(\theta)=A_je^{f_j(\theta)}=M^{-1}exp \left(
\begin{array}{ccc}
 i t^j &  \nu^j\\
\bar{ \nu^j} &  -i t^j
 \end{array}\right)Me^{f_j(\theta)}
$$
with estimate 
\begin{eqnarray}
\label{es5-31}|\nu^j| &\leq& 2 \epsilon_{i_j}^{\frac{15}{16}}e^{-2\pi|n_{i_j}|\tilde{h}},\\
\label{es5-11} \|f_j\|_{\tilde{h}} & \leq &  \epsilon_{i_j}e^{- \frac{\pi}{8} |n_{i_{j+1}}|(h-\tilde{h})},\\
\label{es5-21} |\deg B_j - n_{i_j} | &\leq&  2\epsilon_{i_{j-1}}^{\frac{1}{18\tau}} |n_{i_j}|.
 \end{eqnarray}
 
By the definition of $\delta$,  then 
$$  e^{-\eta_j|n_{i_{j}}|}:=  \|2\rho(E)-\langle n_{i_{1}}+\cdots+n_{i_{j}},\alpha\rangle\|_{\R/\Z} \geq e^{-(\delta+\frac{\e}{2}) |n_{i_{1}}+\cdots+n_{i_{j}}|}
$$
it follows from \eqref{es5-21} that 
\begin{equation}\label{dej}\eta_j\leq \delta+\e.\end{equation}
 Meanwhile,  by  \eqref{es5-11} and Lemma \ref{rol}, we have
\begin{equation*}
\|2\rho(E)-2\xi_{j}-\langle n_{i_{1}}+\cdots+n_{i_{j}},\alpha\rangle\|_{\R/\Z}\leq e^{- \frac{\pi}{16}(h-\tilde{h}) |n_{i_{j+1}}|}, 
\end{equation*}
by \eqref{dej}, this further implies that 
\begin{equation}\label{gro41}
|\xi_j|\geq  \frac{1}{2}  e^{-\eta_j|n_{i_{j}}|} \geq  \frac{1}{2}e^{-(\delta+\e)|n_{i_j}|},
\end{equation}
Now we distinguish the proof into two cases:\\

\smallskip
\noindent \textbf{Case 1: $ 2\pi \tilde{h}> \eta_j $.} In this case, \eqref{es5-31} implies that $ \left|\frac{\nu_j}{\xi_j}\right|\leq \frac{1}{2}$, then by  Lemma \ref{growth} (1),  for any $n\in I_j=[e^{ \frac{\pi}{32} |n_{i_{j}}|(h-\tilde{h})},e^{ \frac{\pi}{32} |n_{i_{j+1}}|(h-\tilde{h})}]$, we have
$$  \|(A_E)_n\|_0\leq C (\frac{ |\ln n|}{h-\tilde{h}})^{4\tau},$$ 
then  \eqref{upper-bound} follows directly. \\

\noindent \textbf{Case 2: $ 2\pi \tilde{h} \leq \eta_j $.} 
Lemma \ref{growth} (see also Remark \ref{grow-rem})  implies that for any $n\in I_j$, 
$$
 \|(A_E)_n\|_0\leq C ( \frac{|\nu_j|}{|\xi_j|}+1) (\frac{ |\ln \xi_j|}{h-\tilde{h}})^{4\tau},
$$
i.e. the maximum is obtained at $n\sim 1/\xi_j$. By \eqref{gro41}, this just means 
$$
\limsup\limits_{n\rightarrow\infty} \frac{\ln\|(A_E)_n\|_0}{\ln |n|}\leq 1-\frac{2\pi \tilde{h}}{\eta_j}+\e
$$
let $\e\rightarrow 0$, we obtain the result. \qed

\section{Explore the structure of the conjugacy}

In Section 4 and Section 5, we leverage quantitative structured almost reducibility to manage an infinite number of resonances, and we show its applications in both dynamical system and spectral theory. This section delves deeper into the quantitative aspects of structured almost reducibility, with a particular emphasis on understanding the intricate structure of the conjugacy.
 It will be crucial for us to establish the optimal quantitative Aubry duality in the next section.

The basic observation is that quantitative  structure of the conjugacy $B\in C^\omega(\T^d,$ $\PSL(2,\R))$ would imply nice estimates of its Fourier coefficients. To precise this, we introduce the 
  following goodness definition of $B$. We denote by
$$
B(\theta)=\begin{pmatrix}b_{11}(\theta)&b_{12}(\theta)\\b_{21}(\theta)&b_{22}(\theta)\end{pmatrix}\in C^\omega(\T^d,\PSL(2,\R)).
$$
Let $\hat{b}_{11}(n)$ and $\hat{b}_{12}(n)$ be the $n$-th Fourier coefficients of $b_{11}(\theta)$ and $b_{12}(\theta)$ respectively. Then we introduce the following:
\begin{Definition}\label{good}
For any $C_1>0$, $C_2>0,$ $\gamma>0$, $\ell\in\Z^d$, $B(\theta)$ is said to be $(C_1,C_
2,\gamma,\ell)$-good, if $(H1)$-$(H2)$ hold where\\
$ \mathbf{(H1)}$:\qquad \qquad \qquad $\|B\|_{0}\leq C_1$,\\
$ \mathbf{(H2)}$: \quad $|\hat{b}_{11}(n)|,|\hat{b}_{12}(n)|\leq C_2(e^{-\gamma |n+\frac{\ell}{2}|}+e^{-\gamma|n-\frac{\ell}{2}|}),$ \quad\text{for any} $n\in\Z^d$.
\end{Definition}

\begin{Remark}
We give an explanation of this definition. Actually $(H1)$ was used to control the $\ell^2$ norm of $\hat{b}_{11}(n)$ and $\hat{b}_{12}(n)$, while $(H2)$ controls the distribution of 
the localization center $\hat{b}_{11}(n)$ and $\hat{b}_{12}(n)$.
\end{Remark}

\begin{Proposition}\label{Amo}
Let $\alpha\in DC_d(\kappa,\tau)$, $V\in C^\omega(\T^d,\R)$ and $\Sigma_{V,\alpha}=\mathcal{AR}_h$. Then for any $E\in \Sigma_{V,\alpha}$ and $\e>0$, there exist sequence of $B_j\in C_{h(1-\e)}^\omega(\T^d, \PSL(2,\R))$, such that 
$$
B_j^{-1}(\theta+\alpha)A_E(\theta)B_j(\theta)=A_je^{f_j(\theta)}=M^{-1}exp \left(
\begin{array}{ccc}
 i t^j &  \nu^j\\
\bar{\nu}^j &  -i t^j
 \end{array}\right)Me^{f_j(\theta)}
$$
with estimates
\begin{itemize}
\item $ \|f_j\|_{h(1-\e)} \leq  e^{- \frac{\pi}{16}\e h |n_{i_{j+1}}|}, $
\item $|\nu_j| \leq e^{-2\pi h(1-\e)|n_{i_j}|}.$
\end{itemize}
Moreover, there exists $C=C(V,\alpha,\e)$ such that
\begin{enumerate}
\item $|\deg B_j - n_{i_j}  |< (\e h)^2 |n_{i_j}|$.
\item $\|B_j\|_r \leq  C(V,\alpha,\e)e^{(\e h)^2|n_{i_j}|} e^{\pi r  |n_{i_j}|} , \qquad  \forall r< h(1-\e)$.
\item   $B_j(\theta)$ is $( C |n_{i_j}|^{2\tau},   Ce^{(\e h )^2|n_{i_j}|}, 2\pi h(1-\e),n_{i_j})- \text{good}$.
\end{enumerate}
\end{Proposition}

\begin{pf}
For any $0< \e <\frac{1}{h}$, one can  take $\eta$ small enough such that
\begin{equation}\label{small1}
\eta  \leq  \epsilon_0\leq  D_0(\kappa,\tau,1)(\frac{\e h}{2})^{C_0\tau},
\end{equation}
where $D_0(\kappa,\tau,1)$ is the constant defined in Proposition  \ref{reducibility}.  
Same as in the proof of Corollary \ref{Amoredu}, by the assumption that $\Sigma_{V,\alpha}=\mathcal{AR}_h$ and $E\in \Sigma_{V,\alpha}$,
 there exists $\Phi_{E}\in C_{h(1-\frac{\e}{2})}^{\omega}(\T^d, \PSL(2,\R))$ such that
\begin{equation*}
\Phi_{E}(\theta+\alpha)^{-1}A_E(\theta)\Phi_{E}(\theta)=R_{\phi(E)}e^{f_E(\theta)},
\end{equation*}
with estimates
\begin{equation}\label{deg}
\|\Phi_E\|_{h}\leq \Gamma=\Gamma(V,\alpha,\eta), \qquad |\deg{\Phi_{E}}|\leq  C(\alpha) |\ln \Gamma|, \qquad \|f_E\|_{h} \leq \eta.
\end{equation}
Moreover, by \eqref{rotation number},
$$
2\rho(\alpha,R_{\phi(E)}e^{f_{E}(\theta)})=2\rho(E)-\langle \deg{\Phi_E},\alpha\rangle.
$$

Note that \eqref{small1} implies that one can apply  Proposition \ref{reducibility}, consequently 
there exists
 $B'_j\in C_{h(1-\e)}^\omega(\T^d, \PSL(2,\R))$ such that
\begin{align*}
B'_j(\theta+\alpha)^{-1}R_{\phi(E)}e^{f_E(\theta)}B'_j(\theta)&=A_je^{f_j(\theta)}=M^{-1}exp \left(
\begin{array}{ccc}
 i t^j &  \nu^j\\
\bar{\nu}^j &  -i t^j
 \end{array}\right)e^{f_j(\theta)}
\end{align*}
with estimates 
 \begin{eqnarray*}
|\nu^j| &\leq&  e^{-2\pi h(1-\e)|n_{i_j}|}, \\
\|f_j\|_{h(1-\e)} &\leq&  e^{- \frac{\pi}{16} \e h |n_{i_{j+1}}|}.
\end{eqnarray*}
Moreover, $B'_j(\theta)$ can be rewritten as  $B'_j(\theta)=\tilde{B}_j'(\theta)R_{\frac{\langle n_{i_j}, \theta\rangle}{2}}e^{Y_j'(\theta)}$ with estimate
  \begin{eqnarray}
  \label{es1-1}
 \|Y_j'\|_{h(1-\e)} &<&   e^{-2\pi  h(1-\e)|n_{i_j}|^4},\\
  \label{es3-1}
\|\tilde{B}_j'\|_{h(1-\e)}&<&C(\alpha)e^{3 \epsilon_{i_{j-1}}^{\frac{1}{18\tau}} |n_{i_j}|},\\
   \label{es4-1}
\|B_j'\|_0&<&C(\alpha)|n_{i_j}|^{2\tau},\\
\label{es5-1} 
|\deg B_j'-n_{i_j} | &\leq& 2\epsilon_{i_{j-1}}^{\frac{1}{18\tau}} |n_{i_j}|.
\end{eqnarray}

If  we let  
\begin{eqnarray*}
B_j(\theta)=\Phi_{E}(\theta)B_j'(\theta) = \Phi_{E}(\theta) \tilde{B}_j'(\theta)R_{\frac{\langle n_{i_j}, \theta\rangle}{2}}e^{Y_j'(\theta)}:= \tilde{B}_j(\theta)R_{\frac{\langle n_{i_j},\theta\rangle}{2}}e^{Y'_j(\theta)}
\end{eqnarray*}
 then one has $$
B_j^{-1}(\theta+\alpha)A_E(\theta)B_j(\theta)=A_je^{f_j(\theta)}=M^{-1}exp \left(
\begin{array}{ccc}
 i t^j &  \nu^j\\
\bar{\nu}^j &  -i t^j
 \end{array}\right)e^{f_j(\theta)}.
$$
In the following, we give the estimates of  the conjugacy $B_j$. First note \eqref{small1} implies that 
\begin{equation}\label{eee3-1} 2\epsilon_{i_{j-1}}^{\frac{1}{18\tau}} < 2 \epsilon_{0}^{\frac{1}{18\tau}} \leq \frac{1}{2}(\e h)^2.
\end{equation}
By \eqref{deg},    \eqref{es3-1} and \eqref{eee3-1}, we have
\begin{align}\label{tb}
 \| \tilde{B}_j \|_{h(1-\e)} \leq \|\Phi_{E}\|_{h(1-\e)}\|\tilde{B}'_j\|_{h(1-\e)} 
 \leq C(V,\alpha,\e)e^{(\e h)^2|n_{i_j}|}.
\end{align}
Consequently, (1) follows from \eqref{deg}, \eqref{es5-1}, and \eqref{eee3-1}. While $(2)$ and $(H1)$ part of $(3)$ follows from \eqref{deg},  \eqref{es4-1} and \eqref{eee3-1}.


We are left to prove $(H2)$ part of $(3)$. To prove this,  rewrite  $B_j(\theta)$ as
\begin{align*}
B_j(\theta)=\tilde{B}_j(\theta) R_{\frac{\langle n_{i_j},\theta\rangle}{2}}e^{Y'_j(\theta)}=\tilde{B}_j(\theta)M^{-1}MR_{\frac{\langle n_{i_j},\theta\rangle}{2}}M^{-1}e^{M Y'_j(\theta)M^{-1}}M,
\end{align*}
 if we denote
\begin{equation*}
e^{M^{-1}Y'_j(\theta)M}=\begin{pmatrix}y_{11}(\theta)&y_{12}(\theta)\\y_{21}(\theta)&y_{22}(\theta)\end{pmatrix}, 
 \qquad 
\tilde{B}_j(\theta)M^{-1}=\begin{pmatrix}\tilde{b}_{11}(\theta)&\tilde{b}_{12}(\theta)\\\tilde{b}_{21}(\theta)&\tilde{b}_{22}(\theta)\end{pmatrix},
\end{equation*}
then direct computation shows that
\begin{align*}
b_{11}(\theta)
&= \frac{1}{2 i}( y_{11}(\theta)+ y_{12}(\theta) ) \tilde{b}_{11}(\theta)e^{2\pi i\frac{\langle n_{i_j},\theta\rangle}{2}}\\ 
&+ \frac{1}{2 i}(y_{21}(\theta)+y_{22}(\theta))\tilde{b}_{12}(\theta)e^{-2\pi i\frac{\langle n_{i_j},\theta\rangle}{2}}.
\end{align*}
By \eqref{es1-1} and \eqref{tb}, we have
$$
|\hat{b}_{11}(n)|\leq Ce^{(\e h)^2|n_{i_j}|}(e^{-2\pi h(1-\e) |n+\frac{n_{i_j}}{2}|}+e^{- 2\pi h(1-\e) |n-\frac{n_{i_j}}{2}|}),
$$
similarly, one obtain estimates of $|\hat{b}_{12}(n)|$ and  finish the whole proof of (3).
\end{pf}

Note in section \ref{kam-res},  we initially demonstrated the link between KAM resonances and rotation number resonances within the perturbative regime. However, our analysis actually extends to the global almost reducibility regime, revealing the connection between these two phenomena that persists even in the  non-perturbative regime.
 Similarly, We order the $\e_0$-resonances of $\rho(E)$ by $0=|\ell_0|<|\ell_1|\leq |\ell_2|\leq \cdots$.
 \begin{Corollary}\label{relation2}
Under the assumption of Proposition  \ref{Amo}, for any $\e_0>0$, there exist $B_j\in C_{\tilde{h}}^\omega(\T^d,\PSL(2,\R))$,   such that 
\begin{equation}\label{conjugacy}
B_j^{-1}(\theta+\alpha)A_E(\theta)B_j(\theta)=A_je^{f_j(\theta)}=M^{-1}exp \left(
\begin{array}{ccc}
 i t^j &  \nu^j\\
\bar{\nu}^j &  -i t^j
 \end{array}\right)Me^{f_j(\theta)}
\end{equation}
with $\ell_i=\deg B_j$ provided $\ell_i$ is large enough (the choice of $j$ depends on $\ell_i$). 
\end{Corollary}
\begin{pf}
Note that by Proposition \ref{Amo}, we have
\begin{equation*}
\Phi_{E}(\theta+\alpha)^{-1}A_E(\theta)\Phi_{E}(\theta)=R_{\phi(E)}e^{f_E(\theta)},
\end{equation*}
and thus
$$
2\rho(\alpha,R_{\phi(E)}e^{f_{E}(\theta)})=2\rho(E)-\langle \deg{\Phi_E},\alpha\rangle.
$$
where by \eqref{deg}, $\deg{\Phi_E}$ is uniformly bounded for $E\in \Sigma_{V,\alpha}$.
Hence for $i$ sufficiently large,   
$$
\cdots<|\ell_i-\deg{\Phi_E}|\leq |\ell_{i+1}-\deg{\Phi_E}|\leq \cdots
$$
are the $\e'_0$-resonances of $\rho(\alpha,R_{\phi(E)}e^{f_{E}(\theta)})$ for some $\e_0'\leq\e_0$. 

By Lemma 3.1, there exists $B_j'$ with $\deg B_j'=\ell_i-\deg{\Phi_E}$ such that
\begin{align*}
B'_j(\theta+\alpha)^{-1}R_{\phi(E)}e^{f_E(\theta)}B'_j(\theta)&=A_je^{f_j(\theta)}
\end{align*}
Let $B_j=\Phi_E B_j'$, the desired result follows.
\end{pf}

\section{Sharp version of quantitative Aubry duality}\label{s1.7}

In the previous section, we established structured quantitative almost reducibility, which provided tight upper bounds on the off-diagonal entries $\nu_j$ within the constant matrix. However, this approach lacked a complementary lower bound analysis, leaving the optimality of these estimates uncertain. This section introduces a novel, sharp version of quantitative Aubry duality, which serves as a crucial stepping stone to address this gap.

The key insight lies in the interplay between the smallness of $\nu_j$ and the structure of the conjugacy matrix $B_j$, particularly its goodness properties (characterized by $(C_1,C_2,\gamma,k)$). When $\nu_j$ is small, this implies the dual almost Mathieu operator $H_{\lambda^{-1},\alpha,0}$ exhibits two approximately linearly independent localized eigenfunctions, centered around $\pm \frac{k}{2}$. However, this observation contradicts the fundamental principle of the Wronskian argument, which prohibits such a configuration.

By leveraging this contradiction, we not only demonstrate the optimality of the estimates derived from quantitative almost reducibility but also open the door to a range of optimal spectral applications.

%
%
%

\begin{Proposition}\label{lowred}
Let $\alpha\in\R\backslash\Q$, $0<|\lambda|<1$, $0<\e h_{\lambda}<\frac{1}{1000}$ and $E\in \Sigma_{\lambda,\alpha}$. 
Then  \begin{eqnarray}\label{estimat1}
B(\theta+\alpha)^{-1}S_{E}^{\lambda}(\theta)B(\theta)=Id
+F(\theta)
\end{eqnarray}
has  no $B\in C_{\frac{h_{\lambda}}{2\pi}(1-\e)}^{\omega}(\T, \PSL(2,\R))$ 
with 
\begin{enumerate}
\item $B(\theta)$ is $( C(\lambda,\alpha,\e) |k|^{2\tau},   C(\lambda,\alpha,\e) e^{(\e h_{\lambda} )^2|k|}, h_{\lambda}(1-\e),k)- \text{good}$,
\item $\|F\|_{\frac{h_{\lambda}}{2\pi}(1-\e)}\leq e^{-h_{\lambda}(1+\e)|k|}. $
\end{enumerate} if $ |k| \geq K(\lambda,\alpha,\e)$ is sufficiently large. 
\end{Proposition}
\pf 
We will argue by contradiction, i.e. assume $(1)$ and $(2)$ hold. 
Without loss of generality, we assume $k>0$, and denote
$$
B(\theta)=\begin{pmatrix}b_{11}(\theta)&b_{12}(\theta)\\ b_{21}(\theta)&b_{22}(\theta)\end{pmatrix}= \left( \overrightarrow{u}(\theta),  \overrightarrow{v}(\theta)\right),\ \ F(\theta)=\begin{pmatrix}
\beta_{1}(\theta) & \beta_{2}(\theta)\\
\beta_{3}(\theta)& \beta_{4}(\theta)
 \end{pmatrix},
$$
then  \eqref{estimat1} implies
 \begin{eqnarray}
\label{new2-1} &&b_{11}(\theta)=b_{21}(\theta+\alpha)+b_{21}(\theta+\alpha)\beta_{1}(\theta)+b_{22}(\theta+\alpha)\beta_{3}(\theta),\\
\label{new2-2}
&&b_{12}(\theta)=b_{21}(\theta+\alpha)\beta_{2}(\theta)+b_{22}(\theta+\alpha)\beta_{4}(\theta)+b_{22}(\theta+\alpha).
 \end{eqnarray}

As the perturbation $F(\theta)$ may not be zero, the Fourier coefficients of  $b_{11}(\theta)$ and $b_{12}(\theta)$, denoted by 
$\hat{b}_{11}(n)$ and $\hat{b}_{12}(n)$, are not
necessarily to be the solutions of the dual almost Mathieu operator $H_{\lambda^{-1},\alpha,0}$:
\begin{equation}\label{mathieu-1}
\widehat{u}(k+1)+\widehat{u}(k-1) +2\lambda^{-1} \cos(
k\alpha)\widehat{u}(k)=\lambda E\widehat{u}(k).
\end{equation}
However, $\hat{b}_{11}(n)$ and $\hat{b}_{12}(n)$ are two approximate solutions of $(\ref{mathieu-1})$ if $F(\theta)$ is sufficiently small. 
To see this, by \eqref{estimat1}, \eqref{new2-1} and \eqref{new2-2}, we have
\begin{eqnarray}\label{block-red1} (E-2\lambda
\cos2\pi\theta)b_{11}(\theta)-b_{11}(\theta-\alpha)-b_{11}(\theta+\alpha)=f(\theta),\end{eqnarray}
\begin{eqnarray}\label{block-red2} (E-2\lambda
\cos2\pi\theta) b_{12}(\theta)-b_{12}(\theta-\alpha)-b_{12}(\theta+\alpha)=g(\theta),\end{eqnarray}
where \begin{eqnarray*} f(\theta)=
b_{11}(\theta+\alpha)\beta_{1}(\theta)+
b_{12}(\theta+\alpha)\beta_{3}(\theta)-
b_{21}(\theta)\beta_{1}(\theta-\alpha)-b_{22}(\theta)\beta_{3}(\theta-\alpha),
\end{eqnarray*}
\begin{eqnarray*} g(\theta)=
b_{11}(\theta+\alpha)\beta_{2}(\theta)+
b_{12}(\theta+\alpha)\beta_{4}(\theta)
-b_{21}(\theta)\beta_{2}(\theta-\alpha)-b_{22}(\theta)\beta_{4}(\theta-\alpha),
\end{eqnarray*}
and the assumption $(1)$ and $(2)$ imply the estimates
\begin{equation}\label{fg} \|f\|_0, \|g\|_0 \leq 4\|B\|_0 \|F\|_0 \leq 4C |k|^{2\tau} e^{-h_{\lambda}(1+\e)|k|} \leq   e^{-h_{\lambda}(1+\frac{4}{5}\e)|k|}. \end{equation}

Setting $\A(\theta)= S_{\lambda^{-1}E}^{\lambda^{-1}}(\theta),$ going to Fourier coefficients, 
it immediately follows from \eqref{block-red1} and \eqref{block-red2} that 
 $\hat{b}_{11}(n)$ and $\hat{b}_{12}(n)$ satisfy
\begin{equation}\label{new4}
\begin{pmatrix}
\hat{b}_{11}(n+1)\\
\hat{b}_{11}(n)\\
\end{pmatrix}= \A(n\alpha) \begin{pmatrix}
\hat{b}_{11}(n)\\
\hat{b}_{11}(n-1)\\
\end{pmatrix}-\begin{pmatrix}\lambda^{-1} \hat{f}(n)\\0\end{pmatrix},
\end{equation}
\begin{equation}\label{new4-1}
\begin{pmatrix}
\hat{b}_{12}(n+1)\\
\hat{b}_{12}(n)\\
\end{pmatrix}=\A(n\alpha)\begin{pmatrix}
\hat{b}_{12}(n)\\
\hat{b}_{12}(n-1)\\
\end{pmatrix}-\begin{pmatrix}\lambda^{-1} \hat{g}(n)\\0\end{pmatrix}.
\end{equation}

The first observation is that assumption $(1)$ implies that  $\hat{b}_{11}(n)$ has two localization center $\pm \frac{k}{2}$, consequently $\hat{b}_{11}(n)$ has relatively large coefficients around   $\pm \frac{k}{2}$. Denote
$$
I_k=\left[\left(1-\frac{\e}{10}\right)\frac{k}{2},\left(1+\frac{\e}{10}\right)\frac{k}{2}\right],
$$
Then we have the following:
\begin{Lemma}\label{initial}
There exists $k_0\in I_k  \cup- I_k  $ and $K_0(\lambda,\alpha,\e)>0$, such that
\begin{equation}\label{z1-estimate-2}
|\hat{b}_{11}(k_0)|\geq \frac{1}{ |k|^{6\tau}} ,
\end{equation}
provided $k>K_0$.
\end{Lemma}
\begin{pf} 
On one hand,  note that 
$\|\overrightarrow{u}\|_{L^2}\|\overrightarrow{v}\|_{L^2}>1$ since $\det B(\theta)=1.$ It follows  that
$$\|\overrightarrow{u}\|_{L^2}> \frac{1}{\|\overrightarrow{v}\|_{L^2}}>\frac{1}{\|B\|_{C^0}}.$$
By \eqref{new2-1}, one has
$$
2\|b_{11}\|_{L^2}+ 2\|B\|^2_0\|F\|_0 \geq \|b_{11}\|_{L^2}+\|b_{21}\|_{L^2}= \|\overrightarrow{u}\|_{L^2} \geq\|B\|_{C^0}^{-1},
$$
consequently by  \eqref{fg} and $(H1)$ part of the assumption $(1)$, we have 
\begin{eqnarray}\label{b11}
\|\hat{b}_{11}\|_{\ell^2}=\|b_{11}\|_{L^2}&\geq&(2\|B\|_{C^0})^{-1}- \|B\|^2_0\|F\|_0 \\
\nonumber &\geq&  2^{-1}C^{-1}  |k|^{-2\tau}- C^2 |k|^{4\tau} e^{-h_{\lambda}(1+\e)|k|}\geq 4^{-1}C^{-1}  |k|^{-2\tau}.
\end{eqnarray}

On the other hand, by  $(H2)$ part of the assumption $(1)$, we have 
$$
|\hat{b}_{11}(n)|\leq C(\lambda,\alpha,\e)e^{(\e h_{\lambda} )^2|k|}\left(e^{-h_{\lambda}(1-\e) |n+\frac{k}{2}|}+e^{-h_{\lambda}(1-\e)|n-\frac{k}{2}|}\right),$$
which implies that
\begin{align*}
\sum\limits_{n\notin I_k \cup - I_k}|\hat{b}_{11}(n)|^2\leq 4Ce^{2(\e h_{\lambda} )^2|k|}e^{-h_{\lambda}(1-\e)\frac{\e}{10}|k|}
\leq e^{-10(\e h_{\lambda} )^2|k|} .
\end{align*}

By \eqref{b11} and the fact $|I_k\cup- I_k|\leq k$, it follows that there exists $k_0\in I_k\cup -I_k$, such that
\begin{align*}
k|\hat{b}_{11}(k_0)|^2 \geq \sum\limits_{n\in I_k\cup- I_k}|\hat{b}_{11}(n)|^2 &=\|b_{11}\|_{\ell^2}^2-\sum\limits_{n\notin I_k\cup- I_k}|\hat{b}_{11}(n)|^2\\
 &\geq  16^{-1}C^{-2} |k|^{-4\tau}-e^{-10(\e h_{\lambda} )^2|k|},
\end{align*}
then the result follows. 
\end{pf}

In the following, $k_0$ is chosen so that \eqref{z1-estimate-2} is satisfied, and denote
$$D_{k_0}=\det{\begin{pmatrix}\hat{b}_{11}(k_0)&\hat{b}_{12}(k_0)\\ \hat{b}_{11}(k_0+1)&\hat{b}_{12}(k_0+1)\end{pmatrix}}.$$
We will show actually $(1)$ and $(2)$ imply that 
 vectors $(\hat{b}_{12}(n))_n$ and $(\hat{b}_{11}(n))_n$ are almost linearly independent.

\begin{Lemma}\label{key1}
 There exists $K_1(\lambda,\alpha,\e)>0$ such that if  $k>K_1$, then 
$$
|D_{k_0}|\geq e^{-h_{\lambda}(1+\frac{\e}{5})k}.
$$
\end{Lemma}
\begin{pf}
We prove this by contradiction. If $|D_{k_0}|\leq e^{-h_{\lambda}(1+\frac{\e}{5})k}$, then by the fact that
\begin{align*}
 \left(
\begin{array}{ccc} \hat{b}_{12}(k_0)\\ \hat{b}_{12}(k_0+1)  \end{array}
\right) = C_{k_0}  \left(
\begin{array}{ccc} \hat{b}_{11}(k_0)\\ \hat{b}_{11}(k_0+1) \end{array}
\right) 
+\frac{D_{k_0}}{\sqrt{|\hat{b}_{11}(k_0)|^2+|\hat{b}_{11}(k_0+1)|^2}}  \left(
\begin{array}{ccc} -\overline{\hat{b}_{11}(k_0+1)} \\ \overline{\hat{b}_{11}(k_0)} \end{array}
\right) ,
\end{align*}
one sees that
\begin{equation}\label{error}
\left|\begin{pmatrix}\hat{b}_{12}(k_0)-C_{k_0}\hat{b}_{11}(k_0)\\ \hat{b}_{12}(k_0+1)-C_{k_0}\hat{b}_{11}(k_0+1)\end{pmatrix}\right|\leq |D_{k_0}|<e^{-h_{\lambda}(1+\frac{\e}{5})k},
\end{equation}
which means  the orthogonal projection of $(\hat{b}_{12}(k_0),\hat{b}_{12}(k_0+1))^{T}$ to the vector $(\hat{b}_{11}(k_0),$ $\hat{b}_{11}(k_0+1))^{T}$ is small. 
Since 
$$
C_{k_0}=\frac{\hat{b}_{12}(k_0)+D_{k_0}\overline{\hat{b}_{11}(k_0+1)}\diagup\sqrt{|\hat{b}_{11}(k_0)|^2+|\hat{b}_{11}(k_0+1)|^2}}{\hat{b}_{11}(k_0)},
$$
  by Lemma \ref{initial} and  $(H1)$ part of the assumption $(1)$, we have 
\begin{equation}\label{gne5}
|C_{k_0}|\leq C  |k|^{12\tau}e^{(\e h_{\lambda})^2k}  \leq e^{4(\e h_{\lambda})^2k},
\end{equation}
provided $k$ is sufficiently large depending on $\lambda,\alpha,\e$.

In the following, we consider
\begin{align*}
\widehat{B}(\theta)= \begin{pmatrix}b_{11}(\theta)&b_{12}(\theta)-C_{k_0} b_{11}(\theta)\\b_{11}(\theta-\alpha)&b_{12}(\theta-\alpha)-C_{k_0} b_{11}(\theta-\alpha)\end{pmatrix}\end{align*}
and aim to estimate its determinant.   First, 
we will show as a consequence of \eqref{error},  $b_{12}(\theta)-C_{k_0} b_{11}(\theta)$ is small. 
To show this,  we have the following:

\begin{Claim} If $|k|\geq K_1(\lambda,\alpha,\e)$ is large enough, then we have:
\begin{enumerate}
\item If $n\geq (\frac{1}{2}+\frac{\e}{20})k$ or $n\leq (-\frac{1}{2}-\frac{\e}{20})k$,  then 
 $$|\hat{b}_{12}(n)-C_{k_0}\hat{b}_{11}(n)| \leq  e^{-\frac{\e h_{\lambda}}{25}|k|}e^{-h_{\lambda}(1-\e)|n-\frac{k}{2}-\frac{\e k}{20}|}.$$
\item If $(-\frac{1}{2}-\frac{\e}{20})k\leq n\leq (\frac{1}{2}+\frac{\e}{20})k$, then 
$$|\hat{b}_{12}(n)-C_{k_0}\hat{b}_{11}(n)| \leq  e^{-\frac{\e h_{\lambda}}{25}|k|} .$$
\end{enumerate}
\end{Claim}
\begin{pf}
By \eqref{gne5} and $(H2)$ part of the assumption $(1)$, we have
\begin{align*}
 |b_{12}(n)|+C_{k_0}|b_{11}(n)|\leq Ce^{6(\e h_{\lambda} )^2|k|}e^{-h_{\lambda}(1-\e)\frac{\e}{10}\frac{|k|}{2}}e^{-h_{\lambda}(1-\e)|n-\frac{k}{2}-\frac{\e k}{20}|},
\end{align*}
the the first item follows.

To prove the second item,  denote by
$$
\tilde{p}_n=\begin{pmatrix}-\lambda^{-1}\hat{g}(n)+\lambda^{-1} C_{k_0}\hat{f}(n)\\0\end{pmatrix},\ \ \tilde{y}_n=\begin{pmatrix}\hat{b}_{12}(n+1)-C_{k_0}\hat{b}_{11}(n+1)\\ \hat{b}_{12}(n)-C_{k_0}\hat{b}_{11}(n)\end{pmatrix}.
$$
Then as a result of \eqref{new4} and \eqref{new4-1}, we have 
$$\tilde{y}_n=\A(n\alpha)\tilde{y}_{n-1}+\tilde{p}_n,$$
which implies that 
$$\tilde{y}_n=\A^{n-k_0}((k_0+1)\alpha)\tilde{y}_{k_0}+\sum\limits_{j=k_0+1}^n\A^{n-j}((j+1)\alpha)\tilde{p}_j$$
where $\A^n(\theta)= \A(\theta+(n-1)\alpha)\cdots \A(\theta).$  

By the assumption and  \eqref{fg}, we have 
\begin{eqnarray*}
 |\tilde{p}_n|  \leq e^{-h_{\lambda}(1+\frac{4}{5}\e)|k|},\qquad 
  |\tilde{y}_{k_0}|=|D_{k_0}|\leq e^{-h_{\lambda}(1+\frac{\e}{5})k},
\end{eqnarray*}
while by Lemma \ref{lem:upperbdd}, we have for any $\theta\in\T$, and for any $n\in \Z$,
\begin{equation}\label{ub}
\|\A^n(\theta)\|\leq C(\lambda,\alpha,\e)e^{h_{\lambda}(1+\frac{\e}{100})|n|}.
\end{equation}
As a result, if $k$ is sufficiently large depending on $\lambda,\alpha,\e$, then one can estimate
\begin{eqnarray*}
 |\tilde{y}_n|&\leq& Ce^{h_{\lambda}(1+\frac{\e}{100})|n-k_0|} \left( e^{-h_{\lambda}(1+\frac{\e}{5})|k|}+ |n-k_0| e^{-h_{\lambda}(1+\frac{4\e}{5})|k|} \right)\\
&\leq& e^{h_{\lambda}(1+\frac{4\e}{25})|k|} e^{-h_{\lambda}(1+\frac{\e}{5})|k|} \leq e^{-\frac{\e h_{\lambda}}{25}|k|},
\end{eqnarray*}
where the second inequality follows from the  choice of $k_0$ that 
$|n-k_0|\leq (1+\frac{\e}{10})k.$

\end{pf}

This claim immediate imply that 
$
\left\|b_{12}(\cdot)-C_{k_0} b_{11}(\cdot)\right\|_0\leq 2e^{-\frac{\e h_{\lambda}}{25}|k|},
$
as a consequence,
\begin{align*}
\det   \widehat{B} \leq  2C e^{(\e h_{\lambda})^2k}e^{-\frac{\e h_{\lambda}}{25}|k|}\leq \frac{1}{100}.
\end{align*}
However by \eqref{new2-1} and \eqref{new2-2}, we have 
\begin{align*}
\det   \widehat{B}
\geq  \det B- 4 \|B\|^2_0 \|F\|_0
\geq \det  B- 4 C(\lambda,\alpha,\e)^2  |k|^{4\tau}e^{-h_{\lambda}(1+\e)|k|}>\frac{1}{2},
\end{align*}
this is a contradiction.
\end{pf}

Let $(u_1(n))$, $(u_2(n))$ be two solutions of the almost Mathieu operator, with the initial datum $u_1(k_0)=b_{11}(k_0)$, $u_1(k_0+1)=b_{11}(k_0+1)$ and $u_2(k_0)=b_{12}(k_0)$, $u_2(k_0+1)=b_{12}(k_0+1)$. Denote
$$
D'_n:=\det{\begin{pmatrix}u_1(n)&u_2(n)\\u_1(n+1)&u_2(n+1)\end{pmatrix}}.
$$
Then, we have the following estimates:
\begin{Lemma}\label{key2}
 There exists $K_2(\lambda,\alpha,\e)>0$ such that if  $k>K_2$, then 
$$
|D'_{(1+\frac{4}{5}\e)k}|\leq 4e^{-h_{\lambda}(1+\frac{1}{4}\e)k}, \ \ {\rm if} \ \ k_0\in I_k,
$$
$$
|D'_{-(1+\frac{4}{5}\e)k}|\leq 4e^{-h_{\lambda}(1+\frac{1}{4}\e)k}, \ \ {\rm if}\ \  k_0\in -I_k.
$$
\end{Lemma}
\begin{pf}  
We only prove the case $k_0\in I_k$, the other case follows exactly the same way. First  by  $(H2)$ part of the assumption $(1)$, we have 
$$
|\hat{b}_{11}(n)|, |\hat{b}_{12}(n)|\leq C(\lambda,\alpha,\e)e^{(\e h_{\lambda} )^2|k|}\left(e^{-h_{\lambda}(1-\e) |n+\frac{k}{2}|}+e^{-h_{\lambda}(1-\e)|n-\frac{k}{2}|}\right),$$
which implies that
\begin{equation}\label{decay}
|\hat{b}_{11}(k+\frac{4}{5}\e k)|, |\hat{b}_{12}(k+\frac{4}{5}\e k)|\leq 2C(\lambda,\alpha,\e)e^{(\e h_{\lambda} )^2|k|}e^{-h_{\lambda}(1+\frac{1}{2}\e)\frac{|k|}{2}} \leq e^{-h_{\lambda}(1+\frac{2}{5}\e)\frac{|k|}{2}}, \end{equation}
without loss of generality,  we only need to prove that 
\begin{equation}\label{decay-1}
|\hat{u}_2(k+\frac{4}{5}\e k)|\leq 2e^{-h_{\lambda}(1+\frac{2}{5}\e)\frac{|k|}{2}}, 
\end{equation}

Denote by
$$
p_n=\begin{pmatrix}-\lambda^{-1}\hat{g}(n)\\0\end{pmatrix}, \ \ y_n=\begin{pmatrix}\hat{u}_2(n+1)-\hat{b}_{12}(n+1)\\ \hat{u}_2(n)-\hat{b}_{12}(n)\end{pmatrix},
$$
 then by the similar argument as in Lemma
\ref{key1}, we know
$$
y_n=\A^{n-k_0}((k_0+1)\alpha)y_{k_0}+\sum\limits_{j=k_0+1}^n \A^{n-j}((j+1)\alpha)p_j,
$$
As the initial datum of  $\hat{u}(k_0)$ and $\hat{b}_{12}(k_0)$
 are equal, i.e. $y_{k_0}=0$, then by \eqref{fg} and \eqref{ub}, we have
\begin{align*}
|y_{(1+\frac{4}{5}\e)k}|\leq 2|k|e^{h_{\lambda}(1+\frac{1}{100}\e)(\frac{1}{2}+\frac{17}{20}\e)|k|}Ce^{(\e h_{\lambda})^2|k|}e^{-h_{\lambda}(1+\e)|k|}
\leq e^{-h_{\lambda}(1+\frac{1}{4}\e)\frac{|k|}{2}}
\end{align*}
provided $k$ is sufficiently large depending on $\lambda,\alpha,\e$. Once we have this, then \eqref{decay-1} follows from \eqref{decay}, and we finish the proof.
\end{pf}

We can now finish the  proof  of Proposition \ref{lowred}. 
Since $(\hat{u}_1(n))$,
$(\hat{u}_2(n))$  are  two linearly independent  solutions of
the almost Mathieu operator $H_{\lambda^{-1},\alpha,0}$, then by
Liouville's theorem, we have
$$
D'_{-(1+\frac{4}{5}\e)k}=D'_{(1+\frac{4}{5}\e)k}=D'_{k_0}=D_{k_0},
$$
however  by  Lemma \ref{key1} and Lemma \ref{key2}, this is a contradiction. \qed


%
%
%

\section{Sharp  spectral applications for  AMO}\label{s1.8}

This section exclusively focuses on the applications of spectral theory of AMO, which constitutes the primary impetus for the development of our structured quantitative  almost reducibility. Key findings encompass estimates on the size of spectral gaps, asymptotic properties of generalized eigenfunctions, and a novel proof of the arithmetic transition conjecture in phase space.

\subsection{Sharp estimates on the length of spectral gaps:}

It was first observed by Moser-P\"oschel \cite{mp} that to estimate the size of each spectral gap, we need to analyze the behavior of Schr\"odinger cocycles at the edge points of the spectral gaps. More precisely, Moser and P\"oschel \cite{mp} showed for an arithmetic defined subset $\mathcal{K}\subset \Z$, if we are at the k-th edge point with $k\in\mathcal{K}$, saying $E_k$ in the spectrum with $N(E_k)=k\alpha(\mod\Z)$, the cocycle is reducible to a  constant parabolic or elliptic cocycle.  i.e., there is $B\in C_{h_*}^\omega(\T,\PSL(2,R))$ such that
\begin{equation}\label{conj2}
B_k^{-1}(x+\alpha)A_{E_k}(x)B_k(x)=\begin{pmatrix}1&c_k\\ 0&1\end{pmatrix}.
\end{equation}
Moreover, the conjugacy $B_k$ is close to identity, consequently 
$$ \frac{|c_k|}{4} \leq  |G_k(V)| \leq 2|c_k|, \quad \forall k\in\mathcal{K} .$$
In general, we have the criterion on the spectral gaps:
\begin{Proposition}[\cite{LYZZ}]\label{criterion}
Let $\alpha\in DC_d(\kappa,\tau)$, $\e\in (0,\frac{1}{4})$ and $V\in C^\omega(\T^d,\R)$ be a non-constant function. Let $E$ be an edge point of the spectral gap $G(V)$. Assume that there are $c \in (0,\frac{1}{2})$ and $B\in C^\omega_r(\T^d,PSL(2,\R))$ such that
$$
B(x+\alpha)^{-1}A_E(x)B(x)=\begin{pmatrix}1&c\\0&1\end{pmatrix}.
$$ 
with estimate
$$
\|B\|_r^{14}c^\e\leq C(d,\kappa,\tau)r^{4(4\tau+1)},
$$
where $C(d,\kappa,\tau)$ is numerical constant,
then we have 
$$ c^{1+\e}  \leq  |G(V)| \leq c^{1-\e}.$$
\end{Proposition}
The  key observation is from Avila-You-Zhou' s solution of non-critical Dry Ten Martini Problem \cite{ayz2}, more precisely, for subcritical AMO, while \eqref{conj2} is true for any $k\in\Z$ (provided by Eliasson's theory \cite{eli}), a lower bound estimate for $c_k$ is indeed achievable. With the aid of \textit{sharp  quantitative Aubry duality}, the precise decay rate of $c_k$ can now be determined.

\begin{Corollary}\label{ut}
Let $\alpha\in DC$, $0<|\lambda|<1$ and  $E\in \Sigma_{\lambda,\alpha}$ with $2\rho(E)=k\alpha\mod \Z$ for $k\in\Z\backslash\{0\}$. Given any $0<\e h_\lambda <\frac{1}{1000}$, if $|k| \geq K(\lambda,\alpha,\e)$ is large enough, there exist $B\in C_{\frac{h_{\lambda}}{2\pi}(1-\e)}^{\omega}(\T,\PSL(2,\R))$ and $c_k\in\R$, such that
$$
B^{-1}(\theta+\alpha)S_E^{\lambda}(\theta)B(\theta)=\begin{pmatrix}1& c_k\\
0&1
\end{pmatrix}.
$$
with estimate
\begin{enumerate}
\item  $\|B\|_{r}\leq C(\lambda,\alpha,\e)e^{2(\e h_{\lambda})^2|k|} e^{2r|k|}, \quad \forall r< \frac{h_{\lambda}}{2\pi}(1-\e)$,
\item  $ e^{-h_{\lambda}(1+2\e)|k|} \leq |c|\leq e^{-h_{\lambda}(1-2\e)|k|}$.
\end{enumerate}
\end{Corollary}

\begin{pf}
By Avila's global theory \cite{avila2015global}, if $|\lambda|<1$, $(\alpha,S_E^\lambda)$ is subcritical in the strip $\{|\Im z|< \frac{h_{\lambda}}{2\pi}\}$, and  for any $\e>0$, $(\alpha,S_E^{\lambda})$ is almost reducible in the strip $\{|\Im z|<\frac{h_{\lambda}}{2\pi}(1- \frac{\e}{2})\}$ \cite{avila,avila2}.
If we assume $2\rho(E)=k\alpha\mod \Z$, then it is well-known there are at most finitely many  KAM resonances \cite{eli,LYZZ}, denote the last resonance by $n_{i_k}$ and $n_{i_{k}+1}=\infty$. Furthermore, 
by Proposition \ref{Amo},   there exists 
 $B_{k}\in C_{\frac{h_{\lambda}}{2\pi}(1-\e)}^\omega(\T, \PSL(2,\R))$ which is $(C |n_{i_k}|^{2\tau},   Ce^{(\e h_{\lambda} )^2|n_{i_k}|}, h_{\lambda}(1-\e),n_{i_k})- \text{good}$, such that 
$$
B_{k}^{-1}(\theta+\alpha)S_E^{\lambda}(\theta)B_{k}(\theta)=A= M^{-1}exp \left(
\begin{array}{ccc}
 i t &  \nu \\
\bar{\nu} &  -i t
 \end{array}\right)M
$$
with estimate
\begin{eqnarray}
\label{nu+}|\nu| &\leq&  e^{-h_{\lambda}(1-\e)|n_{i_k}|}.
\end{eqnarray}
Moreover, by noting $\deg B_{k}=k$, then we have estimate
\begin{eqnarray}
\label{kdiff} |k - n_{i_k}  | &<& (\e h_{\lambda})^2 |n_{i_k}|, \\
\label{bj} \|B_{k}\|_r &\leq &  C(\lambda,\alpha,\e) e^{2(\e h_{\lambda})^2|n_{i_k}|} e^{2\pi r  |n_{k}|}, \qquad  \forall r< \frac{h_{\lambda}}{2\pi}(1-\e).
 \end{eqnarray}

If furthermore $E\in \Sigma_{\lambda,\alpha}$, then $A$ is parabolic, which implies that $t= |\nu|$. Thus if $k$ is large enough (i.e. $n_{i_k}$ is large), then by Proposition \ref{lowred},  we have 
\begin{equation}\label{lbdt} t= |\nu| > e^{-h_{\lambda}(1+\e)|n_{i_k}|}.\end{equation}
Meanwhile, as $A$ is parabolic, there exists $P\in SO(2,\R)$ such that 
$$  P^{-1}AP=P^{-1} M^{-1}exp \left(\begin{array}{ccc}
 i t &  \nu \\
\bar{\nu} &  -i t
 \end{array}\right) MP=  \left(
\begin{array}{ccc}
 1 &  2t \\
0 &  1
 \end{array}\right).  $$
Then $B=B_{k}P$ satisfy our needs, while $(1)$ follows from \eqref{kdiff} and \eqref{bj},  $(2)$ follows from  \eqref{nu+},  \eqref{kdiff}  and \eqref{lbdt}.

\end{pf}

\textbf{Proof of Theorem \ref{optimal gap}:} Choose $r=\frac{\e h_{\lambda}}{300}$, then by Corollary \ref{ut}, if $|k| \geq K(\lambda,\alpha,\e)$ is large enough, then we have 
$$
\|B\|_r^{14}c_k^\e\leq   C(\lambda,\alpha,\e)^{15} e^{\frac{\e h_{\lambda}}{10} |k|}  e^{-\e h_{\lambda}(1-2\e)|k|} \leq  C(\alpha) r^{4(4\tau+1)},
$$
then by Proposition  \ref{criterion} and (2) of Corollary \ref{ut}, we have 
$$
c(\lambda,\alpha,\e)e^{-h_{\lambda}(1+\e)^2|k|}  \leq c_k^{1+\e}\leq G_k(\lambda)\leq c_k^{1-\e} \leq C(\lambda,\alpha,\e)e^{-h_{\lambda}(1-\e)^2|k|}.
$$
Finally, let $\e\rightarrow 0$, we then finish the proof of Theorem \ref{optimal gap}. \qed

\subsection{Optimal  asymptotical growth}

For any $\e_0>0$, let $\{\ell_j\}_{j=1}^\infty$ be a sequence such that
\begin{equation}\label{aaaa}
\|2\rho(E)-\ell_j\alpha\|_{\R/\Z}\leq e^{-\e_0|\ell_j|}.
\end{equation}
By Proposition \ref{Amo} and Corollary \ref{relation2},  for any fixed $j\in \N$ large enough, we have
$$
\ell_j=n_{i_0}+n_{i_1}+\cdots+n_{i_j},\ \ \ell_{j+1}=\ell_j+n_{i_{j+1}}+\cdots+n_{i_{j+m(j)}}.
$$
where $(n_{i_j})_{j\geq 0}$ are the  KAM resonances, and $n_{i_0}\leq C(\lambda,\alpha,\e)$ comes from the global to local reduction, which are uniformly for all $E\in \Sigma_{\lambda,\alpha}$. Meanwhile, let $\eta_{j+k}\in [0,\infty]$ be such that
\begin{equation}\label{gro4}
|\sin2\pi(2\rho(E)-(n_{i_0}+\cdots+n_{i_{j+k}})\alpha)|=e^{-\eta_{j+k}|n_{i_0}+\cdots+n_{i_{j+k}}|}.
\end{equation}
Once we have these,  we can state the following quantitative almost reducibility result for almost Mathieu cocycle:

\begin{Proposition}\label{Amoneed}
Let $\alpha\in DC$, $0<|\lambda|<1$ and  $E\in \Sigma_{\lambda,\alpha}$. Then for any $\e>0$, there exist sequence of $B_j\in C_{\frac{h_{\lambda}}{2\pi}(1-\e)}^\omega(\T, \PSL(2,\R))$ with $\deg{B_j}=\ell_j$, such that for any $0\leq k\leq m(j)-1$,
\begin{equation}\label{gro2}
B_{j+k}^{-1}(\theta+\alpha)S_E^{\lambda}(\theta)B_{j+k}(\theta)=A_{j+k}e^{f_{j+k}(\theta)}=M^{-1}exp \left(
\begin{array}{ccc}
 i t^{j+k} &  \nu^{j+k}\\
\bar{\nu}^{j+k} &  -i t^{j+k}
 \end{array}\right)Me^{f_{j+k}(\theta)}
\end{equation}
with estimates
\begin{enumerate}
\item $ \|f_{j+k}\|_{\frac{h_{\lambda}}{2\pi}(1-\e)} \leq  e^{- \frac{1}{32}\e h_{\lambda} |n_{i_{j+k+1}}|}$,
\item $|\nu_{j+k}| \leq e^{-h_{\lambda}(1-\e)|n_{i_{j+k}}|}$,
\item $|\nu_{j+k}|\geq e^{-h_{\lambda}(1+\e)|n_{i_{j+k}}|}$ if $(1+\e)h_\lambda<\eta_{j+k}$,
\item $|\deg B_{j+k} - n_{i_{j+k}}  |< (\e h_{\lambda})^2 |n_{i_{j+k}}|$,
\item $\|B_{j+k}\|_0 \leq  C|n_{i_{j+k}}|^{\tau}$.
\end{enumerate}
\end{Proposition}
\begin{pf}
By Avila's global theory \cite{avila2015global}, if $|\lambda|<1$, $(\alpha,S_E^\lambda)$ is subcritical in the strip $\{|\Im z|< \frac{h_{\lambda}}{2\pi}\}$, and  for any $\e>0$, $(\alpha,S_E^{\lambda})$ is almost reducible in the strip $\{|\Im z|<\frac{h_{\lambda}}{2\pi}(1- \frac{\e}{2})\}$ \cite{avila2, avila}. Then we only need to prove $(3)$, while the others follow  directly from Proposition \ref{Amo}.  

To prove this, note by the assumption 
$$
|2\rho(\alpha, A_{j+k}e^{f_{j+k}})|=  e^{-\eta_{j+k}|n_{i_0}+\cdots+n_{i_{j+k}}|} \leq e^{-h_{\lambda}(1+ \e)|n_{i_{j+k}}|},
$$
by Lemma \ref{rol}, we have
\begin{equation}\label{last1}
|\xi_{j+k}| =|\rho(\alpha, A_{j+k})| \leq  |\rho(\alpha, A_{j+k}e^{f_{j+k}})|+  e^{- \frac{1}{64}\e h_{\lambda} |n_{i_{j+k+1}}|} < 2 e^{-h_{\lambda}(1+ \e)|n_{i_{j+k}}|}.
\end{equation}
Thus if we assume $|\nu_{j+k}| \leq  e^{-h_{\lambda}(1+ \e)|n_{i_{j+k}}|}$, then \eqref{last1} implies that $|t_{j+k}|  \leq  2 e^{-h_{\lambda}(1+ \e)|n_{i_{j+k}}|} $, consequently we can rewrite
$ A_{j+k}e^{f_{j+k}(\theta)} := \id + \tilde{f}_{j+k}(\theta)$ with estimate
 $$ \| \tilde{f}_{j+k}\|_{\frac{h_{\lambda}}{2\pi}(1-\e)}  \leq e^{-h_{\lambda}(1+ \e)|n_{i_{j+k}}|}.$$
 Meanwhile, Proposition \ref{Amo} shows that  $B_{j+k}(\theta)$ is $( C |n_{i_{j+k}}|^{2\tau},   Ce^{(\e h )^2|n_{i_{j+k}}|},  h_{\lambda}(1-\e),n_{i_{j+k}})- \text{good}$.
 This is a  contradiction by Proposition \ref{lowred}. Then the result follows.

\end{pf}

\noindent{\bf Proof of Theorem \ref{amoasy}:}
 Recall that $\{\ell_j\}_{j=1}^\infty$ is a sequence such that
\begin{equation}\label{aaaa}
\|2\rho(E)-\ell_j\alpha\|_{\R/\Z}\leq e^{-\e_0|\ell_j|}.
\end{equation}
where $\e_0< h_{\lambda}.$ Thus one can take $\e$ sufficiently small such that 
\begin{equation}\label{difference}
h_{\lambda}(1-2\e) > \e_0 (1+\e h_{\lambda}),
\end{equation}
and then choose $j$ sufficiently large  such that 
\begin{eqnarray}
\label{gro11} |(\e h_{\lambda})^2 |n_{i_j}|| &\leq& \frac{1}{512}\e h_\lambda |\ell_j|,\\
\label{gro51} \frac{8\tau \ln t}{t} &\leq& \frac{\e}{4}, \quad \forall t\geq |n_{i_j}|. 
\end{eqnarray}
As direct consequence of  (3) of Proposition \ref{Amoneed} and \eqref{gro11}, we have the following consequence of KAM resonances and rotation number resonances:
\begin{equation}\label{kam-rot} (1-  \e h_{\lambda}) |\ell_j| <  |n_{i_{j}}  |< (1+ \e h_{\lambda}) |\ell_j|,
\end{equation}

Within these parameters, we can finish the whole proof.  First by  Proposition \ref{Amoneed} and Lemma \ref{growth} \footnote{Actually we need an additional global to local reduction transformation by ARC \cite{avila2, avila}, but this will not change the growth of cocycles.}, there is $C(\lambda,\alpha,\e)>0$ such that for any $0\leq k\leq m(j)-1$, we have: \\
{\rm (a)}    If $|\frac{\nu_{j+k}}{\xi_{j+k}}|<\frac{1}{2}$, then for any $e^{\frac{1}{512}\e h_\lambda|n_{i_{j+k}}|}\leq  n<  e^{ \frac{1}{32}\e h_\lambda |n_{i_{j+k+1}}|}$, we have 
$$ C^{-1} (|\ln n|)^{-4\tau}      \leq \|\mathcal{A}_{E}^n\|_0\leq C (|\ln n|)^{4\tau}.$$
{\rm (b)}   $|\frac{\nu_{j+k}}{\xi_{j+k}}| \geq \frac{1}{2}$, then we have the following:
\begin{itemize}
\item If $e^{ \frac{1}{512}\e h_\lambda |n_{i_{j+k}}|}\leq n< \min\{\frac{1}{\xi_{j+k}}, e^{\frac{1}{128}\e h_\lambda |n_{i_{j+k+1}}|}\}$, then 
$$
C^{-1} (|n||\nu_{j+k}|+1)(|\ln n|)^{-4\tau}\leq \|\mathcal{A}_{E}^n\|_0\leq C (|n||\nu_{j+k}|+1) (|\ln n|)^{4\tau}.
$$
\item If $\frac{1}{\xi_{j+k}}  \leq n  < e^{ \frac{1}{128}\e h_\lambda |n_{i_{j+k+1}}|},$ then 
$$ C^{-1} \sqrt{1+ 2\sin^2(n\xi_{j+k}) |\frac{\nu_{j+k}}{\xi_{j+k}}| } (|\ln n|)^{-4\tau}  \leq   \|\mathcal{A}_{E}^n\|_0 \leq  C \sqrt{1+ 2\sin^2(n\xi_{j+k}) |\frac{\nu_{j+k}}{\xi_{j+k}}| } (|\ln n|)^{4\tau}. $$
\end{itemize}

By \eqref{gro2}, Lemma \ref{rol} and (1) of Proposition \ref{Amoneed}, we have
\begin{equation}\label{gro3}
|2\rho(E)-2\xi_{j+k}-(n_{i_{0}}+\cdots+n_{i_{j+k}})\alpha|\leq e^{- \frac{1}{64}\e h_{\lambda} |n_{i_{j+k+1}}|}
\end{equation}
by   (4) of Proposition \ref{Amoneed} and  \eqref{gro4}, we have
\begin{equation}\label{gro31}
\frac{1}{2}e^{-\eta_{j+k}|n_{i_{0}}+\cdots+n_{i_{j+k}}|}\leq |\xi_{j+k}|\leq 2 e^{-\eta_{j+k}|n_{i_{0}}+\cdots+n_{i_{j+k}}|}.
\end{equation}
Then  we distinguish the proof into three cases:

{\bf Case I:} $1\leq k\leq m(j)-1$.  Notice that by the definition of $\ell_i$, in this case, we have $\eta_{j+k}<\e_0$, then by (2) of Proposition \ref{Amoneed}, \eqref{difference} and \eqref{gro31}, we have
$$
\left|\frac{\nu_{j+k}}{\xi_{j+k}}\right|\leq e^{-\e h_\lambda |n_{i_{j+k}}|}<\frac{1}{2}.
$$
By (a) above, for  $e^{\frac{1}{512}\e h_{\lambda}|n_{i_{j+k}}|} <|n|<e^{\frac{1}{128}\e h_{\lambda}|n_{i_{j+k+1}}|} $, we have
$$
C^{-1}(|\ln n|)^{-4\tau}\leq \|\mathcal{A}_{E}^n\|_0\leq C(|\ln n|)^{4\tau}.
$$

{\bf Case II:} $k=0$ and $\eta_{j}\leq (1+\e) h_\lambda$, then by (2) of Proposition \ref{Amoneed}, \eqref{kam-rot} and  \eqref{gro31}, 
we have 
$$
\left|\frac{\nu_{j}}{\xi_{j}}\right|\leq 2e^{-(h_\lambda-\eta_{j}-2\e h_\lambda) |\ell_j|}\leq 2e^{2\e h_\lambda |\ell_j|}.
$$
By (a), (b) above and Remark \ref{grow-rem}, for  $e^{\frac{1}{256}\e h_{\lambda}|\ell_j|} \leq |n|<e^{\frac{1}{128}\e h_{\lambda}|n_{i_{j+1}}|} ,$ we have estimate
$$
C^{-1} e^{-2\e h_\lambda |\ell_j|}(|\ln n|)^{-4\tau}\leq \|\mathcal{A}_{E}^n\|_0\leq Ce^{2\e h_\lambda |\ell_j|} (|\ln n|)^{4\tau}.
$$

{\bf Case III:} $k=0$ and $\eta_j> (1+\e)  h_\lambda$, then by  (2), (3) of Proposition \ref{Amoneed}, \eqref{kam-rot} and  \eqref{gro31},  we have 
$$
\frac{1}{2}e^{-(h_\lambda-\eta_j+2\e h_\lambda) |\ell_j|}
\leq \left|\frac{\nu_j}{\xi_j}\right|\leq 2e^{-(h_\lambda-\eta_j-2\e h_\lambda) |\ell_j|}.
$$
{\bf Case III-1: $e^{\frac{1}{256}\e h_{\lambda}|\ell_j|} \leq |n|<e^{\eta_j|\ell_{j}|} $.} 
By (a),  (b) above and \eqref{kam-rot}, for  $$e^{\frac{1}{512}\e h_{\lambda}|n_{i_j}|}<e^{\frac{1}{256}\e h_{\lambda}|\ell_j|} \leq|n|<e^{\eta_j|\ell_{j}|} ,$$ we have estimate
$$
C^{-1} (|n|e^{-h_\lambda(1+\e)|\ell_j|}+1)  (|\ln n|)^{-4\tau} \leq \|\mathcal{A}_{E}^n\|_0\leq C (|n|e^{-h_\lambda(1-\e)|\ell_j|}+1) (|\ln n|)^{4\tau}.
$$
{\bf Case III-2:   $e^{\eta_j|\ell_{j}|}\leq |n|< e^{\frac{1}{128}\e h_{\lambda}|n_{i_{j+1}}|} $.} In this case,   \eqref{gro3}  implies that 
$$
\|2n\xi_j-n(2\rho(E)-(n_{i_{0}}+\cdots+n_{i_{j+k}})\alpha)\|_{\R/\Z}\leq e^{- \frac{1}{128}\e h_{\lambda} |n_{i_{j+k+1}}|}.
$$
Consequently by $(b)$ above, we have
\begin{align*}
&C^{-1} \sqrt{1+ 2\sin^22\pi(n(\rho(E)-\ell_j\alpha/2)) e^{-2(h_\lambda-\eta_j+2\e h_\lambda) |\ell_j|} }  (|\ln n|)^{-4\tau}\leq   \|\mathcal{A}_{E}^n\|_0 \\
&\leq  C\sqrt{1+ 2\sin^22\pi(n(\rho(E)-\ell_j\alpha/2)) e^{-2(h_\lambda-\eta_j-2\e h_\lambda) |\ell_j|} }  (|\ln n|)^{4\tau}. 
\end{align*}

Combining these cases together, by \eqref{gro51}, these imply that  
\begin{equation}\label{gr09}
|n|^{ f(n)-\e}\leq \|\mathcal{A}_{E}^n\|_0 \leq |n|^{f(n)+ \e},
\end{equation}
where we denote 
\begin{equation*}
f(n)=\begin{cases}\max\{1-\frac{|{\ell_j}|h_\lambda}{\ln|n|},0\}& e^{\frac{1}{256}\e h_{\lambda}|\ell_i|} \leq |n|<e^{\eta_{j}|\ell_j|}\\
\max\{\frac{\ln |\sin 2\pi n(\rho(E)+\ell_j\alpha/2)|+|\ell_j|(\eta_{j}-h_\lambda)}{\ln|n|},0\}& e^{\eta_{j}|\ell_j|}<|n|<e^{\frac{1}{256}\e h_{\lambda} |\ell_{j+1}|}.
\end{cases}
\end{equation*}
As a consequence, \eqref{op-gr1} and \eqref{op-gr2} follow directly from  \eqref{gr09}. 
\qed

\bigskip
{\bf Proof of Corollary \ref{amo reducibility}:} By Avila's global theory \cite{avila2015global}, if $|\lambda|<1$, $(\alpha,S_E^\lambda)$ is subcritical in the strip $\{|\Im z|<-\frac{\ln \lambda}{2\pi}\}$, thus by Theorem \ref{reducibility main2}, if $\alpha\in DC$ and $-\ln \lambda >\delta(\alpha,\rho)$, $(\alpha,S_E^\lambda)$ is reducible. On the other hand, $-\ln \lambda <\delta(\alpha,\rho)$, by \eqref{op-gr1}, $(\alpha,S_E^\lambda)$ is unbounded, which immediately imply $(\alpha,S_E^\lambda)$ is not analytically  reducible (even not $C^0$ reducible). 
\qed

\section*{Acknowledgements} 
 This work was partially supported by National Key R\&D Program of China (2020 YFA0713300) and Nankai Zhide Foundation. L. Ge was partially supported by NSFC grant (12371185) and the Fundamental Research Funds for the Central Universities (the start-up fund), Peking University.  J. You was also partially supported by NSFC grant (11871286). Q. Zhou was supported by NSFC grant (12071232).

\bibliographystyle{spmpsci}
\bibliography{universal}

\begin{thebibliography}{10}
\providecommand{\url}[1]{{#1}}
\providecommand{\urlprefix}{URL }
\expandafter\ifx\csname urlstyle\endcsname\relax
  \providecommand{\doi}[1]{DOI~\discretionary{}{}{}#1}\else
  \providecommand{\doi}{DOI~\discretionary{}{}{}\begingroup
  \urlstyle{rm}\Url}\fi

\bibitem{amor}
Amor, S.H.: H{\"o}lder continuity of the rotation number for quasi-periodic
  co-cycles in $\mathrm{SL} (2,\mathbb{R})$.
\newblock Commun. Math. Phys. \textbf{287}(2), 565--588 (2009)

\bibitem{arnold}
Arnold, V.I.: Small denominators and problems of stability of motion in
  classical and celestial mechanics.
\newblock Russian Mathematical Surveys (1963)

\bibitem{AKN}
Arnold, V.I., Kozlov, V.V., Neishtadt, A.: Mathematical aspects of classical
  and celestial mechanics. Vol. 3.
\newblock Berlin: Springer (2006)

\bibitem{aubryandre}
Aubry, S., Andr{\'e}, G.: Analyticity breaking and {Anderson} localization in
  incommensurate lattices.
\newblock Ann. Israel Phys. Soc. \textbf{3}, 133--164 (1980)

\bibitem{avila}
Avila, A.: The absolutely continuous spectrum of the almost {Mathieu} operator.
\newblock preprint arXiv:0810.2965  (2008)

\bibitem{avila1}
Avila, A.: Almost reducibility and absolute continuity {I}.
\newblock preprint arXiv:1006.0704  (2010)

\bibitem{avila2015global}
Avila, A.: Global theory of one-frequency {Schr{\"o}dinger} operators.
\newblock Acta Math. \textbf{215}, 1--54 (2015)

\bibitem{avila2}
Avila, A.: {KAM}, {Lyapunov} exponents, and the spectral dichotomy for typical
  one-frequency {Schr\"odinger} operators.
\newblock preprint arXiv:2307.11071v2  (2023)

\bibitem{abd}
Avila, A., Bochi, J., Damanik, D.: Cantor spectrum for schr\"odinger operators
  with potentials arising from generalized skew-shifts.
\newblock Duke Math.J. \textbf{146}, 253--280 (2009)

\bibitem{afk}
Avila, A., Fayad, B., Krikorian, R.: A {KAM} scheme for {${\rm SL}(2,\Bbb R)$}
  cocycles with {L}iouvillean frequencies.
\newblock Geom. Funct. Anal. \textbf{21}(5), 1001--1019 (2011)

\bibitem{aj}
Avila, A., Jitomirskaya, S.: The ten martini problem.
\newblock Ann. of Math. \textbf{170}, 303--342 (2009)

\bibitem{aj1}
Avila, A., Jitomirskaya, S.: Almost localization and almost reducibility.
\newblock J. Eur. Math. Soc. \textbf{12}(1), 93--131 (2010)

\bibitem{AK06}
Avila, A., Krikorian, R.: Reducibility or nonuniform hyperbolicity for
  quasiperiodic {{Schr\"odinger}} cocycles.
\newblock Ann. of Math. \textbf{164}(3), 911--940 (2006)

\bibitem{ALYZ}
Avila, A., Last, Y., Shamis, M., Zhou, Q.: On the abominable properties of the
  almost mathieu operator with well approximated frequencies.
\newblock Duke Math. J. \textbf{173}(4), 603--672 (2024)

\bibitem{ayz}
Avila, A., You, J., Zhou, Q.: Sharp phase transitions for the almost {M}athieu
  operator.
\newblock Duke Math. J. \textbf{166}(14), 2697--2718 (2017)

\bibitem{ayz2}
Avila, A., You, J., Zhou, Q.: Dry ten {Martini} problem in the non-critical
  case.
\newblock preprint arXiv:2306.16254v2  (2023)

\bibitem{as}
Avron, J., Simon, B.: Singular continuous spectrum for a class of almost
  periodic {J}acobi matrices.
\newblock Bull. Amer. Math. Soc. (N.S.) \textbf{6}(1), 81--85 (1982)

\bibitem{Av17}
Avron, J.E.: Topological quantum states, the quantum hall effect, and the 2016
  nobel prize.
\newblock IAMP News Bulletin  (2017)

\bibitem{aos}
Avron, J.E., Osadchy, D., Seiler, R.: A topological look at the quantum {Hall}
  effect.
\newblock Physics today. pp. 38--42 (2003)

\bibitem{bs}
Bellissard, J., Simon, B.: Cantor spectrum for the almost mathieu equation.
\newblock J. Funct. Anal. \textbf{48}, 408--419 (1982)

\bibitem{binder2018almost}
Binder, I., Damanik, D., Goldstein, M., Lukic, M.: Almost periodicity in time
  of solutions of the {K}d{V} equation.
\newblock Duke Math. J. \textbf{167}(14), 2633--2678 (2018)

\bibitem{B02}
Bourgain, J.: On the spectrum of lattice {S}chr\"{o}dinger operators with
  deterministic potential. {II}.
\newblock J. Anal. Math. \textbf{88}, 221--254 (2002).
\newblock Dedicated to the memory of Tom Wolff

\bibitem{b1}
Bourgain, J.: Green's function estimates for lattice {S}chr\"odinger operators
  and applications, \emph{Annals of Mathematics Studies}, vol. 158.
\newblock Princeton University Press, Princeton, NJ (2005)

\bibitem{bj02}
Bourgain, J., Jitomirskaya, S.: Absolutely continuous spectrum for {1D}
  quasiperiodic operators.
\newblock Invent. Math. \textbf{148}(3), 453--463 (2002)

\bibitem{CCYZ}
Cai, A., Chavaudret, C., You, J., Zhou, Q.: Sharp {H}\"{o}lder continuity of
  the {L}yapunov exponent of finitely differentiable quasi-periodic cocycles.
\newblock Math. Z. \textbf{291}(3-4), 931--958 (2019)

\bibitem{CSZ1}
Cao, H., Shi, Y., Zhang, Z.: Localization and regularity of the integrated
  density of states for {S}chr\"odinger operators on {$\Bbb Z^d$} with
  {$C^2$}-cosine like quasi-periodic potential.
\newblock Comm. Math. Phys. \textbf{404}(1), 495--561 (2023)

\bibitem{CSZ2}
Cao, H., Shi, Y., Zhang, Z.: On the spectrum of quasi-periodic schr\" odinger
  operators on $\mathbb {Z}^ d $ with $c^2$-cosine type potentials.
\newblock preprint arxiv:2310.07407  (2023)

\bibitem{cey}
Choi, M., Elliott, G., Yui, N.: Gauss polynomials and the rotation algebra.
\newblock Invent. Math. \textbf{99}, 225--246 (1990)

\bibitem{CD}
Chulaevsky, V., Dinaburg, E.: Methods of kam-theory for long-range
  quasi-periodic oper- ators on $\mathbb{Z}^{\mu}$. pure point spectrum.
\newblock Commun. Math. Phys \textbf{153}, 559--577 (1993)

\bibitem{dg}
Damanik, D., Goldstein, M.: On the inverse spectral problem for the
  quasi-periodic {S}chr{\"o}dinger equation.
\newblock Publ. Math. Inst. Hautes {\'E}tudes Sci. \textbf{119}, 217--401
  (2014)

\bibitem{Deift1}
Deift, P.: Some open problems in random matrix theory and the theory of
  integrable systems.
\newblock In: Integrable systems and random matrices, \emph{Contemp. Math.},
  vol. 458, pp. 419--430. Amer. Math. Soc., Providence, RI (2008)

\bibitem{Deift2}
Deift, P.: Some open problems in random matrix theory and the theory of
  integrable systems. {II}.
\newblock SIGMA Symmetry Integrability Geom. Methods Appl. \textbf{13}, Paper
  No. 016, 23 (2017)

\bibitem{DS}
Dinaburg, E., Sinai., Y.: The one dimensional schr\"odinger equation with a
  quasi-periodic potential.
\newblock Funct. Anal. Appl. pp. 279--289 (1975)

\bibitem{S14}
Dumas, H.S.: The {KAM} story: A friendly introduction to the content, History,
  and significance Of classical {K}olmogorov-{A}rnold-{M}oser theory.
\newblock World Scientific Publishing Company (2014)

\bibitem{eli}
Eliasson, L.: Floquet solutions for the 1-dimensional quasi-periodic
  {Schr{\"o}dinger} equation.
\newblock Commun. Math. Phys. \textbf{146}(3), 447--482 (1992)

\bibitem{FK}
Fayad, B., Krikorian, R.: Rigidity results for quasiperiodic {${\rm SL}(2,\Bbb
  R)$}-cocycles.
\newblock J. Mod. Dyn. \textbf{3}(4), 497--510 (2009)

\bibitem{f}
Fedotov, A.: Monodromization method in the theory of almost-periodic equations.
\newblock St. Petersburg Mathematical Mathematics \textbf{179}(2), 153--196
  (1997)

\bibitem{Fe}
Fedotov, A.: A series of spectral gaps for the almost {Mathieu} operator with a
  small coupling constant.
\newblock preprint arXiv:2012.03356  (2020)

\bibitem{gj1}
Ge, L., Jitomirskaya, S.: Hidden subcriticality, symplectic structure, and
  universality of sharp arithmetic spectral results for type i operators.
\newblock preprint arXiv:2407.08866  (2024)

\bibitem{gjy}
Ge, L., Jitomirskaya, S., You, J.: Kotani theory, puig's argument, and
  stability of the ten martini problem.
\newblock preprint arXiv:2308.09321  (2023)

\bibitem{gjyz}
Ge, L., Jitomirskaya, S., You, J., Zhou, Q.: Multiplicative jensen's formula
  and quantitative global theory of one-frequency schr{\"o}dinger operators.
\newblock preprint arXiv:2306.16387  (2023)

\bibitem{gk}
Ge, L., Kachkovskiy, I.: Ballistic transport for one-dimensional quasiperiodic
  {S}chr{\"o}dinger operators.
\newblock Comm. Pure Appl. Math. \textbf{76}(10), 2577--2612 (2023)

\bibitem{gy}
Ge, L., You, J.: Arithmetic version of {A}nderson localization via
  reducibility.
\newblock Geom. Funct. Anal. \textbf{30}(5), 1370--1401 (2020)

\bibitem{gyzh}
Ge, L., You, J., Zhao, X.: H\"older regularity of the integrated density of
  states for quasi-periodic long-range operators on {$\ell^2(\Bbb Z^d)$}.
\newblock Comm. Math. Phys. \textbf{392}(2), 347--376 (2022)

\bibitem{gyz}
Ge, L., You, J., Zhou, Q.: Exponential dynamical localization: Criterion and
  applications.
\newblock Ann. Sci. \'{E}c. Norm. Sup\'{e}r. (4) \textbf{56}(1), 91--126 (2023)

\bibitem{gor}
Gordon, A.: The point spectrum of one-dimensional {Schr\"doinger} operators.
\newblock Uspekhi Mat. Nauk. \textbf{31}, 257--258 (1976)

\bibitem{gjls}
Gordon, A., Jitomirskaya, S., Last, Y., Simon, B.: Duality and singular
  continuous spectrum in the almost mathieu equation.
\newblock Acta Math. \textbf{178}, 169--183 (1997)

\bibitem{H}
Halperin, B.: Quantized hall conductance, current-carrying edge states, and the
  existence of extended states in a two-dimensional disordered potential.
\newblock Phys. Rev. B \textbf{25}, 2185 (1982)

\bibitem{harper}
Harper, P.: Single band motion of conducion electrons in a unifor magnetic
  field.
\newblock Pro. Phys. Soc. \textbf{68}(10), 874--878 (1955)

\bibitem{hs}
Helffer, B., Sj\"ostrand, J.: Semi-classical analysis for {Haper's} equation.
  iii: Cantor structure of the spectrum.
\newblock Mem. Soc. Math. France \textbf{39}, 1--124 (1989)

\bibitem{houyou}
Hou, X., You, J.: Almost reducibility and non-perturbative reducibility of
  quasi-periodic linear systems.
\newblock Invent. Math. \textbf{190}(1), 209--260 (2012)

\bibitem{BJ}
J.Bourgain, Jitomirskaya, S.: Continuity of the lyapunov exponent for
  quasiperiodic operators with analytic potential.
\newblock J. Stat. Phys. pp. 1203--1218 (2002)

\bibitem{JLZ}
Jiang, K., Li, S., Zhang, P.: Numerical methods and analysis of computing
  quasiperiodic systems.
\newblock SIAM J. Numer. Anal. \textbf{62}(1), 353--375 (2024)

\bibitem{JZ}
Jiang, K., Zhang, P.: Numerical mathematics of quasicrystals.
\newblock In: Proceedings of the {I}nternational {C}ongress of
  {M}athematicians---{R}io de {J}aneiro 2018. {V}ol. {IV}. {I}nvited lectures,
  pp. 3591--3609. World Sci. Publ., Hackensack, NJ (2018)

\bibitem{J}
Jitomirskaya, S.: Almost everything about the almost {M}athieu operator. {II}.
\newblock In: X{I}th {I}nternational {C}ongress of {M}athematical {P}hysics
  ({P}aris, 1994), pp. 373--382. Int. Press, Cambridge, MA (1995)

\bibitem{J99}
Jitomirskaya, S.: Metal-insulator transition for the almost {M}athieu operator.
\newblock Ann. of Math. (2) \textbf{150}(3), 1159--1175 (1999)

\bibitem{JICM}
Jitomirskaya, S.: One-dimensional quasiperiodic operators: global theory,
  duality, and sharp analysis of small denominators.
\newblock In: I{CM}---{I}nternational {C}ongress of {M}athematicians. {V}ol. 2.
  {P}lenary lectures, pp. 1090--1120. EMS Press, Berlin (2023)

\bibitem{jk2}
Jitomirskaya, S., Kachkovskiy, I.: $l^2$-reducibility and localization for
  quasiperiodic operators.
\newblock Math. Res. Lett. \textbf{23}, 431--444 (2016)

\bibitem{JLiu}
Jitomirskaya, S., Liu, W.: Universal hierarchical structure of quasiperiodic
  eigenfunctions.
\newblock Ann. of Math. \textbf{187}(3), 721--776 (2018)

\bibitem{JLiu2}
Jitomirskaya, S., Liu, W.: Universal reflective-hierarchical structure of
  quasiperiodic eigenfunctions and sharp spectral transition in phase.
\newblock To appear in J. Eur. Math. Soc.  (2018)

\bibitem{JS94}
Jitomirskaya, S., Simon, B.: Operators with singular continuous spectrum.
  {III}. {A}lmost periodic {S}chr\"odinger operators.
\newblock Comm. Math. Phys. \textbf{165}(1), 201--205 (1994)

\bibitem{JM}
Johnson., R., Moser, J.: The rotation number for almost periodic potentials.
\newblock Comm. Math. Phys. \textbf{84}(3), 403--438 (1982)

\bibitem{Kac}
Kac, M.: Public commun.
\newblock AMS Annual Meeting  (1981)

\bibitem{TP}
Kappeler, T., P\"oschel, J.: Kd{V} \& {KAM}.
\newblock Springer-Verlag, Berlin (2013)

\bibitem{KXZ}
Karaliolios, N., Xu, X., Zhou, Q.: Anosov-{K}atok constructions for
  quasi-periodic {${\rm SL}(2,{\Bbb R})$} cocycles.
\newblock Peking Math. J. \textbf{7}(1), 203--245 (2024)

\bibitem{KDP}
Klitzing K.~V., D.G..P.N.: New method for high-accuracy determination of the
  fine-structure constant based on quantized hall resistance.
\newblock Phys. Rev. Lett. \textbf{45}, 494--497 (1980)

\bibitem{KJYZ}
Krikorian, R., Wang, J., You, J., Zhou, Q.: Linearization of quasiperiodically
  forced circle flows beyond brjuno condition.
\newblock Communications in Mathematical Physics \textbf{358}, 81--100 (2018)

\bibitem{GYZh2}
L.~Ge, J.Y., Zhao, X.: The arithmetic version of the frequency transition
  conjecture: new proof and generalization.
\newblock Peking. Math. J. \textbf{5}(2), 349--364 (2021)

\bibitem{last}
Last, Y.: Zero measure spectrum for the almost mathieu operator.
\newblock Commun. Math. Phys. \textbf{164}, 421--432 (1994)

\bibitem{L1}
Last, Y.: Spectral theory of {S}turm-{L}iouville operators on infinite
  intervals: a review of recent developments.
\newblock In: Sturm-{L}iouville theory, pp. 99--120. Birkh\"{a}user, Basel
  (2005)

\bibitem{LYZZ}
Leguil, M., You, J., Zhao, Z., Zhou, Q.: {Asymptotics of spectral gaps of
  quasi-periodic Schr\"odinger operators}.
\newblock Cambridge Journal of Mathematics,
  https://dx.doi.org/10.4310/CJM.241128001445  (2017)

\bibitem{LYZ}
Li, X., You, J., Zhou, Q.: Exact local distribution of the absolutely
  continuous spectral measure.
\newblock https://arxiv.org/abs/2407.09278  (2024)

\bibitem{ly}
Liu, W., Yuan, X.: Anderson localization for the almost {M}athieu operator in
  the exponential regime.
\newblock J. Spectr. Theory \textbf{5}(1), 89--112 (2015)

\bibitem{gs}
M.~Goldstein, W.S.: On resonances and the formation of gaps in the spectrum of
  quasi-periodic schr\"odinger equations.
\newblock Ann. of Math. \textbf{173}, 337--475 (2011)

\bibitem{He79}
M.Herman: Sur la conjugaison differentiable des diffeomorphismes du cercle a
  des rotations.
\newblock Inst. Hautes Etudes Sci. Publ. Math (5-233) (1979)

\bibitem{mp}
Moser, J., P\"{o}schel, J.: An extension of a result by {Dinaburg and S}inai on
  quasiperiodic potentials.
\newblock Comment. Math. Helv. \textbf{59}(1), 39--85 (1984)

\bibitem{ntw}
Niu, Q., Thouless, D., Wu, Y.: Quantized {Hall} conductance as a topological
  invariant.
\newblock Phys. Rev. B \textbf{31}, 3372 (1985)

\bibitem{oa}
Osadchy, D., Avron, J.E.: Hofstadter butterfly as quantum phase diagram.
\newblock J. Math Phys. \textbf{42}, 5665--5671 (2001)

\bibitem{Peierls}
Peierls, R.: Zur theorie des diamagnetismus von leitungselektronen.
\newblock Z. Phy. \textbf{80}(11), 763--791 (1933)

\bibitem{PuigC}
Puig, J.: Cantor spectrum for the almost mathieu operator.
\newblock Comm. Math. Phys. \textbf{244}, 297--309 (2004)

\bibitem{puig}
Puig, J.: A nonperturbative eliasson's reducibility theorem.
\newblock Nonlinearity \textbf{19}(2), 355--376 (2005)

\bibitem{R}
Rauh, A.: Degeneracy of landau levels in crystals.
\newblock Phys. Status Solidi B \textbf{65}, 131--135 (1974)

\bibitem{Barry}
Simon, B.: Almost periodic {S}chr{\"o}dinger operators: a review.
\newblock Advances in Applied Mathematics \textbf{3}(4), 463--490 (1982)

\bibitem{simon}
Simon, B.: Localization in general one-dimensional random systems. {I}.
  {J}acobi matrices.
\newblock Comm. Math. Phys. \textbf{102}(2), 327--336 (1985)

\bibitem{Sinai}
Sinai, Y.G.: Anderson localization for one-dimensional difference
  {Schr{\"o}dinger} operator with quasiperiodic potential.
\newblock J. Stat. Phys. \textbf{46}(5-6), 861--909 (1987)

\bibitem{Yo84}
Yoccoz, J.C.: Conjugaison differentiable des diffeomorphismes du cercle dont le
  nombre de rotation verifie une condition diophantienne.
\newblock Ann. Sci. Ecole Norm. Sup. (333-359) (1984)

\bibitem{YouICM}
You, J.: Quantitative almost reducibility and its applications.
\newblock In: Proceedings of the {I}nternational {C}ongress of
  {M}athematicians---{R}io de {J}aneiro 2018. {V}ol. {III}. {I}nvited lectures,
  pp. 2113--2135. World Sci. Publ., Hackensack, NJ (2018)

\bibitem{yzhou}
You, J., Zhou, Q.: Embedding of analytic quasi-periodic cocycles into analytic
  quasi-periodic linear systems and its applications.
\newblock Comm. Math. Phys. \textbf{323}(3), 975--1005 (2013)

\bibitem{ZhouW12}
Zhou, Q., Wang, J.: Reducibility results for quasiperiodic cocycles with
  {L}iouvillean frequency.
\newblock J. Dynam. Differential Equations \textbf{24}(1), 61--83 (2012)

\end{thebibliography}

\end{document}